\documentclass[]{interact}

\usepackage{epstopdf}
\usepackage[caption=false]{subfig}
\usepackage{microtype}
\usepackage{hyperref}

\usepackage[numbers,sort&compress]{natbib}
\bibpunct[, ]{[}{]}{,}{n}{,}{,}

\usepackage{amsmath}
\usepackage{amsthm}
\usepackage{amssymb}
\usepackage{graphicx}
\usepackage{color}
\usepackage{ifpdf}
\usepackage{url}
\usepackage{paralist}
\usepackage{enumitem}

\usepackage{algorithm}
\usepackage{algorithmic}
\usepackage[T1]{fontenc} 

\usepackage{multicol} 
\usepackage{booktabs} 
\usepackage{float} 
\renewcommand{\thesection}{ \arabic{section}}
\renewcommand{\thesubsection}{\arabic{section}.\arabic{subsection}}
\renewcommand{\thesubsubsection}{\arabic{section}.\arabic{subsection}.\alph{subsubsection}}
 \theoremstyle{plain}
 \newtheorem{theorem}{Theorem}[section]
 \newtheorem{corollary}{Corollary}[section]
 \newtheorem{lemma}{Lemma}[section]
 \newtheorem{proposition}{Proposition}[section]
 \theoremstyle{definition}
 \newtheorem{definition}{Definition}[section]
 \newtheorem{assumption}{Assumption}[section]
 \theoremstyle{remark}
 \newtheorem{remark}{Remark}[section]

\newcommand{\Id}{\mathbb{I}}

\newcommand{\R}{\mathbb{R}}

\newcommand{\set}[1]{\left\{#1\right\}}
\newcommand{\sets}[1]{\{#1\}}

\newcommand{\norms}[1]{\Vert#1\Vert}

\newcommand{\Eproof}{\hfill $\square$}

\newcommand{\dom}[1]{\mathrm{dom}(#1)}

\newcommand{\zero}[1]{{\boldsymbol{0}}}

\newcommand{\Expsk}[2]{\mathbb{E}_{#1}\big[#2\big]}
\newcommand{\Expk}[1]{\mathbb{E}\big[#1\big]}

\newcommand{\Prob}[1]{\mathbf{Prob}\left(#1\right)}

\newcommand{\zer}[1]{\mathrm{zer}(#1)}
\newcommand{\zerM}[1]{\mathrm{zer}_M(#1)}

\newcommand{\gra}[1]{\mathrm{gra}(#1)}

\newcommand{\xopt}{x^{\star}}
\newcommand{\uopt}{u^{\star}}

\newcommand{\Hc}{\mathcal{H}}

\newcommand{\Bc}{\mathcal{B}}

\newcommand{\Sc}{\mathcal{S}}

\newcommand{\Dc}{\mathcal{D}}
\newcommand{\Lc}{\mathcal{L}}
\newcommand{\Qc}{\mathcal{Q}}
\newcommand{\Gc}{\mathcal{G}}

\newcommand{\Tc}{\mathcal{T}}

\newcommand{\Fc}{\mathcal{F}}

\newcommand{\Vc}{\mathcal{V}}
\newcommand{\Cc}{\mathcal{C}}
\newcommand{\Ec}{\mathcal{E}}
\newcommand{\Nc}{\mathcal{N}}
\newcommand{\Pc}{\mathcal{P}}

\newcommand{\Rc}{\mathcal{R}}

\newcommand{\iprods}[1]{\langle #1\rangle}
\newcommand{\Exp}[1]{\mathbb{E}\left[#1\right]}

\newcommand{\BigO}[1]{\mathcal{O}\left(#1\right)}
\newcommand{\BigOs}[1]{\mathcal{O}\big(#1\big)}
\newcommand{\BigOlog}[1]{\widetilde{\mathcal{O}}\big(#1\big)}
\newcommand{\SmallOs}[1]{o\big(#1\big)}
\newcommand{\SmallO}[1]{o\left(#1\right)}

\newcommand{\p}{\mathbf{p}}
\newcommand{\q}{\mathbf{q}}

\newcommand{\mbf}[1]{\mathbf{#1}}

\newcommand{\beforesec}{\vspace{-3ex}}
\newcommand{\aftersec}{\vspace{-2ex}}
\newcommand{\beforesubsec}{\vspace{-3.5ex}}
\newcommand{\aftersubsec}{\vspace{-1.25ex}}

\newcommand{\myeq}[2]{
\vspace{-0.5ex}
\begin{equation}\label{#1}
{#2}
\vspace{-1ex}
\end{equation}
}
\newcommand{\myeqn}[1]{
\vspace{-0.5ex}
\begin{equation*}
{#1}
\vspace{-1ex}
\end{equation*}
}

\renewcommand{\thesection}{\arabic{section}}
\allowdisplaybreaks

\begin{document}
	
	
	\title{Optimistic Gradient Algorithms with Randomized Coordinate and Variance Reduction Strategies for Composite Inclusions}
	
	\author{
		\name{Nghia Nguyen-Trung\textsuperscript{a} and Ion Necoara\textsuperscript{b} and Quoc Tran-Dinh\textsuperscript{a}\thanks{Contact: Quoc Tran-Dinh. Email: quoctd@email.unc.edu}}
		\affil{\textsuperscript{a}Department of Statistics and Operations Research, University of North Carolina at Chapel Hill, Chapel Hill, NC 27599-3260, USA; \textsuperscript{b}Automatic Control and Systems Engineering Department, Politehnica Bucharest, 060042 Bucharest, Romania}
	}
	
	\maketitle
	
	\begin{abstract}
		In this paper, we consider a class of possibly nonmonotone finite-sum composite inclusions in which evaluating all components or updating all coordinates at every iteration can be computationally expensive.
		We develop two optimistic gradient frameworks that combine variance reduction with randomized block-coordinate updates.
		The first framework applies under a blockwise average Lipschitz continuity condition and the existence of a weak-Minty solution.
		It achieves an $\BigOs{1/k}$ convergence rate for the expected squared residual, and its iterates converge almost surely under mild additional conditions.
		The second framework relies on a backward-forward reformulation under average co-coercivity and co-hypomonotonicity.
		It achieves faster $\BigOs{1/k^2}$ and $\SmallOs{1/k^2}$ rates for the expected squared residual, as well as $\SmallOs{1/k^2}$ almost-sure rates.
		For both frameworks, we provide a generic condition on the stochastic estimator that covers mini-batch, loopless SVRG, and a generalized JacSketch-type estimator, and we derive the corresponding oracle complexities under importance-based block sampling.
		In particular, in highly heterogeneous regimes, the L-SVRG and JacSketch variants attain $\BigOs{n + n^{2/3}\epsilon^{-2}}$ and $\BigOs{n + n^{2/3}\epsilon^{-1}}$ oracle complexities for the nonaccelerated and accelerated methods, respectively.
		Numerical experiments illustrate the effects of variance reduction and block-coordinate sampling, and compare the proposed methods with existing baselines.
	\end{abstract}
	
	
	\beforesec
	\section{Introduction}\label{sec:ch2:intro}
	\aftersec
	The goal of this paper is to develop two novel stochastic algorithms to approximate solutions of the following large finite-sum \textit{composite inclusion} (also called \textit{generalized equation}) in a high-dimensional space $\R^p$:
	\myeq{eq:FNI}{
		\begin{array}{c}
			\textrm{Find $x^\star \in \R^p$ such that:} ~ 0 \in \Psi \xopt := Gx^\star + Tx^\star, \quad \text{where} \quad Gx := \frac{1}{n}\sum_{i=1}^n G_i x.
		\end{array}
		\tag{CI}
	}
	Here, $G: \R^p \to \R^p$ and $G_i: \R^p \to \R^p$ are single-valued operators, $i \in [n] := \sets{1, \cdots, n}$, and  $T: \R^p \rightrightarrows \R^p$ is a possibly multivalued operator.
	
	Problem~\eqref{eq:FNI} provides a unified template to cover a broad range of applications, including the optimality conditions of constrained and unconstrained optimization, the optimality conditions of minimax and saddle point problems,  systems of equations and inequalities, fixed-point problems, and variational inequalities; see, e.g., \cite{Bauschke2011,Facchinei2003,ryu2016primer,tran2024revisiting}.
It also finds diverse applications in optimization, game theory, mechanics, economics, and machine learning, see, e.g., \cite{banach1922sur,von1928theorie,hartman1966some,rockafellar1976monotone,sinha2018certifying,goodfellow2014generative,gidel2019variational,madry2018towards,li2019robust,deng2020distributionally}.
	
	The setting we study in this paper has three distinctive features as follows.
\begin{compactitem}
	\item 
	First, the operator $\Psi = G+T$ may be nonmonotone, as occurs, for example, in nonconvex-nonconcave minimax problems  \cite{lin2020gradient, diakonikolas2021efficient} and statistical learning with nonconvex regularizers \cite{tran2024accelerated}, so the classical convergence theory for monotone inclusions is not directly applicable.
	
	\item 
	Second, the finite-sum structure may involve a large number of components ($n\gg1$), making a full evaluation of $G$ expensive.
	
	\item
	Third, when the ambient dimension is large ($p\gg1$), updating all coordinates at every iteration can itself be computationally prohibitive.
\end{compactitem}	
	These features motivate methods that exploit both the finite-sum and block-coordinate structures while retaining convergence guarantees beyond the monotone setting.

	\beforesubsec
	\subsection{Our solution and main contributions}\label{subsec:ch2:contribution}
	\aftersubsec
	Variance reduction and randomized block-coordinate methods have each been studied extensively in the literature,  see Subsection~\ref{subsec:ch2:related_work} for more details, but combining them within optimistic gradient and extragradient-type methods raises several issues that do not appear in their full-dimensional or non-variance-reduced counterparts.
	To overcome these issues, we rely on the following two techniques.
	First, variance reduction is applied to the optimistic gradient (OG) direction (see Subsection~\ref{subsubsec:ch2:vr_estimators} for its precise definition), which couples two consecutive iterates and allows a varying blockwise parameter $\gamma^k$, rather than directly to a single finite-sum operator $Gx^k$.
	Second, the backward step involving the resolvent of $T$ must preserve the computational advantage of block-coordinate updates.
	In the accelerated framework, two further issues arise: the backward-forward reformulation maintains the finite-sum structure but changes the component operators, while the inertial term is intrinsically full-dimensional.
	We handle these issues through a common-shift estimator construction and a block-coordinate implementation in which only the sampled block needs to be modified at each iteration.
	
	\vspace{0.75ex}
	\noindent\textbf{$\mathrm{(a)}$~The non-accelerated framework.}
	We first develop a novel \underline{v}ariance-reduced randomized block-\underline{c}oordinate \underline{o}ptimistic \underline{g}radient method (VCOG) for possibly nonmonotone \eqref{eq:FNI}, using a general class of unbiased estimators of the optimistic direction that includes mini-batch, L-SVRG, and generalized JacSketch estimators.
	Under block-wise average Lipschitz continuity, the existence of a weak-Minty solution, and a Lipschitz regularity condition on the resolvent of $T$, we establish an $\BigOs{1/k}$ rate for the averaged expected squared residual, together with almost-sure residual and iterate convergence under mild additional conditions.
	When using importance-based block sampling, the L-SVRG and generalized JacSketch variants attain a component-oracle complexity of $\BigOs{n+n^{2/3}\epsilon^{-2}}$ in the heterogeneous-block regime.
	
	\vspace{0.75ex}
	\noindent\textbf{$\mathrm{(b)}$~The accelerated framework.}
	Under the average co-coercivity of $G$ and block-separable co-hypomonotonicity of $T$, which also allows nonmonotonicity, we develop a new \underline{v}ariance-reduced randomized block-\underline{c}oordinate \underline{f}ast \underline{OG} method (VCFOG) for \eqref{eq:FNI}, using a backward-forward reformulation.
	We establish $\BigOs{1/k^2}$ expected squared-residual bounds and $\SmallO{1/k^2}$ asymptotic rates, together with almost-sure convergence of the residual and the iterates.
	When using importance-based sampling, the L-SVRG and generalized JacSketch variants attain a component-oracle complexity of $\BigOs{n+n^{2/3}\epsilon^{-1}}$ in the heterogeneous-block regime.
	This complexity is no worse than that of the full-block update regime and matches the complexity of SVRG methods in the literature.
	
	\beforesubsec
	\subsection{Related work}\label{subsec:ch2:related_work}
	\aftersubsec
	The theory and methods for \eqref{eq:FNI} and its special cases are well-developed, particularly under monotonicity assumptions \citep{Bauschke2011,Facchinei2003}.
	In this work, we focus on nonmonotone settings and exploit recent advances in variance-reduction and randomized coordinate schemes.
	
	\vspace{0.5ex}
	\noindent$\mathrm{(a)}$~\textit{Variance-reduction methods.}
	Variance-reduction techniques such as SAGA \cite{Defazio2014}, SVRG \cite{SVRG}, SARAH \cite{nguyen2017sarah}, and related control-variate schemes have been widely developed for finite-sum optimization and subsequently extended to variational inequalities, saddle-point problems, and nonlinear equations; see, e.g., \cite{davis2022variance,alacaoglu2022beyond,cai2023variance,bot2019forward,carmon2019variance,chavdarova2019reducing,huang2022accelerated,palaniappan2016stochastic,yu2022fast,cui2021analysis,tran2024accelerated} and the references therein.
	Of particular relevance, \cite{TranDinh2024} develops variance-reduced splitting schemes for possibly nonmonotone generalized equations by constructing unbiased estimators of an optimistic search direction.
	The methods there remain full-dimensional and use a fixed control parameter.
	Accelerated variance-reduced OG schemes for nonmonotone generalized equations are further studied in \cite{tran2025vfog}.
	Our estimator framework follows this direction but incorporates randomized block-coordinate updates and allows a varying blockwise parameter $\gamma^k$ to accommodate Nesterov-type acceleration.
	
	\vspace{0.5ex}
	\noindent$\mathrm{(b)}$~\textit{Block-coordinate methods.}
	Randomized block-coordinate methods have a long history in large-scale optimization
	\cite{BeckTetruashvili2013,Nesterov2012,richtarik2016parallel,wright2015coordinate}
	and have been extended to fixed-point, primal-dual, variational-inequality, and saddle-point problems
	\cite{DangLan2015,chakrabarti2024block,combettes2015stochastic,PesquetRepetti2015,peng2016arock}.
	The closest optimistic-gradient counterpart is \cite{tran2022accelerated}, which develops nonaccelerated and accelerated randomized block-coordinate OG methods for root-finding problems.
	Our focus in this paper is different: we explicitly exploit the finite-sum structure through variance-reduced estimators of the OG direction while preserving block-coordinate computation of the composite updates.
		
	\vspace{0.5ex}
	\noindent$\mathrm{(c)}$~\textit{Variance reduction with coordinate updates.}
	Methods combining variance reduction with coordinate strategies have also been studied for finite-sum composite optimization
	\cite{cai2023cyclic,hanzely2020variance}.
	More recently, \cite{diakonikolas2025block} develops a randomized block-coordinate variance-reduced method for generalized variational inequalities of Minty type and demonstrates the benefit of nonuniform block sampling.
	Our work is complementary to this line: we study OG schemes for composite inclusions under two possibly nonmonotone regimes, namely a weak-Minty setting for VCOG and an average-co-coercive/co-hypomonotone setting for VCFOG, together with almost-sure convergence and importance-based component-oracle guarantees.

	\vspace{0.75ex}
	\noindent\textbf{Paper organization.}
	Section~\ref{sec:ch2:OG_operator} collects the notation, block structure, sampling rule, and stochastic estimator framework used throughout the paper.
	Sections~\ref{sec:ch2:VCOG_alg} and \ref{sec:ch2:VACOG} develop the
	nonaccelerated VCOG and accelerated VCFOG frameworks, respectively, together with their convergence and oracle-complexity guarantees.
	Finally, Section~\ref{sec:ch2:num_experiments} presents representative numerical examples for validation and compares our results against existing methods.
	
	\beforesec
	\section{Preliminaries and Stochastic Optimistic Gradient Estimator} \label{sec:ch2:OG_operator}
	\aftersec
	\subsection{Preliminaries and block-coordinate structure}
	\aftersubsec
	\noindent\textbf{Notation and terminology.}
	For an operator $G$, $\dom{G}$, $\gra{G}$, and $J_{G}$ denote the domain, graph, and resolvent of $G$, respectively (see \cite{Bauschke2011}).
	Let $\sets{\Fc_k}_{k \geq 0}$ be the filtration generated by the algorithmic history, where $\Fc_k$ contains all information available just before the fresh random variables of iteration $k$ are drawn.
	In particular, the current and past iterates and the estimator state carried over from iteration $k-1$ are $\Fc_k$-measurable.
	In contrast, the snapshot $w^k$ in \eqref{eq:loopless_SVRG_snapshot} and the table $\sets{H_i^k}$ in \eqref{eq:JacSketch_table} are refreshed at iteration $k$, so they are not $\Fc_k$-measurable in general.
	We denote by $\mathbb{E}_k[ \cdot ] := \mathbb{E}[~\cdot\mid\Fc_k]$ the conditional expectation w.r.t.\ $\Fc_k$ and by $\mathbb{E}[\cdot]$ the total expectation.
	At iteration $k$, all fresh random variables used to construct the estimator are drawn independently of $\Fc_k$ and, whenever applicable, independently of each other.
	
	\vspace{0.75ex}
	\noindent\textbf{Block-coordinate structure.}
	Throughout this paper, we assume that $x \in \R^p$ can be decomposed into $m$ blocks as $x = [x_1, \cdots, x_m]$, where $x_j \in \R^{p_j}$, $j \in [m] := \sets{1, \cdots, m}$, $1 \leq p_j \leq p$, and $\sum_{j=1}^m p_j = p$.
	For the operator $G$ and its components $G_i$ in \eqref{eq:FNI} and $x \in \dom{G}$, let $[Gx]_j$ and $[G_i x]_j$ denote the $j$-th block of $Gx$ and $G_i x$, respectively, so that $Gx = [ [Gx]_1, \cdots, [Gx]_m]$ and $G_i x = [[G_i x]_1, \cdots, [G_i x]_m]$.
	We also denote
	\myeqn{
		\arraycolsep=0.2em
		\begin{array}{lclcl}
			G_{[j]}x &=& [0, \cdots, 0, [Gx]_j, 0, \cdots, 0] \quad \text{and} \quad
			G_{i, [j]} x &=& [0, \cdots, 0, [G_i x]_j, 0, \cdots, 0].
		\end{array}
		\vspace{-1ex}
	}
	For given $\alpha,\ \beta \in \R^m$, a scalar $\nu \in \R$, and block vectors $u, v \in \R^p$, we use the following block-coordinate-wise conventions: $\alpha^\nu := (\alpha_1^{\nu}, \dots, \alpha_m^{\nu})$, $\alpha \beta := (\alpha_1 \beta_1, \cdots, \alpha_m \beta_m)$, $\alpha \circ u := \left[\alpha_1 u_1, \cdots, \alpha_m u_m \right]$, $\alpha \circ \iprods{u,v} := \sum_{j = 1}^m \alpha_j \iprods{u_j, v_j}$, and $\alpha \circ \norms{u}^2 := \sum_{j=1}^m \alpha_j \norms{u_j}^2$.
	If $\beta \in \R^m_{++}$, we also define the weighted norm as $\norms{u}_{\beta}^2 := \sum_{j=1}^m \beta_j \norms{u_j}^2$.

	A multi-valued operator $T: \R^p \rightrightarrows \R^p$ is $\rho$-co-hypomonotone if there exists $\rho \in \R^m_+$ such that $\iprods{u - v, x - y} \geq -\rho \circ \norms{u - v}^2$ for all $(x,u), (y,v) \in \gra{T}$.
	$T$ is said to be maximally $\rho$-co-hypomonotone if it is $\rho$-co-hypomonotone and $\gra{T}$ is not contained properly in the graph of any $\rho$-co-hypomonotone operator.
	Note that a co-hypomonotone operator can be nonmonotone; see, e.g., \cite{tran2024accelerated,tran2025vfog} for concrete examples.
	
	\vspace{0.75ex}
	{\noindent\textbf{Coordinate-friendly computational setting.}
		For oracle-complexity purposes, we call the chosen partition block-coordinate friendly if evaluating the $m$ block outputs
		$[G_i(x^{(1)})]_1, \cdots, [G_i(x^{(m)})]_m$ at arbitrary admissible points has the same order of cost as one full evaluation of $G_i$.
		When the blocks have comparable computational sizes, one selected-block evaluation is therefore counted as $1/m$ of a full-component evaluation, on average.
		If $T$ is block-separable, i.e., $T(x) = (T_1 x_1, \cdots, T_m x_m)$, then so is the resolvent, i.e., $J_{\lambda \circ T}(u) = (J_{\lambda_1 T_1}(u_1), \cdots, J_{\lambda_m T_m}(u_m))$, thus a selected resolvent block can be evaluated independently.
		In Section~\ref{sec:ch2:VACOG}, we additionally require that
		$[G_i(J_{\lambda\circ T}u)]_j$ can be formed using selected-block operator/resolvent work without materializing the full resolvent vector.
		These conventions affect only the computational interpretation, not the	convergence analysis.
	}
	
	\vspace{0.75ex}
	\noindent\textbf{Randomized block sampling.}
	Let $\set{j_k}$ be an i.i.d. random sequence on $[m]$ representing the block selected at the $k$-th iteration, sampled according to
	\myeq{eq:VACOG_prob}{
		\arraycolsep=0.2em
		\begin{array}{lcl}
			\Prob{j_k = j} &=& \p_j > 0 \quad \text{for all } j \in [m],
		\end{array}
	}
	with $\sum_{j=1}^m \p_j = 1$.
	Thus, $\p \in \Delta_m$, the probability simplex in $\R^m$.
	Let $\underline{\p} := \min_{j \in [m]} \p_j$ and $\bar{\p} := \max_{j \in [m]} \p_j$.
	Throughout this paper, we assume $j_k$ is independent of the fresh randomness used to construct the estimator, conditional on $\mathcal{F}_k$.
	
	\beforesubsec
	\subsection{\textls[-20]{Optimistic gradient direction and its unbiased variance-reduced estimators}}
	\label{subsubsec:ch2:vr_estimators}
	\aftersubsec
	For a given sequence $\sets{x^k}_{k\geq 0}$ generated by an underlying algorithm, we define the \textit{optimistic gradient direction} as follows:
	\myeq{eq:Sk_operator}{
		\arraycolsep=0.2em
		\begin{array}{lcl}
			S^k & := & Gx^k - \gamma^k \circ Gx^{k-1},
		\end{array}
	}
	where $\gamma^k \in (0,1)^m$ is a given vector of coordinate parameters, determined later.
	Here, $S^k$ inherits the block-coordinate structure from $Gx$.
	
	Next, let $\widetilde{S}^k := \widetilde{S}^k(\Bc_k)$ be a stochastic estimator of $S^k$ in \eqref{eq:Sk_operator} constructed from a mini-batch $\Bc_k \subseteq [n]$.
	We introduce the following definition, which covers a class of  variance-reduced unbiased estimators of $S^k$.
	
	\begin{definition}\label{def:VR_estimator}
		For a given sequence $\sets{x^k}$ and a filtration $\sets{\Fc_k}$
		generated by an algorithm, a stochastic estimator $\widetilde{S}^k$ is called a \textit{variance-reduced unbiased estimator} of $S^k$ if there exist three constants $\kappa \in (0, 1]$, $\Theta > 0$, and $\hat{\Theta} \geq 0$, a nonnegative random sequence $\sets{\Delta_k}$, a nonnegative deterministic  sequence $\delta_k$, and a fixed vector $\beta \in \R^m_{++}$ such that for all $k \geq 1$, the following conditions hold almost surely:
		\myeq{eq:VR_estimator}{
			\hspace{-3ex}
			\left\{
			\arraycolsep=0.2em
			\begin{array}{lcl}
				\Expsk{k}{\widetilde{S}^k  - S^k} &=& 0, \vspace{0.5ex}\\
				\Expsk{k}{\norms{\widetilde{S}^k  - S^k }_{\beta}^2} &\leq& \Expsk{k}{\Delta_k}, \vspace{0.5ex} \\
				\Expsk{k}{\Delta_k} &\leq& (1 - \kappa)\norms{\nu^k}_{\infty}^2 \Delta_{k-1} + \Theta U_k + \hat{\Theta}\norms{\nu^k}_{\infty}^2 U_{k-1} +  \norms{1-\gamma^k}_{\infty}^2 \delta_k V_{k-1},
			\end{array}
			\right.
			\hspace{-5ex}
		}
		where $U_k := \frac{1}{n}\sum_{i=1}^n \norms{G_i x^k - G_i x^{k-1}}_{\beta}^2$, $V_k := \frac{1}{n}\sum_{i = 1}^n \norms{G_ix^k - Gx^k}_{\beta}^2$, and $\nu^k \in \R^m_{++}$ with $\nu_j^k := \frac{1 - \gamma_j^k}{1 - \gamma_j^{k-1}}$.
		Here, we initialize with $x^{-1} = x^0$ and $\widetilde{S}^0 := S^0 = (1 - \gamma^0)\circ Gx^0$, and set $\Delta_{-1} = \Delta_0 := 0$, $U_{-1} = U_0 := 0$, and $\delta_0 = 0$.
	\end{definition}
	The first condition in \eqref{eq:VR_estimator} is the unbiasedness of $\widetilde{S}^k$, while the last two conditions express its variance-reduction property.
	Two special choices of $\gamma^k$ will be used later.
	If $\gamma_j^k := \gamma_j \in (0, 1)$ for all $k \geq 0$, then the third line of \eqref{eq:VR_estimator} reduces to
	\myeq{eq:VR_estimator2}{
		\arraycolsep=0.2em
		\begin{array}{lcl}
			\Expsk{k}{\Delta_k} & \leq & (1 - \kappa) \Delta_{k-1} + \Theta U_k + \hat{\Theta}  U_{k-1} + \norms{1-\gamma}_{\infty}^2\delta_k V_{k-1}.
		\end{array}
	}
	If instead $\gamma_j^k := \gamma_k \in (0,1)$ for all $j \in [m]$, then it reduces to
	\myeq{eq:VR_estimator3}{
		\arraycolsep=0.2em
		\begin{array}{lcl}
			\Expsk{k}{\Delta_k} & \leq & (1 - \kappa) \nu_k^2 \Delta_{k-1} + \Theta U_k + \hat{\Theta} \nu_k^2 U_{k-1} + (1-\gamma_k)^2\delta_k V_{k-1}, \quad \nu_k := \frac{1-\gamma_k}{1 - \gamma_{k-1}}.
		\end{array}
	}
	Definition~\ref{def:VR_estimator} covers a class of variance-reduced unbiased estimators.
	Now, we present three estimators that satisfy this definition; their proofs are given in
	Appendix~\ref{apdx:sec:ch2:OG_operator}.
	
	\vspace{0.75ex}
	\noindent\textbf{$\mathrm{(a)}$~Mini-batch stochastic estimator.}
	Let $\Bc_k \subseteq [n]$ be sampled uniformly without replacement, independently of the history $\Fc_k$, with batch size $b_k := |\Bc_k|$.
	Define $G_{\Bc_k}x := \frac{1}{b_k} \sum_{i \in \Bc_k} G_i x$ and construct the mini-batch estimator of $S^k$:
	\myeq{eq:minibatch_estimator}{
		\arraycolsep=0.2em
		\begin{array}{lcl}
			\widetilde{S}^k &:=& G_{\Bc_k}x^k - \gamma^k \circ G_{\Bc_k}x^{k-1}.
		\end{array}
	}
	To quantify the finite-population sampling error, we define $\vartheta_k := \frac{n - b_k}{b_k(n-1)}$ if $n \geq 2$ and $\vartheta_k := 0$ if $n = 1$.
	Notice that $\vartheta_k = 0$ when $b_k = n$, in which case $\widetilde{S}^k = S^k$ (full-batch).
	\begin{lemma}\label{le:imb_estimator}
		Let $\sets{x^k}$  be updated by our algorithm and $\widetilde{S}^k$ be defined by \eqref{eq:minibatch_estimator}.
		Then, there exists $\sets{\Delta_k}\subset [0, +\infty)$ such that $\widetilde{S}^k$ satisfies Definition~\ref{def:VR_estimator} with $\kappa = 1$, $\Theta = 2$, $\hat{\Theta} = 0$, and $\delta_k := 2\vartheta_k$.
	\end{lemma}
	
	\noindent\textbf{$\mathrm{(b)}$~Loopless SVRG estimator.}
	Let $\Bc_k \subseteq [n]$ be a mini-batch of fixed size $b := |\Bc_k|$, sampled uniformly without replacement, independently across iterations and independently of the snapshot update \eqref{eq:loopless_SVRG_snapshot}.
	Given $\sets{x^k}$ and $S^k$ from \eqref{eq:Sk_operator}, we define
	\myeq{eq:loopless_SVRG}{
		\arraycolsep=0.2em
		\begin{array}{lcl}
			\widetilde{S}^k &:=& (1 - \gamma^k) \circ ( Gw^k - G_{\Bc_k} w^k ) + G_{\Bc_k} x^k - \gamma^k \circ G_{\Bc_k} x^{k-1},
		\end{array}
		\tag{L-SVRG}
	}
	where the reference point $w^k$ is initialized at $w^0 := x^0$ and then updated as
	\myeq{eq:loopless_SVRG_snapshot}{
		\arraycolsep=0.2em
		\begin{array}{lcl}
			w^k := \left\{
			\begin{array}{lcl}
				x^{k-1} \quad \text{with probability } \q, \vspace{0.5ex}\\
				w^{k-1} \quad \text{with probability } 1 - \q.
			\end{array}
			\right.
		\end{array}
	}
	Here, $\q \in (0, 1)$ is a given probability.
	Then, we have the following lemma.
	
	\begin{lemma}\label{lem:ch2:SVRG_estimator}
		Let $\sets{x^k}$  be updated by our algorithm, $S^k$ be defined by \eqref{eq:Sk_operator}, $\widetilde{S}^k$ be constructed by \eqref{eq:loopless_SVRG}, and $w^k$ be given by \eqref{eq:loopless_SVRG_snapshot}.
		Then, there exists $\sets{\Delta_k}\subset [0, +\infty)$ such that $\widetilde{S}^k$ satisfies Definition \ref{def:VR_estimator} with $\kappa := \frac{\q}{2} \in (0,1)$, $\Theta = \hat{\Theta} := \frac{4}{b\q}$, and $\delta_k = 0$.
	\end{lemma}

		\noindent$\mathrm{(c)}$~\textbf{A generalized JacSketch-type estimator.}
		Let $S_i^k := G_i x^k - \gamma^k \circ G_i x^{k-1}$ and $S^k = \frac{1}{n}\sum_{i=1}^n S_i^k$.
		Let $\mathcal{B}_k \subseteq [n]$ be a mini-batch of fixed size $b := |\mathcal{B}_k|$, sampled uniformly without replacement and independently across iterations.
		We maintain a table $\sets{H_i^k}_{i=1}^n$ of component directions, initialized at $H_i^0 := S_i^0$ for $i \in [n]$.
		All entries of the table are refreshed simultaneously: at each iteration, we draw a single Bernoulli random variable with success probability $\q \in (0,1)$, independently of $\Bc_k$, and set
		\myeq{eq:JacSketch_table}{
			\arraycolsep=0.2em
			H_i^k := \left\{\begin{array}{ll}
				S_i^k \ \text{ for all } i \in [n] & \text{with probability } \q, \vspace{1ex}\\
				\nu^k \circ H_i^{k-1} \ \text{ for all } i \in [n] & \text{with probability } 1 - \q.
			\end{array}
			\right.
		}
		Define $\bar{H}^k := \frac{1}{n}\sum_{i=1}^n H_i^k$,
		$H_{\Bc_k}^k := \frac{1}{b}\sum_{i\in\Bc_k}H_i^k$,
		$S_{\Bc_k}^k := \frac{1}{b}\sum_{i\in\Bc_k}S_i^k$.
		Then, we construct
		\myeq{eq:JacSketch_estimator}{
			\tag{GJS}
			\arraycolsep=0.2em
			\begin{array}{lcl}
				\widetilde S^k &:=& \bar{H}^k + S_{\Bc_k}^k - H_{\Bc_k}^k.
			\end{array}
		}
		A cumulative scaling vector can be used between table refreshes, so the entries need not be rescaled individually; when $\gamma^k\equiv\gamma$, we have $\nu^k = 1$, so no rescaling is needed.
		
		\begin{lemma}\label{le:GJS_estimator}
			Let $\{x^k\}$ be updated by our algorithm, $S^k$ be defined by \eqref{eq:Sk_operator}, and $\widetilde{S}^k$ be constructed by \eqref{eq:JacSketch_estimator}.
			Then, there exists $\sets{\Delta_k} \subset [0,+\infty)$ such that $\widetilde{S}^k$ satisfies Definition~\ref{def:VR_estimator} with $\kappa := \frac{\mbf{q}}{2}\in(0,1)$, $\Theta = \hat{\Theta} := \frac{4}{b\mbf{q}}$, and $\delta_k = 0$.
		\end{lemma}
	
	\beforesec
	\section{Variance-Reduced Block-Coordinate OG Method for Weak-Minty \eqref{eq:FNI}}\label{sec:ch2:VCOG_alg}
	\aftersec
	We first develop a \textit{\textbf{V}ariance-reduced Block-\textbf{C}oordinate \textbf{O}ptimistic \textbf{G}radient method} (VCOG) to solve a class of \eqref{eq:FNI} under the following assumptions.
	
	\begin{assumption}\label{as:VCOG}
		The composite inclusion~\eqref{eq:FNI} satisfies the following conditions:
		\begin{enumerate}
		\itemsep=0.35em
			\item[$\mathrm{(i)}$] (\textbf{\textit{Weak-Minty solution}})
			There exist $\rho \in \R^m_{+}$ and a solution $\xopt \in \zer{\Psi}$ such that $\iprods{Gx + \xi, x - \xopt} \geq -\rho \circ \norms{Gx + \xi}^2$ for all $x \in \dom{G}$ and $\xi \in Tx$.
			
			\item[$\mathrm{(ii)}$] (\textbf{\textit{Lipschitz continuity}})
			$G$ in $\eqref{eq:FNI}$ is $L$-block-wise average Lipschitz continuous with $L \in \R^m_{++}$, i.e., for all $x, y\in \dom{G}$, we have
			\myeq{eq:VCOG_Lipschitz_continuity}{
				\arraycolsep=0.2em
				\begin{array}{lcl}
					\frac{1}{n} \sum_{i=1}^n \sum_{j=1}^m \frac{1}{L_j^2} \big\Vert [G_ix]_j - [G_iy]_j \big\Vert^2 &\leq& \norms{x - y}^2.
				\end{array}
			}
			
			\item[$\mathrm{(iii)}$] (\textbf{\textit{Lipschitz resolvent}})
			For the parameter vector $\gamma\eta := (\gamma_1\eta_1, \cdots, \gamma_m\eta_m)$ used in Algorithm~\ref{alg:A1}, the resolvent $J_{\gamma\eta \circ T} := (\Id + \gamma\eta \circ T)^{-1}$ is single-valued on $\R^p$, and there exists $\mu > 0$ such that
			\myeq{eq:VCOG_Lipschitz_resolvent}{
				\arraycolsep=0.2em
				\begin{array}{lcl}
					\norms{J_{\gamma\eta \circ T}(u) - J_{\gamma\eta \circ T}(v)} &\leq& \mu \norms{u-v}, \qquad \forall u,v\in\R^p.
				\end{array}
			}
		\end{enumerate}
	\end{assumption}
	Assumption~\ref{as:VCOG}(i) is a block-wise extension of the standard
	weak-Minty solution condition~\cite{diakonikolas2021efficient}; when
	$\rho_1=\cdots=\rho_m=\rho$, it reduces to its usual scalar form.
	Let $\zerM{\Psi}$ denote the set of $\rho$-weak-Minty solutions of
	\eqref{eq:FNI}. Clearly, $\zerM{\Psi}\subseteq\zer{\Psi}$.
	If $\Psi$ is $\rho$-co-hypomonotone,  then every solution is $\rho$-weak-Minty; hence $\zerM{\Psi} = \zer{\Psi}$.
	Thus, co-hypomonotonicity provides a sufficient condition under which the weak-Minty solution used in our analysis can be chosen as an arbitrary solution of \eqref{eq:FNI}.
	The case $\rho=0$ recovers the standard weak-Minty (star-monotone) condition and, in particular, includes monotone $\Psi$.
	Several concrete nonmonotone examples can be found in~\cite{TranDinh2025a}.
	
	By Jensen's inequality, Assumption~\ref{as:VCOG}(ii) implies
	\myeqn{
		\arraycolsep=0.2em
		\begin{array}{lcl}
			\sum_{j=1}^m \frac{1}{L_j^2} \norms{[Gx ]_j - [Gy]_j}^2 &\leq& \norms{x - y}^2, \quad \forall x, y \in \dom{G},
		\end{array}
	}
	so $G$ is also $L$-block-wise Lipschitz continuous.
	Such average-Lipschitz conditions are standard in variance-reduced
	block-coordinate optimization; see, e.g.,
	\cite{cai2023cyclic,hanzely2020variance}.
	
	Assumption~\ref{as:VCOG}(iii) is mild.
	It is a regularity condition on the resolvent of $T$ at the step-size vector used by the algorithm, and it does not require $T$ or $\Psi$ to be monotone.
	When $T$ is block-separable as in Assumption~\ref{as:VCOG_T_separable} below, $J_{\gamma\eta\circ T}$ acts blockwise through $J_{\gamma_j\eta_j T_j}$.
	Below are three concrete examples in which \eqref{eq:VCOG_Lipschitz_resolvent} holds.
	\begin{compactitem}
		\item[(a)] If each $T_j$ is maximally monotone, then \eqref{eq:VCOG_Lipschitz_resolvent} holds with $\mu = 1$ for every $\gamma\eta$.
		\item[(b)] If each $T_j$ is maximally $\hat{\rho}_j$-co-hypomonotone and $\gamma_j\eta_j > \hat{\rho}_j$, then \eqref{eq:VCOG_Lipschitz_resolvent} holds with $\mu = \max_j\max\big\{1, \frac{\hat{\rho}_j}{\gamma_j\eta_j - \hat{\rho}_j}\big\}$, and $\mu = 1$ if $\gamma_j\eta_j \geq 2\hat{\rho}_j$ for all $j \in [m]$.
		The step-sizes of Theorem~\ref{th:convergence_of_Alg1} below satisfy $\gamma_j\eta_j \geq \frac{3\tau}{\p_j} > 6\rho_j$, so $\mu = 1$ whenever $\hat{\rho}_j \leq 3\rho_j$.
		\item[(c)] If each $T_j$ is maximally $\hat{\rho}_j$-hypomonotone (i.e., $\iprods{u - v, x - y} \geq -\hat{\rho}_j\norms{x - y}^2$ for all $(x,u), (y,v) \in \gra{T_j}$) and $\gamma_j\eta_j\hat{\rho}_j < 1$, then \eqref{eq:VCOG_Lipschitz_resolvent} holds with $\mu = \max_j \frac{1}{1 - \gamma_j\eta_j\hat{\rho}_j}$.
		In particular, if $\hat{\rho}_j \leq \sqrt{2}L_j$ for all $j \in [m]$, then Theorem~\ref{th:convergence_of_Alg1} can be applied with $\mu := 2$.
		Indeed, its step-size conditions with $\mu = 2$ give $\gamma_j\eta_j \leq \frac{\tau}{\lambda\p_j} \leq \frac{1}{4\sqrt{2}L_j}$, hence $\gamma_j\eta_j\hat{\rho}_j \leq \frac{1}{4}$, and $J_{\gamma\eta\circ T}$ is $\frac{4}{3}$-Lipschitz.
	\end{compactitem}
	Without block-separability, a maximally monotone $T$ still satisfies \eqref{eq:VCOG_Lipschitz_resolvent} with $\mu = \big(\frac{\max_j\gamma_j\eta_j}{\min_j\gamma_j\eta_j}\big)^{1/2}$, since $J_{\gamma\eta\circ T}$ is nonexpansive in the norm weighted by $(\gamma\eta)^{-1}$.
	
	\beforesubsec
	\subsection{The proposed algorithm and its implementable version}\label{subsec:ch2:proposed_alg}
	\aftersubsec
	Throughout this section, we instantiate the variance-reduction conditions in Definition~\ref{def:VR_estimator} with $\beta_j := \frac{1}{L_j^2}$, $j \in [m]$.
	For the optimistic gradient direction $S^k$ defined by \eqref{eq:Sk_operator} and its estimator $\widetilde{S}^k$ satisfying Definition~\ref{def:VR_estimator}, we define $\widetilde{S}^k_{[j]}  :=  [0, \cdots, 0, [\widetilde{S}^k]_j, 0, \cdots, 0]$,
	then Algorithm~\ref{alg:A1} presents our first method for solving \eqref{eq:FNI}.
	\setlength{\textfloatsep}{10pt plus 1.0pt minus 3.0pt}
	\begin{algorithm}[hpt!]\caption{(\textit{\textbf{The VCOG framework for solving \eqref{eq:FNI} under Assumption~\ref{as:VCOG}}})}\label{alg:A1}
		\normalsize
		\begin{algorithmic}[1]
			\STATE\label{step:A1_i0}{\bfseries Initialization:}  
			Choose $x^0 \in \R^p$, $\eta \in \R^m_{++}$, $\gamma \in (0,1)^m$, and $k_{\max} > 0$.
			Set  $x^{-1} := x^0$, compute $\widetilde{S}^0 := (1-\gamma) \circ Gx^0$, and $y^0 := x^0 + \gamma \eta \circ \xi^0$, where $\xi^0 \in Tx^0$.
			\STATE\hspace{0ex}\label{step:A1_loop}{\bfseries For $k = 0$ to $k_{\max}$ do}
			\vspace{0.25ex}
			\STATE\hspace{1.5ex}\label{step:A1_o1}Sample a block $j_k \in [m]$ following the probabilistic rule \eqref{eq:VACOG_prob}.
			\STATE\hspace{1.5ex}\label{step:A1_o2}If $k \geq 1$, construct $\widetilde{S}^k$ satisfying Definition~\ref{def:VR_estimator},  independent of $j_k$.
			\STATE\hspace{1.5ex}\label{step:A1_o3}
			Update the iterates
			\myeq{eq:VCOG4}{
				\arraycolsep=0.2em
				\left\{\begin{array}{lcl}
					y^{k+1} & := & y^k - \frac{\eta_{j_k}}{\p_{j_k}} \widetilde{S}^k_{[j_k]} - \frac{(1 - \gamma_{j_k})\eta_{j_k}}{\p_{j_k}} \xi^k_{[j_k]}, \vspace{1ex} \\
					x^{k+1} & := & J_{\gamma \eta \circ T}(y^{k+1}), \vspace{1ex}\\
					{}[\xi^{k+1}]_j & := &  \frac{y_j^{k+1} - x_j^{k+1}}{\gamma_j \eta_j}, \quad j \in [m].
				\end{array}\right.
			}
			\STATE\hspace{0ex}{\bfseries End For}
		\end{algorithmic}
	\end{algorithm}
	
	To motivate Scheme~\eqref{eq:VCOG4}, we first consider the case $T = 0$.
	Replacing the deterministic OG step by a randomized block update with the stochastic direction $\widetilde S^k$ gives
	\myeq{eq:VCOG3}{
		\arraycolsep=0.2em
		\begin{array}{lcl}
			x^{k+1} := x^k - \eta_{j_k}\p_{j_k}^{-1}\widetilde{S}^k_{[j_k]}.
		\end{array}
	}
	Then, to accommodate the case $T\neq 0$, we adopt the idea of forward-backward splitting methods and modify \eqref{eq:VCOG3} as \eqref{eq:VCOG4},
	where the term $\xi^k \in Tx^k$ on the first line is added to handle the operator $T$, and the second line is a backward step.
	
	\vspace{0.75ex}
	\noindent\textbf{Implementable version.}
	The latter step is full-dimensional in general.
	However, a complete block-coordinate implementation is available under the following assumption.
	
	\begin{assumption}[\textbf{Separability of $T$}]\label{as:VCOG_T_separable}
		The operator $T$ in \eqref{eq:FNI} is block-separable, i.e.,
		$Tx = \big[ T_1x_1, T_2x_2, \cdots, T_m x_m \big]$, for any $x = [ x_1,  \cdots, x_m]$.
	\end{assumption}
	\noindent
	Note that this assumption is required only for the block-coordinate implementation, not for convergence analysis.
	It covers standard separable regularizers and constraints, including $\ell_1$, $\ell_{2,1}$, and SCAD regularizers, and the indicator of box constraints.
	Then, the block-coordinate implementation of Algorithm~\ref{alg:A1} can be derived as follows.

	Let us define the full vectors $y^k := x^k + \gamma \eta \circ \xi^k$ and $x^k = J_{\gamma \eta \circ T}(y^k)$.
	Using the first line of \eqref{eq:VCOG4} and Assumption~\ref{as:VCOG_T_separable}, we obtain the following scheme.
	\begin{compactenum}
		\item[1.] \textit{Starting from $x^0 \in \R^p$, we set $x^{-1} := x^0$, $\xi^0 \in Tx^0$, and $y^0 := x^0 + \gamma \eta \circ \xi^0$.}
		\item[2.] \textit{At each iteration $k \geq 0$, we sample a block $j_k$ according to \eqref{eq:VACOG_prob}, then update}
		\myeqn{
			\arraycolsep=0.2em
			\left\{\begin{array}{lcl}
				y^{k+1}_{j_k} &:=& y^k_{j_k} - \frac{\eta_{j_k}}{\p_{j_k}}[\widetilde{S}^k]_{j_k} - \frac{(1-\gamma_{j_k})\eta_{j_k}}{\p_{j_k}}[\xi^k]_{j_k}
				\vspace{1ex}\\
				x^{k+1}_{j_k} &:=& J_{\gamma_{j_k} \eta_{j_k} T_{j_k}}(y^{k+1}_{j_k}),
				\vspace{1ex}\\
				\xi^{k+1}_{j_k} &:=& \frac{1}{\gamma_{j_k} \eta_{j_k}} (y^{k+1}_{j_k} - x^{k+1}_{j_k}),
			\end{array}\right.
		}
		\textit{and keep $(y_j^{k+1},x_j^{k+1},\xi_j^{k+1}) = (y_j^k,x_j^k,\xi_j^k)$ unchanged for $j\neq j_k$.}
		\vspace{0.75ex}
	\end{compactenum}
	This formulation makes the block-coordinate structure explicit: at each iteration, only the sampled block $[\widetilde S^k]_{j_k}$ is needed, so the full estimator $\widetilde S^k$ need not be formed.
	Moreover, under Assumption~\ref{as:VCOG_T_separable}, we also need to evaluate only one block $J_{\gamma_{j_k}\eta_{j_k} T_{j_k}}(y_{j_k}^{k+1})$.
	This substantially reduces the computational cost when the problem dimension and the number of blocks are large.
	
	\beforesubsec
	\subsection{Convergence analysis}
	\aftersubsec
	For convenience of convergence analysis, we eliminate $y^k$ from \eqref{eq:VCOG4} and obtain the following equivalent scheme:
	\myeq{eq:VCOG}{
		\arraycolsep=0.2em
		\begin{array}{lcl}
			x^{k+1} &:=& x^k - \eta \gamma \circ (\xi^{k+1} - \xi^k) - \frac{\eta_{j_k}}{\p_{j_k}}\widetilde{S}^k_{[j_k]}  - \frac{(1-\gamma_{j_k}) \eta_{j_k}}{\p_{j_k}}\xi^k_{[j_k]}.
		\end{array}
		\tag{VCOG} 
	}
	Then, we are ready to analyze the convergence of Algorithm~\ref{alg:A1}.
	
	\vspace{0.75ex}
	\noindent\textbf{$\mathrm{(a)}$~Lyapunov function.}
	For a given $\xi^k \in Tx^k$, we first introduce:
	\myeq{eq:VCOG_w_notations}{
		\arraycolsep=0.2em
		\begin{array}{l}
			w^k :=  Gx^k + \xi^k,  \qquad \hat{w}^k :=  Gx^{k-1} + \xi^k, \qquad \text{and} \qquad e^k := \widetilde{S}^k - S^k.
		\end{array}
	}
	Next, for $\kappa$, $\hat{\Theta}$, $\Delta_k$, and $U_k$ in Definition~\ref{def:VR_estimator}, we define the following functions:
	\myeq{eq:VCOG_potential_func}{
		\arraycolsep=0.2em
		\begin{array}{lcl}
			\Lc_k &:=& \norms{x^k - \xopt + \gamma \eta \circ (Gx^{k-1} + \xi^k)}^2 + \frac{1}{\mu^2}\norms{x^k - x^{k-1}}^2, \vspace{1ex}\\
			\Ec_k &:=& \Lc_k + \frac{2M_1(1-\kappa)}{\kappa} \Delta_{k-1} + \frac{2M_1\hat{\Theta}}{\kappa} U_{k-1},
		\end{array}
	}
	where $\gamma \in (0, 1)^m$ and $\eta \in \R^m_{++}$ are given, determined later, and $M_1 := \max_{j \in [m]} \frac{L_j^2 \eta_j^2}{\p_j}$.
	The function $\Ec_k$ can be viewed as a Lyapunov function for Algorithm~\ref{alg:A1}.
	
	For later use, we define $V(x) := \frac{1}{n}\sum_{i = 1}^n \norms{G_ix - Gx}_{\beta}^2 \geq 0$ and $\sigma_{\star}^2 := V(x^{\star}) < +\infty$.
	Denoting $\hat{e}_i(x) := G_ix - Gx$, then using Minkowski's inequality, we have
	\myeqn{
		\arraycolsep=0.2em
		\begin{array}{lcl}
			\sqrt{V(x)}  &\leq& \big(\frac{1}{n}\sum_{i = 1}^n \norms{\hat{e}_i(x) - \hat{e}_i(x^{\star})}_{\beta}^2\big)^{1/2} + \big(\frac{1}{n} \sum_{i=1}^{n} \norms{\hat{e}_i(x^{\star})}_{\beta}^2 \big)^{1/2}.
		\end{array}
	}
	The second term is exactly $\sigma_{\star}$.
	Using $\beta_j = \frac{1}{L_j^2}$ for all $j \in [m]$ and Assumption~\ref{as:VCOG}(ii), the first term can be upper bounded as
	\myeqn{
		\arraycolsep=0.2em
		\begin{array}{lcl}
			\frac{1}{n}\sum_{i = 1}^n \norms{\hat{e}_i(x) - \hat{e}_i(x^{\star})}_{\beta}^2 &=& \frac{1}{n}\sum_{i = 1}^n \norms{(G_ix - G_ix^{\star}) - (Gx - Gx^{\star})}_{\beta}^2 \vspace{1ex}\\
			&=& \frac{1}{n}\sum_{i=1}^n \norms{G_ix - G_ix^{\star}}_{\beta}^2 - \norms{Gx - Gx^{\star}}_{\beta}^2 \vspace{1ex}\\
			&\leq& \frac{1}{n}\sum_{i=1}^n \sum_{j=1}^m \frac{1}{L_j^2} \norms{[G_ix]_j - [G_ix^{\star}]_j}^2 \vspace{1ex}\\
			&\leq& \norms{x - x^{\star}}^2.
		\end{array}
	}
	Thus, we have $\sqrt{V(x)} \leq \norms{x - x^{\star}} + \sigma_{\star}$, leading to
	\myeq{eq:VCOG_Vk_bound}{
		\arraycolsep=0.2em
		\begin{array}{lcl}
			V_k = V(x^k) &\leq& 2\norms{x^k - x^{\star}}^2 + 2\sigma_{\star}^2 \qquad \text{for all } k \geq 0.
		\end{array}
	}
		
	\noindent\textbf{$\mathrm{(b)}$~Key estimates.}
	First, we state the following key lemma, whose proof is deferred to Appendix~\ref{apdx:VCOG_convergence}.
	\begin{lemma}\label{lem:ch2:VCOG_key_lemma1}
		Suppose that Assumption~\ref{as:VCOG} holds for \eqref{eq:FNI}.
		Let $\sets{x^k}$ be generated by Algorithm~\ref{alg:A1} and $\Lc_k$ be defined by \eqref{eq:VCOG_potential_func} for an arbitrary $x^{\star} \in \zerM{\Psi}$.
		Given $\tau > 0$, if
		\myeq{eq:VCOG_para_cond1}{
			\arraycolsep=0.2em
			\begin{array}{lcl}
				(1 - \gamma_j)\eta_j &=& \tau \quad \text{and} \quad 2\mu^2 L_j^2 \eta_j^2 \gamma_j^2 (2-\p_j) < \p_j \quad \text{for all } j \in [m],
			\end{array}
		}
		then, almost surely, we have
		\myeq{eq:VCOG_key_lemma1}{
			\arraycolsep=0.2em
			\hspace{-2ex}
			\begin{array}{lcl}
				\Lc_k - \Expsk{k}{\Lc_{k+1}} & \geq &  \frac{1}{2\mu^2}\norms{x^k - x^{k-1}}^2 - \sum_{j=1}^m \frac{2\eta_{j}^2}{\p_{j}}  \Expsk{k}{\norms{[e^k]_{j}}^2} + \sum_{j=1}^m \psi_j \norms{[w^k]_j}^2,
			\end{array}
			\hspace{-2ex}
		}
		where $\psi_j := \frac{2(1-\gamma_j)\eta_j^2[\gamma_j \p_j - (1 - \gamma_j)]}{\p_j} - 2\tau \rho_j - \frac{4(1-\gamma_j)^2 \gamma_j^2 \eta_j^4}{\p_j^2 \phi_j}$ and $\phi_j := \frac{1}{2\mu^2 L_j^2} - \frac{(2-\p_j)\gamma_j^2 \eta_j^2}{\p_j}$.
	\end{lemma}
	Given this key estimate, we can next prove a quasi-descent property of the Lyapunov function $\Ec_k$ defined by \eqref{eq:VCOG_potential_func}.
	The proof of this lemma can be found in Appendix~\ref{apdx:VCOG_convergence}.
	\begin{lemma}\label{lem:ch2:VCOG_key_lemma2}
		Under the hypotheses of Lemma~\ref{lem:ch2:VCOG_key_lemma1}, let $\Ec_k$ be defined by \eqref{eq:VCOG_potential_func} using an arbitrary $x^{\star} \in \zerM{\Psi}$, and let $\widetilde{S}^k$ be a stochastic estimator of $S^k$ satisfying Definition~\ref{def:VR_estimator} with $\gamma^k \equiv \gamma \in (0,1)^m$.
		Let  $\hat{\delta}_k := \frac{128\mu^2 M_1 \norms{1 - \gamma}_{\infty}^2}{\kappa} \delta_k$.
		Then, almost surely, we have
		\myeq{eq:VCOG_key_lemma2}{
			\hspace{-2ex}
			\arraycolsep=0.2em
			\begin{array}{lcl}
				\textls[-20]{$\Expsk{k}{\Ec_{k+1}} \leq (1 + \hat{\delta}_k)\Ec_k -  \Big[\frac{1}{2\mu^2} - \frac{2M_1(\Theta + \hat{\Theta})}{\kappa}\Big] \norms{x^k - x^{k-1}}^2 - \sum\limits_{j=1}^m \psi_j \norms{[w^k]_j}^2 + \frac{\sigma_{\star}^2}{32\mu^2}\hat{\delta}_k.$}
			\end{array}
			\hspace{-2ex}
		}
	\end{lemma}
	
	\noindent$\mathrm{(c)}$~\textbf{Main result 1.}
	Given $\kappa$, $\Theta$, and $\hat{\Theta}$ in Definition~\ref{def:VR_estimator}, let $\omega^2 := \min\set{1, \frac{\kappa}{\Theta + \hat{\Theta}}}$.
	For a given sampling distribution $\p$ defined in \eqref{eq:VACOG_prob}, we define $L_{\p} := \norms{\p^{-3/2}L}_{\infty}$, and $\rho_{\p} = \norms{\p\rho}_{\infty}$.
	Then, the convergence of Algorithm~\ref{alg:A1} is given in the following theorem.
	
	\begin{theorem}\label{th:convergence_of_Alg1}
		For \eqref{eq:FNI}, suppose that Assumption~\ref{as:VCOG} holds.
		Let $\sets{x^k}$ be generated by Algorithm~\ref{alg:A1} using $\gamma_j = 1 - \lambda \p_j$ and $\eta_j = \frac{\tau}{\lambda \p_j}$, where
		\myeq{eq:VCOG_param_conditions}{
			\arraycolsep=0.2em
			\begin{array}{lll}
				\lambda \in \big(0, \frac{1}{4}\big], \quad
				L_{\p} \rho_{\p} < \frac{\lambda \omega}{4\sqrt{2}\mu}, \quad \text{and} \quad
				2\rho_{\p} < \tau \leq \frac{\lambda \omega}{2\sqrt{2}\mu L_{\p}}.
			\end{array}
		}
		Then, the following statement holds.
		\begin{compactitem}
		
		\item[$\mathrm{(a)}$] \textbf{\textit{$($Bounds in expectation$)$}}
		For $\xi^k \in Tx^k$ generated by Algorithm~\ref{alg:A1}, we have
		\myeq{eq:VCOG_convergence_2}{
				{\frac{1}{K+1}\cdot\sum_{k=0}^K \Expk{\norms{Gx^k + \xi^k}^2} } \leq \frac{\bar{\p} \Rc_K^2}{\tau^2 (K+1)},
		}
		where $\Rc_K^2 := \big(\Rc_0^2 + \frac{\sigma_{\star}^2}{32\mu^2}\big) \exp(\Sc_K) - \frac{\sigma_{\star}^2}{32\mu^2}$ with $\Rc_0^2 := 2\norms{x^0 - \xopt}^2 + 2\gamma^2 \eta^2 \circ \norms{Gx^0 + \xi^0}^2$ and $\Sc_K := \sum_{k=0}^K \hat{\delta}_k$.
		
		\item[] 
		In particular, if $\sum_{k=0}^{\infty} \delta_k < +\infty$, then
		\myeq{eq:VCOG_convergence_2.1}{
				{\frac{1}{K+1}\cdot\sum_{k=0}^K \Expk{\norms{Gx^k + \xi^k}^2} } \leq \frac{\bar{\p}\Rc_{\infty}^2}{\tau^2 (K+1)},
		}
		where $\Rc_{\infty}^2 := \big(\Rc_0^2 + \frac{\sigma_{\star}^2}{32\mu^2}\big) \exp(\Sc_{\infty}) - \frac{\sigma_{\star}^2}{32\mu^2}$ and $\Sc_{\infty} := \sum_{k=0}^{\infty} \hat{\delta}_k$.

		\item[$\mathrm{(b)}$]\textbf{\textit{$($Almost-sure convergence$)$}}
		If $\sum_{k=0}^{\infty} \delta_k < +\infty$, then almost surely,
		\myeq{eq:VCOG_convergence_3}{
			\arraycolsep=0.2em
			\begin{array}{lcl}
				\lim\limits_{k \to \infty} \norms{Gx^k + \xi^k} &=& 0, \quad \xi^k \in Tx^k.
			\end{array}
		}
		If, in addition, $\gra{\Psi}$ is closed and $\zer{\Psi} = \zerM{\Psi}$ $($in particular, if $\Psi$ is $\rho$-co-hypomonotone$)$, then $\sets{x^k}$ converges almost surely to a $\zer{\Psi}$-valued random variable, a solution of \eqref{eq:FNI}.
		\end{compactitem}
	\end{theorem}
	
	\begin{proof}
		We first show that the parameters are well defined, and then prove the convergence in expectation and the almost-sure convergence.
		\vspace{0.75ex}
		
		\noindent\textbf{$\mathrm{(i)}$~\textit{Parameter bounds}.}
		First, let us choose $\gamma_j$ and $\eta_j$ such that $1 - \gamma_j = \lambda \p_j$ and $\eta_j = \frac{\tau}{\lambda \p_j}$ for some $\lambda \in \big(0, \frac{1}{4}\big]$.
		Let $\omega^2 := \min\set{1, \frac{\kappa}{\Theta + \hat{\Theta}}}$, and we impose $M_1 \leq \frac{\omega^2}{8\mu^2}$.
		Using this bound, $\lambda \in \big(0, \frac{1}{4}\big]$, and $\gamma_j \in (0,1)$, we can compute that
		\myeq{eq:VCOG_theorem_proof1}{
			\arraycolsep=0.2em
			\left\{\begin{array}{lcl}
				\chi &:=& \frac{1}{2\mu^2} - \frac{2M_1(\Theta + \hat{\Theta})}{\kappa} \geq \frac{1}{2\mu^2} - \frac{1}{4\mu^2} = \frac{1}{4\mu^2}, \vspace{1ex}\\
				
				\phi_j &=& \frac{1}{2\mu^2 L_j^2} - \frac{(2-\p_j)\gamma_j^2 \eta_j^2}{\p_j}
				\geq \frac{1}{2\mu^2 L_j^2} - \frac{2 \eta_j^2}{\p_j}
				\geq \frac{1}{2\mu^2 L_j^2} - \frac{2M_1}{L_j^2} \geq \frac{1}{4\mu^2 L_j^2}, \vspace{1ex}\\
				\psi_j &=& \frac{2\tau^2}{\lambda \p_j}[1 - \lambda(\p_j + 1)] - 2\tau \rho_j - \frac{4\tau^2 \gamma_j^2 \eta_j^2}{\p_j^2 \phi_j}
				\geq \frac{4\tau^2}{\p_j} - 2\tau \rho_j - \frac{4\tau^2 M_1}{\p_j \phi_j L_j^2}
				\geq \frac{\tau^2}{\p_j}.
			\end{array}
			\right.
		}
		In the last line of \eqref{eq:VCOG_theorem_proof1}, we have used $\frac{4M_1}{\phi_jL_j^2} \leq 2$, which follows from $M_1 \leq \frac{1}{8\mu^2}$ and $\phi_j \geq \frac{1}{4\mu^2L_j^2}$.
		This gives $\psi_j \geq \frac{2\tau^2}{\p_j} - 2\tau\rho_j$, and the right-hand side is at least $\frac{\tau^2}{\p_j}$ whenever $\tau \geq 2\p_j\rho_j$.
		Moreover, since $M_1 = \max_{j \in [m]} \frac{L_j^2 \eta_j^2}{\p_j} = \frac{\tau^2}{\lambda^2} \max_{j \in [m]} \frac{L_j^2}{\p_j^3}$, the condition $M_1 \leq \frac{\omega^2}{8\mu^2}$ holds if $\tau \leq \frac{\lambda \omega}{2\sqrt{2}\mu} \frac{\p_j^{3/2}}{L_j}$ for all $j \in [m]$.
		This leads to the choice $2\rho_{\p} < \tau \leq \frac{\lambda \omega}{2\sqrt{2}\mu L_{\p}}$.
		To guarantee the existence of a common $\tau$, we require $L_{\p} \rho_{\p} < \frac{\lambda \omega}{4\sqrt{2}\mu}$.
				
		\vspace{0.75ex}
		\noindent\textbf{$\mathrm{(ii)}$~\textit{Convergence in expectation}.}
		Now, using the bounds in \eqref{eq:VCOG_theorem_proof1} and $\bar{\p} := \max_{j \in [m]} \p_j$, we obtain from \eqref{eq:VCOG_key_lemma2} that
		\myeq{eq:VCOG_theorem_proof2}{
			\arraycolsep=0.2em
			\begin{array}{lcl}
				\Expsk{k}{\Ec_{k+1}} &\leq& (1+\hat{\delta}_k)\Ec_k - \frac{1}{4\mu^2} \norms{x^k - x^{k-1}}^2 - \frac{\tau^2}{\bar{\p}} \norms{Gx^k + \xi^k}^2 + \frac{\sigma_{\star}^2}{32\mu^2}\hat{\delta}_k.
			\end{array}
		}
		Taking the total expectation on both sides of this inequality, we obtain
		\myeqn{
			\arraycolsep=0.2em
			\begin{array}{lcl}
				\Expk{\Ec_{k+1}} &\leq& (1+\hat{\delta}_k)\Expk{\Ec_k} -  \frac{1}{4\mu^2} \Expk{\norms{x^k - x^{k-1}}^2} - \frac{\tau^2}{\bar{\p}} \Expk{\norms{Gx^k + \xi^k}^2} + \frac{\sigma_{\star}^2}{32\mu^2}\hat{\delta}_k.
			\end{array}
		}
		
		Let $A_0 := 1$ and $A_{k+1} = A_k(1+\hat{\delta}_k) = \prod_{l=0}^{k}(1+\hat{\delta}_l) > 0$.
		Since $\frac{1}{A_k} - \frac{1}{A_{k+1}} = \frac{\hat{\delta}_k}{A_{k+1}}$, dividing the last inequality by $A_{k+1}$ and defining $\Qc_k := \Ec_k + \frac{\sigma_{\star}^2}{32\mu^2}$, we obtain
		\myeqn{
			\arraycolsep=0.2em
			\begin{array}{lcl}
				\frac{\Expk{\Qc_{k+1}}}{A_{k+1}} &\leq& \frac{\Expk{\Qc_k}}{A_k} - \frac{\Expk{\norms{x^k - x^{k-1}}^2}}{4\mu^2 A_{k+1}} - \frac{\tau^2 \Expk{\norms{Gx^k + \xi^k}^2}}{\bar{\p}A_{k+1}}.
			\end{array}
		}
		Since $\Qc_{K+1} \geq \frac{\sigma_{\star}^2}{32\mu^2}$ and $A_{k+1} \leq A_{K+1}$ for $0 \leq k \leq K$, summing this inequality from $k=0$ to $k=K$ and multiplying the result by $A_{K+1}$ gives
		\myeq{eq:VCOG_theorem_proof2.1}{
			\arraycolsep=0.2em
			\begin{array}{lcl}
				\sum_{k = 0}^K \Expk{\norms{x^k - x^{k-1}}^2} &\leq& 4\mu^2 \big[\Ec_0 A_{K+1} + \frac{\sigma_{\star}^2}{32\mu^2} (A_{K+1} - 1)\big], \vspace{1ex}\\
				\sum_{k = 0}^K \Expk{\norms{Gx^k + \xi^k}^2} &\leq& \frac{\bar{\p}}{\tau^2}\big[\Ec_0 A_{K+1} + \frac{\sigma_{\star}^2}{32\mu^2} (A_{K+1} - 1)\big].
			\end{array}
		}
		By the initialization, $x^{-1} = x^0$ and $\Delta_{-1} = U_{-1} = 0$.
		Hence, by Young's inequality,
		\myeqn{
			\arraycolsep=0.2em
			\begin{array}{lcl}
				\Ec_0 &=& \norms{x^0 - \xopt + \gamma \eta \circ (Gx^0 + \xi^0)}^2 \leq 2\norms{x^0 - \xopt}^2 + 2\gamma^2 \eta^2 \circ \norms{Gx^0 + \xi^0}^2 := \Rc_0^2.
			\end{array}
		}
		Moreover, denoting $\Sc_k := \sum_{l=0}^k \hat{\delta}_l$, we also have
		\myeqn{
			\arraycolsep=0.2em
			\begin{array}{lcl}
				A_{K+1} = \prod_{k=0}^{K}(1+\hat{\delta}_k) \leq \prod_{k=0}^K \exp(\hat{\delta}_k) = \exp(\Sc_K)
			\end{array}
		}
		Substituting the last two bounds into \eqref{eq:VCOG_theorem_proof2.1} and defining $\Rc_K^2 := \big(\Rc_0^2 + \frac{\sigma_{\star}^2}{32\mu^2}\big) e^{\Sc_K} - \frac{\sigma_{\star}^2}{32\mu^2}$, we obtain
		\myeq{eq:VCOG_theorem_proof3}{
			\arraycolsep=0.2em
			\begin{array}{lcl}
				\frac{1}{K+1}\sum_{k = 0}^K \Expk{\norms{x^k - x^{k-1}}^2} &\leq& \frac{4\mu^2 \Rc_K^2}{K+1}, \vspace{1ex}\\
				\frac{1}{K+1}\sum_{k = 0}^K \Expk{\norms{Gx^k + \xi^k}^2} &\leq& \frac{\bar{\p} \Rc_K^2}{\tau^2 (K+1)}. \vspace{1ex}\\
			\end{array}
		}
		The second line of \eqref{eq:VCOG_theorem_proof3} is exactly \eqref{eq:VCOG_convergence_2}.
		If additionally $\sum_{k = 0}^\infty \delta_k < +\infty$, then $\Sc_K \leq \Sc_{\infty} < +\infty$, so $\Rc_K$ is uniformly bounded, and \eqref{eq:VCOG_convergence_2.1} follows directly from  \eqref{eq:VCOG_convergence_2}.
		
		\vspace{0.5ex}
		\noindent\textbf{$\mathrm{(iii)}$~\textit{Almost sure convergence}.}
		If $\sum_{k = 0}^\infty \delta_k < +\infty$, then $\Sc_K \leq \Sc_{\infty} = \sum_{k = 0}^\infty \hat{\delta}_k = \frac{128\mu^2 M_1 \norms{1 - \gamma}_{\infty}^2}{\kappa} \sum_{k = 0}^\infty \delta_k < +\infty$ and $\Rc_K^2 \leq \Rc_{\infty}^2 := \big(\Rc_0^2  + \frac{\sigma_{\star}^2}{32\mu^2}\big)\exp(\Sc_{\infty}) - \frac{\sigma_{\star}^2}{32\mu^2} < +\infty$.
		Hence, applying Lemma~\ref{lem:ch2:supermartingale} to \eqref{eq:VCOG_theorem_proof2}, we conclude that $\lim_{k \to \infty} \Ec_k$ exists almost surely. Moreover, almost surely, we have
		\myeq{eq:VCOG_summability_proof3}{
			\arraycolsep=0.2em
			\begin{array}{lclcl}
				\frac{\tau^2}{\bar{\p}} \sum_{k=0}^\infty \norms{Gx^k + \xi^k}^2 &<& +\infty
				\quad \text{and} \quad
				\frac{1}{4\mu^2} \sum_{k=0}^\infty \norms{x^{k} - x^{k-1}}^2 &<& + \infty.
			\end{array}
		}
		This implies \eqref{eq:VCOG_convergence_3}.
		Consequently,
		\myeq{eq:VCOG_summability_proof3.1}{
			\arraycolsep=0.2em
			\begin{array}{lclcl}
				\lim_{k \to \infty} \norms{Gx^k + \xi^k}^2 = 0 \quad \text{and} \quad \lim_{k \to \infty} \norms{x^{k} - x^{k-1}}^2 = 0 \quad \text{almost surely}.
			\end{array}
		}
		By Young's inequality and \eqref{eq:VCOG_Lipschitz_continuity}, we have
		\myeqn{
			\arraycolsep=0.2em
			\begin{array}{lcl}
				\gamma^2 \eta^2 \circ \norms{Gx^{k-1} + \xi^k}^2 &\leq& 2\gamma^2\eta^2 \circ \norms{Gx^k + \xi^k}^2 + 2\gamma^2\eta^2 \circ \norms{Gx^k - Gx^{k-1}}^2 \vspace{1ex}\\
				&\leq& 2\gamma^2\eta^2 \circ \norms{Gx^k + \xi^k}^2 + 2M_2^2 \sum_{j = 1}^m \frac{1}{L_j^2} \norms{[Gx^k]_j - [Gx^{k-1}]_j}^2 \vspace{1ex}\\
				&\leq& 2\gamma^2\eta^2 \circ \norms{Gx^k + \xi^k}^2 + 2M_2^2 \norms{x^k - x^{k-1}}^2,
			\end{array}
		}
		where $M_2 := \max_{j \in [m]} \sets{\gamma_j \eta_j L_j}$.
		Combining this relation and \eqref{eq:VCOG_summability_proof3.1}, we get
		\myeq{eq:VCOG_summability_proof3.2}{
			\arraycolsep=0.2em
			\begin{array}{lclcl}
				\lim_{k \to \infty} \gamma^2 \eta^2 \circ \norms{Gx^{k-1} + \xi^k}^2 = 0 \quad  \text{almost surely}.
			\end{array}
		}
		Moreover, using \eqref{eq:VCOG_Lipschitz_continuity} from Assumption~\ref{as:VCOG}(ii), we have
		\myeqn{
			\arraycolsep=0.2em
			\begin{array}{lcl}
				U_k & := & \frac{1}{n} \sum\limits_{i=1}^n \norms{G_ix^k - G_ix^{k-1}}_{\beta}^2 = \frac{1}{n} \sum\limits_{i=1}^n \sum\limits_{j=1}^m \frac{1}{ L_j^2 } \norms{[G_i x^k]_j - [G_i x^{k-1}]_j}^2
				\overset{\tiny\eqref{eq:VCOG_Lipschitz_continuity} }{ \leq } \norms{x^k - x^{k-1}}^2.
			\end{array}
		}
		Combining the second bound in \eqref{eq:VCOG_summability_proof3} and the last inequality, we obtain that
		\myeq{eq:VCOG_theorem_proof4a}{
			\arraycolsep=0.2em
			\begin{array}{lcl}
				\sum_{k=0}^\infty U_k < + \infty \quad \text{and} \quad \lim_{k \to \infty} U_k = 0 \quad \text{ almost surely.}
			\end{array}
		}
		Next, since $\Ec_k$ converges almost surely, we can show from \eqref{eq:VCOG_key_lemma2_proof4} that $\sup_k V_k < +\infty$ almost surely, which implies $\sum_{k = 1}^\infty \delta_k V_{k-1} < +\infty$ almost surely.
		Then, taking the conditional expectation $\mathbb{E}_k[\cdot]$ of \eqref{eq:VR_estimator2}, we obtain
		\myeqn{
			\arraycolsep=0.2em
			\begin{array}{lcl}
				\Expsk{k}{\Delta_k} & \leq & (1 - \kappa) \Delta_{k-1} + \Theta U_k + \hat{\Theta}  U_{k-1} + \norms{1 - \gamma}_{\infty}^2 \delta_k V_{k-1}.
			\end{array}
		}
		Since $\sum_{k=0}^\infty U_k < + \infty$ and $\sum_{k = 1}^\infty \delta_k V_{k-1} < +\infty$ almost surely, applying Lemma~\ref{lem:ch2:supermartingale}, we conclude that $\sets{\Delta_k}$ converges almost surely and $\sum_{k=0}^\infty \Delta_k < +\infty$ almost surely, which implies
		\myeq{eq:VCOG_theorem_proof4}{
			\arraycolsep=0.2em
			\begin{array}{lcl}
				\lim_{k \to \infty} \Delta_k = 0 \quad \textrm{almost surely}.
			\end{array}
		}
		Thus, from the definition of $\Ec_k$ in \eqref{eq:VCOG_potential_func}, the almost sure convergence of $\Ec_k$, \eqref{eq:VCOG_summability_proof3.1}, \eqref{eq:VCOG_theorem_proof4a}, and \eqref{eq:VCOG_theorem_proof4}, we conclude that $\lim_{k\to \infty} \norms{x^k - \xopt + \gamma \eta \circ (Gx^{k-1} + \xi^k)}^2$ exists almost surely.
		
		Next, by \eqref{eq:VCOG_key_lemma2_proof3.1} in the proof of Lemma~\ref{lem:ch2:VCOG_key_lemma2} and the almost-sure convergence of $\Ec_k$, there exists a random variable $M$, finite almost surely, such that $\norms{x^k - \xopt} \leq M$ for all $k \geq 0$.
		By \eqref{eq:VCOG_summability_proof3.2}, we have $\norms{\gamma\eta \circ (Gx^{k-1} + \xi^k)} \to 0$ almost surely.
		Hence, almost surely,
		\myeqn{
			\arraycolsep=0.2em
			\begin{array}{lcl}
				\big\vert \norms{x^k - \xopt + \gamma\eta \circ (Gx^{k-1} + \xi^k)}^2 - \norms{x^k - \xopt}^2 \big\vert &\leq& 2M\norms{\gamma\eta \circ (Gx^{k-1} + \xi^k)} \vspace{1ex}\\
				&& + {~} \norms{\gamma\eta \circ (Gx^{k-1} + \xi^k)}^2 \longrightarrow 0.
			\end{array}
		}
		Since $\lim_{k \to \infty} \norms{x^k - \xopt + \gamma\eta \circ (Gx^{k-1} + \xi^k)}^2$ exists almost surely, so does $\lim_{k \to \infty} \norms{x^k - \xopt}^2$.
		As $\xopt \in \zerM{\Psi}$ is arbitrary, for every $\xopt \in \zerM{\Psi}$, the limit $\lim_{k \to \infty} \norms{x^k - \xopt}$ exists almost surely.
		Under the additional assumption $\zerM{\Psi} = \zer{\Psi}$, this holds for every $\xopt \in \zer{\Psi}$.
		Finally, since $\lim_{k \to \infty} \norms{Gx^k + \xi^k} = 0$ almost surely by \eqref{eq:VCOG_convergence_3} and $\gra{\Psi}$ is closed, Lemma~\ref{lem:ch2:Opial} implies that $\sets{x^k}$ converges almost surely to a $\zer{\Psi}$-valued random variable $\bar{x}^{\star}$, which is a solution of \eqref{eq:FNI}.
		\end{proof}
	
	{\begin{remark}[Mini-batch schedules]
			For the mini-batch estimator, the condition
			$\sum_{k=0}^\infty\delta_k<\infty$ can be ensured by choosing $b_k := \left\lceil\frac{n}{1+(n-1)\vartheta_k}\right\rceil$.
			Some possibilities are $\vartheta_k = \frac{1}{K+1}\mbf{1}_{\set{k \leq K}}$ (fixed budget), $\vartheta_k = \frac{c}{(k+2)^{1+\zeta}}$ (polynomial), and $\vartheta_k = \frac{c}{(k+2)[\log(k+2)]^{1+\zeta}}$ (nearly linear), where $c, \zeta > 0$ and $K$ is a prescribed iteration budget.
			Then, $\delta_k \leq 2\vartheta_k$ in all three cases.
			The first choice uses a constant batch size up to iteration $K$ and	a full batch thereafter, thus if $K \gg n$, it forces the use of almost full batch at the very beginning.
			The other two options provide budget-independent increasing batch sizes, which can be controlled by the constant parameters $c$ and $\zeta$, chosen appropriately.
		\end{remark}
	}
	
	\beforesubsec
	\subsection{Oracle complexity with importance-based sampling strategy}\label{subsec:ch2:oracle_complexity_of_VCOG}
	\aftersubsec
	{In this subsection, let $\p_j := \frac{L_j^{2/3}}{\sum_{l=1}^m L_l^{2/3}} = \frac{L_j^{2/3}}{\norms{L^{2/3}}_1}$, leading to $\frac{L_j}{\p_j^{3/2}} = \norms{L^{2/3}}_1^{3/2}$.
	Therefore, 
	\myeqn{
		\arraycolsep=0.2em
		\begin{array}{lcl}
			\bar{\p} &=& \frac{\norms{L^{2/3}}_{\infty}}{\norms{L^{2/3}}_1}, \qquad
			L_{\p} =  \norms{L^{2/3}}_1^{3/2}, \quad \text{and} \quad 
			\rho_{\p} =  \frac{1}{\norms{L^{2/3}}_1} \max_{j \in [m]} L_j^{2/3} \rho_j,
		\end{array}
	}
	and hence $L_{\p} \rho_{\p} = \norms{L^{2/3}}_1^{1/2} \max_{j \in [m]} L_j^{2/3} \rho_j$.
	
	To output a single point, let $I$ be drawn uniformly from $\sets{0, \cdots, K}$, independently of the algorithmic randomness, and set $(\hat{x}^K, \hat{\xi}^K) := (x^I, \xi^I)$.
	Then $\hat{\xi}^K \in T\hat{x}^K$, and $\Expk{\norms{G\hat{x}^K + \hat{\xi}^K}^2}$ equals the left-hand side of \eqref{eq:VCOG_convergence_2}.
	We now specify the oracle complexity of Algorithm~\ref{alg:A1} for the mini-batch, L-SVRG, and GJS estimators.
	For the mini-batch estimator, we fix $c > 0$ and $\zeta > 0$ and use the batch sizes
	\myeq{eq:VCOG_batch_schedule}{
		\arraycolsep=0.2em
		\begin{array}{lcl}
			b_k &:=& \Big\lceil \frac{n}{1 + c(n-1)/\big((k+2)[\log(k+2)]^{1+\zeta}\big)} \Big\rceil, \qquad k \geq 1,
		\end{array}
	}
	while the initial estimator is computed exactly, as in Algorithm~\ref{alg:A1}.
	
	\begin{corollary}\label{co:VCOG_complexity}
		Suppose that Assumptions~\ref{as:VCOG} and \ref{as:VCOG_T_separable}
		hold for \eqref{eq:FNI}.
		Assume that $L_{\p} \rho_{\p} < \frac{\omega}{16\sqrt{2}\mu}$.
		Let $\sets{x^k}$ be generated by Algorithm~\ref{alg:A1} using
		$\tau := \frac{\omega}{8\sqrt{2}\mu L_{\p}}$.
		Let $\Rc_K$ be defined in \eqref{eq:VCOG_convergence_2}.
		Then, the following statements hold:
		\begin{compactitem}
		\item[\textbf{$\mathrm{(i)}$~$(\text{Mini-batch estimator})$}]
		If the estimator \eqref{eq:minibatch_estimator} is used, then we have
		\myeq{eq:VCOG_complexity_convergence1}{
				\frac{1}{K+1}\sum_{k=0}^K \Expk{\norms{Gx^k + \xi^k}^2} \leq \frac{256\mu^2 \norms{L^{2/3}}_{\infty} \norms{L^{2/3}}_1^2}{K+1} \Rc_K^2.
		}
		For $\epsilon > 0$, if we use the batch sizes \eqref{eq:VCOG_batch_schedule}, then the number of evaluations of $G_i$ and $J_{\gamma \eta \circ T}$ that Algorithm~\ref{alg:A1} requires to achieve $\Expk{\norms{G\hat{x}^K + \hat{\xi}^K}^2} \leq \epsilon^2$
		 are at most $\Tc_{G_i} := n + \BigOlog{\min\sets{\frac{\mu^4 \norms{L^{2/3}}_{\infty}^2 \norms{L^{2/3}}_1^4}{m \epsilon^4}, \frac{n \mu^2 \norms{L^{2/3}}_{\infty} \norms{L^{2/3}}_1^2}{m \epsilon^2}}}$  and $\Tc_{T} := \BigO{\frac{\mu^2\norms{L^{2/3}}_{\infty} \norms{L^{2/3}}_1^2}{m \epsilon^2 }}$, respectively.
		
		\item[\textbf{$\mathrm{(ii)}$ $($L-SVRG and Generalized JacSketch$)$}]
		If \eqref{eq:loopless_SVRG} or \eqref{eq:JacSketch_estimator} is used and $b\q^2 \leq 1$, then we have
		\myeq{eq:VCOG_complexity_convergence2}{
				\frac{1}{K+1}\sum_{k=0}^K \Expk{\norms{Gx^k + \xi^k}^2} \leq \frac{2048\mu^2 \norms{L^{2/3}}_{\infty} \norms{L^{2/3}}_1^2 \Rc_0^2}{b\q^2(K+1)}.
		}
		For $\epsilon > 0$, if we set $\q := n^{-1/3}$ and $b := \lfloor n^{2/3} \rfloor$, then Algorithm~\ref{alg:A1} requires at most $\Tc_{G_i} := n + \BigO{\mu^2 \norms{L^{2/3}}_{\infty} \norms{L^{2/3}}_1^2 n^{2/3}\epsilon^{-2}}$ evaluations of $G_i$ and $\Tc_{T} := \BigO{\frac{\mu^2\norms{L^{2/3}}_{\infty} \norms{L^{2/3}}_1^2}{m \epsilon^2 }}$ evaluations of $J_{\gamma \eta \circ T}$  to achieve $\Expk{\norms{G\hat{x}^K + \hat{\xi}^K}^2} \leq \epsilon^2$. 
		\end{compactitem}
	\end{corollary}
	
		\begin{proof}
			First, if we choose $\lambda := \frac{1}{4}$, then the second condition in \eqref{eq:VCOG_param_conditions} becomes $L_{\p} \rho_{\p} < \frac{\omega}{16\sqrt{2}\mu}$.
			Next, we choose $\tau := \frac{\omega}{8\sqrt{2}\mu L_{\p}}$, then we get  $\frac{\bar{\p}}{\tau^2} = \frac{128\mu^2}{\omega^2}\norms{L^{2/3}}_{\infty} \norms{L^{2/3}}_1^2$.
			Using this expression, we obtain from \eqref{eq:VCOG_convergence_2} that
			\myeq{eq:VCOG_complexity_convergence}{
				\arraycolsep=0.2em
				\begin{array}{lcl}
					\frac{1}{K+1}\sum_{k=0}^K \Expk{\norms{Gx^k + \xi^k}^2} &\leq& \frac{128\mu^2 \norms{L^{2/3}}_{\infty} \norms{L^{2/3}}_1^2}{\omega^2 (K+1)} \Rc_K^2.
				\end{array}
			}
			Now, we consider the three estimators.
			
			\vspace{0.5ex}
			\noindent\textbf{$\mathrm{(i)}$~Mini-batch estimator.} For \eqref{eq:minibatch_estimator}, we have $\kappa = 1$, $\Theta = 2$, $\hat{\Theta} = 0$, and $\delta_k = \frac{2(n - b_k)}{b_k(n-1)}$, leading to $\omega^2 = \frac{1}{2}$.
			Substituting this into \eqref{eq:VCOG_complexity_convergence}, we obtain \eqref{eq:VCOG_complexity_convergence1}.
			
			Hence, to guarantee $\frac{1}{K+1}\sum_{k=0}^K \Expk{\norms{Gx^k + \xi^k}^2} \leq \epsilon^2$, from  \eqref{eq:VCOG_convergence_2} we can impose $\frac{256\mu^2 \norms{L^{2/3}}_{\infty} \norms{L^{2/3}}_1^2}{K+1} \Rc_{K}^2 \leq \epsilon^2$.
			To keep $\Rc_K = \BigOs{1}$, we use \eqref{eq:VCOG_batch_schedule}, which gives $b_k = \BigOs{k(\log k)^{1+\zeta}}$.
			Then, $\delta_k \leq \frac{2c}{(k+2)[\log(k+2)]^{1+\zeta}}$ so that $\sum_{k = 0}^{\infty} \delta_k < +\infty$ and hence $\sup_{K \geq 0} \Rc_K < +\infty$.
			Therefore, it suffices to take $K := \Big\lfloor \frac{256 \mu^2 \norms{L^{2/3}}_{\infty} \norms{L^{2/3}}_1^2\Rc_{\infty}^2}{\epsilon^2} \Big\rfloor$.
			Initially, we compute $Gx^0$, which requires $n$ evaluations of $G_i$.
			At each iteration, Algorithm~\ref{alg:A1} with \eqref{eq:minibatch_estimator} evaluates two mini-batches on one randomly selected block.
			Then, the total complexity is
			\myeqn{
				\arraycolsep=0.2em
				\begin{array}{lcl}
					\Tc_{G_i} &:=& n + \frac{1}{m} \sum_{k=1}^K 2b_k
					= n + \BigOs{\min\sets{\frac{K^2 [\log(K + 2)]^{1+\zeta}}{m}, \frac{nK}{m}}} \vspace{1ex}\\
					&=& n + \BigOlog{\min\sets{\frac{\mu^4 \norms{L^{2/3}}_{\infty}^2 \norms{L^{2/3}}_1^4}{m \epsilon^4}, \frac{n \mu^2 \norms{L^{2/3}}_{\infty} \norms{L^{2/3}}_1^2}{m \epsilon^2}}}.
				\end{array}
			}
			If additionally Assumption~\ref{as:VCOG_T_separable} holds, only one resolvent block is evaluated per iteration.
			Therefore, the total evaluation of $J_{\gamma \eta \circ T}$ is $\Tc_{T} := \frac{K}{m} = \BigO{\frac{\mu^2\norms{L^{2/3}}_{\infty} \norms{L^{2/3}}_1^2}{m \epsilon^2 }}$.
			
			{\vspace{0.5ex}
				\noindent\textbf{$\mathrm{(ii)}$~L-SVRG estimator.}
				For \eqref{eq:loopless_SVRG}, we have $\kappa := \frac{\q}{2}$, $\Theta = \hat{\Theta} := \frac{4}{b\q}$, and $\delta_k = 0$, leading to $\omega^2 = \frac{b\q^2}{16}$.
				In \eqref{eq:VCOG_convergence_2}, since $\delta_k = 0$, we have $\Rc_K^2 = \Rc_0^2$.
				Substituting this into \eqref{eq:VCOG_complexity_convergence}, we obtain \eqref{eq:VCOG_complexity_convergence2}.
				
				Hence, to guarantee $\frac{1}{K+1}\sum_{k=0}^K \Expk{\norms{Gx^k + \xi^k}^2} \leq \epsilon^2$, from  \eqref{eq:VCOG_convergence_2}, we need to impose $\frac{2048\mu^2 \norms{L^{2/3}}_{\infty} \norms{L^{2/3}}_1^2 \Rc_0^2}{b\q^2 (K+1)} \leq \epsilon^2$, which holds if $K := \Big\lfloor \frac{2048\mu^2 \norms{L^{2/3}}_{\infty} \norms{L^{2/3}}_1^2 \Rc_0^2}{b\q^2 \epsilon^2} \Big\rfloor$.
				At each iteration $k$, Algorithm~\ref{alg:A1} with \eqref{eq:loopless_SVRG} evaluates $3$ mini-batches of size $b$ on one block and refreshes the full $Gw^k$ with probability $\q$.
				Since Algorithm~\ref{alg:A1} stores $x^k$ explicitly, these refreshes do not require any additional evaluation of $J_{\gamma\eta\circ T}$.
				Thus, the expected component-oracle cost per iteration is $n\q + \frac{3b}{m}$, and therefore, the total complexity is
				\myeqn{
					\arraycolsep=0.2em
					\begin{array}{lcl}
						\Tc_c &:=& K \big(n\q + \frac{3b}{m}\big) = \frac{2048\mu^2 \norms{L^{2/3}}_{\infty} \norms{L^{2/3}}_1^2 \Rc_0^2}{b\q^2 \epsilon^2} \big(n\q + \frac{3b}{m}\big).
					\end{array}
				}
				If we choose $b := \lfloor n^{2/3} \rfloor$ and $\q := n^{-1/3}$, then  $\Tc_c = \BigOs{\frac{\mu^2 \norms{L^{2/3}}_{\infty} \norms{L^{2/3}}_1^2 \Rc_0^2 n^{2/3}}{\epsilon^2}}$.
				Accounting for the initial full evaluation of $Gx^0$, the total evaluations of $G_i$ is $\Tc_{G_i} = n + \Tc_c = n + \BigO{\frac{\mu^2 \norms{L^{2/3}}_{\infty} \norms{L^{2/3}}_1^2 \Rc_0^2 n^{2/3}}{\epsilon^2}}$.
				If additionally Assumption~\ref{as:VCOG_T_separable} holds, then the total evaluation of $J_{\gamma \eta \circ T}$ is again $\Tc_{T} := \BigO{\frac{\mu^2\norms{L^{2/3}}_{\infty} \norms{L^{2/3}}_1^2}{m \epsilon^2 }}$.
				
				\vspace{0.5ex}
				\noindent\textbf{$\mathrm{(iii)}$~Generalized JacSketch-type estimator.}
				At each iteration, Algorithm~\ref{alg:A1} with \eqref{eq:JacSketch_estimator} evaluates $2$ mini-batches of size $b$ on one block and refreshes the full table $\sets{H_i^k}_{i=1}^n$ with probability $\q$.
				Thus, as for \eqref{eq:loopless_SVRG}, the total complexity is
				\myeqn{
					\arraycolsep=0.2em
					\begin{array}{lcl}
						\Tc_c &:=& K \left(2n\q + \frac{2b}{m}\right) = \frac{4096\mu^2 \norms{L^{2/3}}_{\infty} \norms{L^{2/3}}_1^2 \Rc_0^2}{b\q^2 \epsilon^2} \big(n\q + \frac{b}{m}\big).
					\end{array}
				}
				If we choose $b := \lfloor n^{2/3} \rfloor$ and $\q := n^{-1/3}$, then we obtain the same asymptotic oracle complexities as for \eqref{eq:loopless_SVRG}.
			}
		\end{proof}
	}
	
		\begin{remark}[\textbf{Effect of block heterogeneity}]
			The oracle complexities depend on the block smoothness through $\norms{L^{2/3}}_{\infty}\norms{L^{2/3}}_1^2$.
			If the blocks are highly heterogeneous so that	$\norms{L^{2/3}}_1=\BigOs{L_{\max}^{2/3}}$,
			where $L_{\max}:=\max_j L_j$, then
			$\norms{L^{2/3}}_{\infty}\norms{L^{2/3}}_1^2
			=\BigOs{L_{\max}^2}$.
			Consequently, for \ref{eq:loopless_SVRG} and
			\ref{eq:JacSketch_estimator}, importance sampling removes any explicit dependence on $m$, yielding $n + \BigOs{n^{2/3} \epsilon^{-2}}$ oracle complexity up to
			problem-dependent constants.
			In contrast, if $L_j = \bar L$ for all $j \in [m]$, then
			$\p_j=1/m$ and 	$\norms{L^{2/3}}_{\infty} \norms{L^{2/3}}_1^2 = m^2\bar L^2$,	leading to	$n + \BigOs{m^2 n^{2/3} \epsilon^{-2}}$.
			Thus, importance sampling is most effective when the block Lipschitz constants are highly nonuniform, consistently with
			\cite{diakonikolas2025block}.
			Moreover, in the uniform case $L_j = \bar{L}$, \eqref{eq:VCOG_complexity_convergence} and Jensen's inequality give
			$\min_{0\le k\le K}	\Expk{\norms{Gx^k+\xi^k}} \leq
			\frac{8\sqrt{2}\,m\mu\bar L\Rc_K}{\omega\sqrt{K+1}}$,
			so the bound grows linearly in $m$, as also observed for coordinate methods in optimization and root-finding; see, e.g.,	\cite{Necoara2016,Nesterov2012,richtarik2016parallel,wright2015coordinate,tran2022accelerated}.
	\end{remark}

	\beforesec
	\section{\textls[-0]{Variance-Reduced Block-Coordinate Fast OG Method for Solving \eqref{eq:FNI}}}\label{sec:ch2:VACOG}
	\aftersec
	Algorithm~\ref{alg:A1} attains an $\BigO{1/k}$ rate under relatively mild assumptions.
	Our next goal is to develop a new \textit{\textbf{V}ariance-reduced Block-\textbf{C}oordinate \textbf{F}ast  \textbf{O}ptimistic \textbf{G}radient method} (VCFOG) for \eqref{eq:FNI}, which achieves an $\BigO{1/k^2}$ rate.
	However, this requires narrowing the problem class by imposing the following co-coercive-type assumption on $G$.
	\begin{assumption}\label{as:VACOG}
		The composite inclusion \eqref{eq:FNI} satisfies the following conditions:
		\begin{enumerate}[label=(\roman*), leftmargin=*, nosep, widest=ii, align=left]
			\itemsep=-0.5em
			\item[$\mathrm{(i)}$] (\textbf{\textit{Co-coercivity}})
			$G$ is $L^{-1}$-average co-coercive, i.e., there exists a vector $L \in \R^{mn}_{++}$ such that for any $x, y \in \dom{G}$:
			\myeq{eq:FNI_assumption_cocoercive}{
					\iprods{Gx - Gy, x - y}  \geq \frac{1}{n} \sum_{i=1}^n \sum_{j=1}^m \frac{1}{L_{ij}}  \norms{[G_i x]_j - [G_i y]_j}^2.
			}
			\item[$\mathrm{(ii)}$] (\textbf{\textit{Co-hypomonotonicity}})
			$T$ is block-separable, i.e.,
			$Tx = \big[ T_1x_1, \cdots, T_mx_m \big]$, and maximally $\rho$-co-hypomonotone with $\rho \in \R^m_{+}$
			such that $\norms{\bar{L}\rho}_{\infty} < 1$, where $\bar{L} := (\bar{L}_1, \cdots, \bar{L}_m)$ with $\bar{L}_j := \max_{i \in [n]} L_{ij}$ for $j \in [m]$.
		\end{enumerate}
	\end{assumption}
	
	Assumption~\ref{as:VACOG}(i) is stronger than the Lipschitz continuity in Assumption~\ref{as:VCOG}(ii) but is standard in accelerated methods \cite{attouch2020convergence,bot2022bfast,mainge2021accelerated,tran2024accelerated,TranDinh2025a,tran2022accelerated},
	while Assumption~\ref{as:VACOG}(ii) still allows nonmonotonicity.
	Throughout this section, we instantiate Definition~\ref{def:VR_estimator} with $\beta_j := \p_j^{-1}$ for all $j \in [m]$, where $\p_j$ is defined in \eqref{eq:VACOG_prob}.
	
	\beforesubsec
	\subsection{Equivalent reformulation of \eqref{eq:FNI}}\label{subsec:ch2:FNI_reform}
	\aftersubsec
	For a given $\lambda \in \R^m_{++}$, define the block-wise resolvent $J_{\lambda \circ T}(u) := [J_{\lambda_1 T_1}(u_1), \dots, J_{\lambda_m T_m}(u_m)]$.
	Let $x := J_{\lambda \circ T}(u)$ and $A_{\lambda \circ T}u := \lambda^{-1} \circ (u - x)$.
	Then, the backward-forward splitting (BFS) operator associated with \eqref{eq:FNI} is defined as follows:
	\myeq{eq:BFS_operators}{
		\arraycolsep=0.2em
		\begin{array}{lclcl}
			\Gc_{i, \lambda}u &:=& G_ix + A_{\lambda \circ T}u \quad \text{and} \quad
			\Gc_{\lambda}u &:=& \frac{1}{n}\sum\limits_{i=1}^n \Gc_{i, \lambda}u = Gx + A_{\lambda \circ T}u.
		\end{array}
	}
	The following key properties of the BFS operator are required to develop our algorithm, whose proof can be found in Appendix~\ref{apdx:subsec:ch2:VACOG_technical_lemmas}.
	
	\begin{lemma}\label{lem:ch2:BFS_FNI_equivalence}
		For \eqref{eq:FNI}, suppose that Assumption~\ref{as:VACOG} holds and that $\lambda_j > \rho_j$ for all $j \in [m]$, so that each block resolvent $J_{\lambda_j T_j}$ is single-valued and has full domain.
		Then, the following statements hold.
		\begin{compactenum}
			\item[$\mathrm{(i)}$] $\zer{G + T} = J_{\lambda \circ T}(\zer{\Gc_\lambda})$.
			Moreover, for $x := J_{\lambda \circ T}u$ and $\xi := \lambda^{-1} \circ (u - x) \in Tx$, we have $\norms{\Gc_{\lambda}u} = \norms{Gx + \xi}$.
			\item[$\mathrm{(ii)}$]  If $\lambda_j \geq 2\rho_j$ for all $j \in [m]$, then $J_{\lambda \circ T}$ is nonexpansive.
			\item[$\mathrm{(iii)}$]  If $\lambda_j \in \Big( \frac{2 - 2\sqrt{1 - \bar{L}_j\rho_j}}{\bar{L}_j}, \frac{2 + 2\sqrt{1 - \bar{L}_j\rho_j}}{\bar{L}_j}\Big)$,
			then $\Gc_{\lambda}$ is $\bar{\tau}$-average co-coercive, i.e.,
			\myeq{eq:BFS_cocoercive}{
				\arraycolsep=0.2em
				\begin{array}{lcl}
					\iprods{\Gc_{\lambda}u - \Gc_{\lambda}v, u - v} \geq \frac{1}{n} \sum\limits_{i=1}^n \sum\limits_{j=1}^m \bar{\tau}_{ij} \norms{[\Gc_{i, \lambda} u]_j - [\Gc_{i, \lambda} v]_j}^2, \quad  \forall u, v \in \R^p,
				\end{array} \vspace{-2ex}
			}
			where $\bar{\tau}_{ij} := \frac{4(\lambda_j - \rho_j) - \lambda_j^2 L_{ij}}{4(1 - \rho_j L_{ij})}$.
			Consequently, $\Gc_{\lambda}$ is continuous.
		\end{compactenum}
	\end{lemma}
	Now, let $\Lambda_0 := \min\limits_{i,j} \p_j \bar{\tau}_{ij}$. Then, \eqref{eq:BFS_cocoercive} implies
	\myeqn{
		\arraycolsep=0.2em
		\begin{array}{lcl}
			\iprods{\Gc_{\lambda}u - \Gc_{\lambda}v, u - v} &\geq& \frac{\Lambda_0}{n} \sum\limits_{i=1}^n \norms{\Gc_{i, \lambda} u - \Gc_{i, \lambda} v}_{\beta}^2, \quad  \forall u, v \in \R^p.
		\end{array}
	}
	On the other hand, using Jensen's inequality, we have
	\myeqn{
		\arraycolsep=0.2em
		\begin{array}{lcl}
			\norms{\Gc_{\lambda}u - \Gc_{\lambda}v}^2 &\leq& \frac{1}{n}\sum\limits_{i=1}^n \sum\limits_{j=1}^m \p_j \cdot \frac{1}{\p_j}\norms{[\Gc_{i, \lambda}u]_j - [\Gc_{i, \lambda}v]_j}^2
			\leq \frac{\bar{\p}}{n}\sum\limits_{i = 1}^n \norms{\Gc_{i, \lambda}u - \Gc_{i, \lambda}v}_{\beta}^2.
		\end{array}
	}
	Combining the last two relations and using the Cauchy--Schwarz inequality, we conclude that $\Gc_{\lambda}$ is $\frac{\bar{\p}}{\Lambda_0^2}$-average Lipschitz continuous, i.e., $\frac{1}{n} \sum\limits_{i=1}^n \norms{\Gc_{i, \lambda} u - \Gc_{i, \lambda} v}_{\beta}^2 \leq \frac{\bar{\p}}{\Lambda_0^2} \norms{u-v}^2.$
	If we define $\Vc_{\lambda}(u) := \frac{1}{n} \sum_{i = 1}^n \norms{\Gc_{i,\lambda}u - \Gc_{\lambda} u}_{\beta}^2 \geq 0$ and $\sigma_{\lambda, \star}^2 := \Vc_{\lambda}(u^{\star}) < +\infty$, then, as for \eqref{eq:VCOG_Vk_bound} in Section~\ref{sec:ch2:VCOG_alg}, we obtain
	\myeq{eq:VACOG_Vk_bound}{
		\arraycolsep=0.2em
		\begin{array}{lcl}
			V_k = \Vc_{\lambda}(u^k) \leq \frac{2\bar{\p}}{\Lambda_0^2}\norms{u^k - u^{\star}}^2 + 2\sigma_{\lambda, \star}^2.
		\end{array}
	}
	
	\beforesubsec
	\subsection{The proposed method and its implementable version}\label{subsec:ch2:main_alg}
	\aftersubsec
	We now develop the accelerated framework for the BFS equation
	$\Gc_\lambda u^{\star} = 0$.
	The main idea is to combine a fast optimistic-gradient step with randomized block-coordinate updates and variance reduction, while recovering a solution of \eqref{eq:FNI} through
	$x^k := J_{\lambda\circ T}u^k$.
	More precisely, we apply the accelerated OG scheme to the direction
	$\Sc_{\lambda}^k := \Gc_\lambda u^k - \gamma_k \Gc_\lambda u^{k-1}$, but update only the randomly selected block $j_k$ and replace the corresponding block of $\Sc_{\lambda}^k$ by a variance-reduced stochastic estimator.
	The remaining question is how to construct such an estimator without working directly with the BFS components $\Gc_{i,\lambda}$.
	For this purpose, using \eqref{eq:BFS_operators} and $x^k=J_{\lambda\circ T}u^k$, we obtain
	\myeqn{
		\arraycolsep=0.2em
		\begin{array}{lcl}
			\Sc_{i, \lambda}^k
			& = &  S_i^k + \lambda^{-1} \circ \big( u^k - x^k - \gamma_k(u^{k-1} - x^{k-1}) \big), \vspace{1ex}\\
			\Sc_{\lambda}^k
			& = &  S^k + \lambda^{-1} \circ \big( u^k - x^k - \gamma_k(u^{k-1} - x^{k-1}) \big),
		\end{array}
	}
	where $S_i^k := G_ix^k - \gamma_k G_ix^{k-1}$ and $S^k := G{x}^k - \gamma_k G{x}^{k-1}$.
	Suppose that $\widetilde{S}^k$ is a variance-reduced estimator of $S^k$ generated by a control-variate construction rule $\mathfrak{C}$ satisfying the condition of Lemma~\ref{le:control_variate_compability}.
	Then, by Lemma~\ref{le:control_variate_compability}, we have
	\myeq{eq:BFS_estimator_equivariance}{
		\arraycolsep=0.2em
		\begin{array}{lcl}
			\widetilde{\Sc}_{\lambda}^k
			& := &  \widetilde{S}^k + \lambda^{-1} \circ \big( u^k - x^k - \gamma_k(u^{k-1} - x^{k-1}) \big), 
		\end{array}
	}
	is a variance-reduced estimator of $\Sc_{\lambda}^k$ satisfying Definition~\ref{def:VR_estimator} (with $U_k$ and $V_k$ defined in terms of $\Gc_{i, \lambda}$).
	Therefore, we can construct $\widetilde{S}^k$ first and obtain $\widetilde{\Sc}_{\lambda}^k$ from \eqref{eq:BFS_estimator_equivariance}, avoiding a direct construction based on the BFS components.
	
	Putting the above steps together, we obtain Algorithm~\ref{alg:A2}.
	Note that when $\norms{\bar{L}\rho}_{\infty} < 1$, we can verify that the range of $\lambda$ in Algorithm~\ref{alg:A2} is nonempty.
	
	\setlength{\textfloatsep}{10pt plus 1.0pt minus 3.0pt}
	\begin{algorithm}[hpt!]\caption{(\textit{\textbf{The VCFOG framework for solving \eqref{eq:FNI} under Assumption~\ref{as:VACOG}}})}\label{alg:A2}
		\normalsize
		\begin{algorithmic}[1]
			\STATE\label{step:A2_i0}{\bfseries Initialization:}
			Choose $u^{-1} = u^0 \in \R^p$,  $\lambda \in \R^m_{++}$ where $\lambda_j \in \Big[2\rho_j, \frac{2+2\sqrt{1-\bar{L}_j\rho_j}}{\bar{L}_j}\Big)$.
			Compute $x^{-1} = x^0 := J_{\lambda \circ T}u^0$ and $\widetilde{S}^0 := (1-\gamma_0)Gx^0$.
			Choose $k_{\max} > 0$.
			\STATE\hspace{0ex}\label{step:A2_loop}{\bfseries For $k := 0,\cdots, k_{\max}$ do}
			\vspace{0.25ex}
			\STATE\hspace{1.5ex}\label{step:A2_o0}Update $\gamma_k$, $\theta_k$, and $\eta_k$ according to \eqref{eq:th1_param}.
			\STATE\hspace{1.5ex}\label{step:A2_o1}\textls[0]{Sample a block $j_k \in [m]$ following the rule \eqref{eq:VACOG_prob}.}
			\STATE\hspace{1.5ex}\label{step:A2_o2}If $k\geq 1$, construct $\widetilde{S}^k$ of $S^k$ independent of $j_k$ and satisfying Definition~\ref{def:VR_estimator}.
			\STATE\hspace{1.5ex}\label{step:A2_o3}\textls[-4]{Compute $\widetilde{\Sc}_{\lambda}^k := \widetilde{S}^k + \lambda^{-1} \circ \big( u^k - x^k - \gamma_k(u^{k-1} - x^{k-1}) \big)$.}
			\STATE\hspace{1.5ex}\label{step:A2_o4}For each block $j \in [m]$,  update
			\myeqn{
				\arraycolsep=0.2em
				u^{k+1}_j := \vspace{-1ex}\left\{ \begin{array}{ll}
					u^k_j + \theta_k (u^k_j - u^{k-1}_j) - \frac{\eta_k}{\p_j}[\widetilde{\Sc}_{\lambda}^k]_j & \text{if } j = j_k, \vspace{1ex} \\
					u^k_j + \theta_k (u^k_j - u^{k-1}_j) & \text{otherwise}.
				\end{array}\right.
			}
			\STATE\hspace{1.5ex}\label{step:A2_o5}Compute $x^{k+1} := J_{\lambda \circ T} u^{k+1}$.
			\STATE\hspace{0ex}{\bfseries End For}
		\end{algorithmic}
	\end{algorithm}
	Algorithm~\ref{alg:A2} remains conceptual because the inertial step in Step~7 updates all blocks of $u^{k+1}$ and Step~8 requires the full evaluation of $J_{\lambda\circ T}$.
	On the other hand, only the block of the stochastic
	estimator is needed.
	In particular, by \eqref{eq:BFS_estimator_equivariance},
	$[\widetilde{\Sc}_{\lambda}^k]_{j_k} = [\widetilde S^k]_{j_k} + \frac{1}{\lambda_{j_k}} \big(u^k_{j_k} - x^k_{j_k} - \gamma_k (u^{k-1}_{j_k} - x^{k-1}_{j_k})\big)$,
	so neither $\widetilde S^k$ nor $\widetilde{\Sc}_{\lambda}^k$ needs to be formed in full.
	It remains to eliminate the full-dimensional inertial and
	resolvent updates.
	To this end, we exploit the block-separability of $T$ together
	with the implementation technique of~\cite{fercoq2015accelerated}, leading to the following equivalent block-coordinate scheme.
	
	\begin{proposition}[\textbf{\textit{Block-coordinate implementation of  Algorithm~\ref{alg:A2}}}]
		Given $u^0 \in \dom{\Gc_{\lambda}}$, set $z^0 = z^{-1} := u^0$ and $w^0 = w^{-1} := 0$, and at each iteration $k \ge 0$, perform the following steps:
		\begin{compactenum}
			\item[$\mathrm{(1)}$] sample a block index $j_k \in [m]$ $($i.i.d. across $k$$)$ according to the rule \eqref{eq:VACOG_prob}.
			\item[$\mathrm{(2)}$] construct a variance-reduced estimator $\widetilde{\Sc}_{\lambda}^k$ of $\Sc_{\lambda}^k := \Gc_{\lambda}(z^k + c_k w^k) - \gamma_k \Gc_{\lambda}(z^{k-1} + c_{k-1}w^{k-1})$, independent of $j_k$.
			\item[$\mathrm{(3)}$] update one block of $w^{k+1}$ and $z^{k+1}$ accordingly 
				\myeq{eq:VACOG_practical}{
				\arraycolsep=0.2em
				\left\{\begin{array}{lcl}
					w^{k+1}_{j_k} &:=& w^k_{j_k} -  \frac{\eta_k}{\p_{j_k}\tau_{k+1}}[\widetilde{\Sc}_{\lambda}^k]_{j_k}, \vspace{0.75ex} \\
					z^{k+1}_{j_k} &:=& z^k_{j_k} + \frac{c_k \eta_k}{\p_{j_k}\tau_{k+1}}[\widetilde{\Sc}_{\lambda}^k]_{j_k},
				\end{array}\right.
			}
			and keep $w^{k+1}_j = w^k_j$ and $z^{k+1}_j = z^k_j$ for $j \neq j_k$.
			Here, $c_0 = c_{-1} := 0$, $\tau_0 := 1$, and $\tau_{k+1} := \tau_k \theta_k$ and $c_{k+1} := c_k + \tau_{k+1}$ for $k \geq 0$, where $\theta_k$ is given by \eqref{eq:th1_param}.
		\end{compactenum}
		Then, the sequence $\sets{u^k}$ defined by $u^k := z^k + c_k w^k$ generated by \eqref{eq:VACOG_practical} is identical to the one generated by Algorithm~\ref{alg:A2}.
	\end{proposition}
	
	\begin{proof}
		We run both schemes with the same block samples and the same estimator randomness.
		Since $\theta_k > 0$, all $\tau_k$ are positive, so \eqref{eq:VACOG_practical} is well defined.
		Let $u^k := z^k + c_kw^k$.
		By \eqref{eq:VACOG_practical} and $c_{k+1} = c_k + \tau_{k+1}$, we have
		\myeqn{
		\arraycolsep=0.2em
		\begin{array}{lcl}
			u^{k+1} = u^k + \tau_{k+1}w^k - \frac{\eta_k}{\p_{j_k}}\widetilde{\Sc}^k_{\lambda,[j_k]} \qquad \text{and} \qquad
			\tau_{k+1}w^{k+1} = \tau_{k+1}w^k - \frac{\eta_k}{\p_{j_k}}\widetilde{\Sc}^k_{\lambda,[j_k]}.
		\end{array}
		}
		Hence, $u^{k+1} - u^k = \tau_{k+1}w^{k+1}$ for all $k \geq 0$ (and also for $k = -1$ because $u^{-1} = u^0$ and $w^0 = 0$).
		Consequently, $\theta_k(u^k - u^{k-1}) = \theta_k\tau_kw^k = \tau_{k+1}w^k$, and therefore $u^{k+1} = u^k + \theta_k(u^k - u^{k-1}) - \frac{\eta_k}{\p_{j_k}}\widetilde{\Sc}^k_{\lambda,[j_k]}$, which is exactly \eqref{eq:VACOG}.
		Since both schemes start from $u^{-1} = u^0$, induction on $k$ shows that they generate the same sequence $\sets{u^k}$.
	\end{proof}
	Since \eqref{eq:VACOG_practical} only requires one block $[\widetilde{\Sc}^k_\lambda]_{j_k}$, we do not need to form the entire estimator $\widetilde{\Sc}^k_\lambda$, but only its block $j_k$.
	Similarly, we do not need to compute the entire $J_{\lambda \circ T}$ when forming $\Gc_{\lambda}$, but only one block $J_{\lambda_{j_k} T_{j_k}}$.
	Therefore, when $p$ is large (usually requiring a large $m$), our method can significantly reduce the evaluation cost of $G$ and $T$.
	
	\beforesubsec
	\subsection{Convergence analysis}
	\aftersubsec
	For the analysis, let us denote $\widetilde{\Sc}_{\lambda, [j]}^k := [0, \cdots, 0, [\widetilde{\Sc}_{\lambda}^k ]_j, 0, \cdots, 0]$, then the main update of Algorithm~\ref{alg:A2} (line 7) can be expressed as
	\myeq{eq:VACOG}{
		\arraycolsep=0.2em
		\begin{array}{lcl}
			u^{k+1} &:=& u^k + \theta_k(u^k - u^{k-1}) - \frac{\eta_k}{\p_{j_k}} \widetilde{\Sc}_{\lambda, [j_k]}^k.
		\end{array}
		\tag{VCFOG}
	}
	This scheme can be viewed as an optimistic gradient variant for solving the BFS equation $\Gc_{\lambda}u^\star = 0$, integrating variance reduction, block-coordinate sampling, and Nesterov's momentum for acceleration.
	
	\vspace{0.75ex}
	\noindent$\mathrm{(a)}$~\textbf{Lyapunov function.}
	Let $\tau_{j} \in (0, \underline{\bar{\tau}}_j)$ be chosen arbitrarily, where $\underline{\bar{\tau}}_j := \min_{i \in [n]} \bar{\tau}_{ij}$ and $\bar{\tau}_{ij}$ is given in Lemma~\ref{lem:ch2:BFS_FNI_equivalence}.
	Let $e^k_j$, $\Delta_k$, $U_k$, $\kappa$, $\Theta$, and $\hat{\Theta}$ be given in Definition~\ref{def:VR_estimator} associated with the stochastic estimator $\widetilde{\Sc}^k_{\lambda}$ of $\Sc^k_{\lambda}$.
	We define the following functions:
	\myeq{eq:VACOG_Lyapunov_funcs}{
		\arraycolsep=0.2em
		\begin{array}{lcl}
			\Qc_k &:=& \norms{r(u^{k-1} - \uopt) + t_k(u^k - u^{k-1})}^2 + r \norms{u^{k-1} - \uopt}^2, \vspace{1ex}\\
			\Lc_k &:=& \Qc_k +  2r t_k \eta_{k-1} \left[\iprods{\Gc_{\lambda}u^{k-1}, u^{k-1} - \uopt} - \sum_{j=1}^m \tau_{j} \norms{[\Gc_{\lambda} u^{k-1}]_j}^2 \right], \vspace{1ex}\\
			\Pc_k &:=& \Lc_k + \frac{\eta^2 (1-\kappa)}{\kappa} (t_{k-1} - 1)^2 \Delta_{k-1} + \frac{\eta^2 \hat{\Theta}}{\kappa} (t_{k-1} - 1)^2 U_{k-1},
		\end{array}
	}
	where $r > 0$, $t_k > 0$, $\eta > 0$, and $\eta_k > 0$ are parameters specified later.
	At $k = 0$, we use the conventions $t_{-1} := r+1$ and $\eta_{-1} := \frac{\eta r}{r+2}$, which are consistent with \eqref{eq:th1_param} and give $t_0\eta_{-1} = r\eta$, together with $\Delta_{-1} = U_{-1} = 0$.
	
	Under the $\bar{\tau}$-average co-coercivity of $\Gc_{\lambda}$ from \eqref{eq:BFS_cocoercive}, it is easy to see that $\Qc_k \geq 0$, $\Lc_k \geq 0$, and hence, $\Pc_k \geq 0$.
	To analyze the convergence of Algorithm~\ref{alg:A2}, we use the \eqref{eq:VACOG} scheme, which can be viewed as a variance-reduced randomized block-coordinate fast OG method for solving the BFS equation $\Gc_{\lambda}u^\star = 0$, and it indirectly solves \eqref{eq:FNI} via a backward-forward splitting scheme.
	
	\vspace{0.75ex}
	\noindent\textbf{$\mathrm{(b)}$~Key estimates.}
	We first prove the following two key lemmas, which will be used in the proof of the main theorems.
	The proofs of these lemmas can be found in Appendix~\ref{apdx:subsec:ch2:VACOG_technical_lemmas}.
	The first lemma provides a lower bound for $\Lc_k - \Expsk{k}{\Lc_{k+1}}$.
	
	\begin{lemma}\label{le:VACOG_key_est1}
		Suppose that $\zer{\Psi} \ne \emptyset$ and Assumption \ref{as:VACOG} holds for \eqref{eq:FNI}.
		Let $\sets{u^k}$ be generated by Algorithm~\ref{alg:A2} and $\Lc_k$ be defined by \eqref{eq:VACOG_Lyapunov_funcs}.
		Let $\theta_k$ and $\gamma_k$ be updated by
		\myeqn{
			\arraycolsep=0.2em
			\begin{array}{lcl}
				\theta_k &:=& \frac{t_k - r - 1}{t_{k+1}} \quad \text{and} \quad \gamma_k := \frac{t_{k+1}\theta_k}{t_{k+1}\theta_k + r} \in (0,1),
			\end{array}
		}
		where $r > 0$ and $t_k > r + 1$ for all $k \geq 0$.
		Then, for all $k \geq 0$, we have
		\myeq{eq:VACOG_key_estimate1}{
			\arraycolsep=0.2em
			\begin{array}{lcl}
				\Lc_k - \Expsk{k}{\Lc_{k+1}} &\geq& [(t_k - r)^2 - t_{k+1}^2 \theta_k^2 +  r]\norms{u^k - u^{k-1}}^2  - t_{k+1}^2 \eta_k^2 \Expsk{k}{\norms{e^k}_{\beta}^2} \vspace{1ex}\\
				&& + {~} \frac{2t_{k+1}^2 \theta_k \eta_k}{n} \sum\limits_{i=1}^n \sum\limits_{j=1}^m (\bar{\tau}_{ij} - \tau_{j}) \norms{[\Gc_{i, \lambda} u^k]_j - [\Gc_{i, \lambda} u^{k-1}]_j}^2  \vspace{0ex}\\
				&& + {~} t_{k+1}\eta_k \sum\limits_{j=1}^m \frac{2\tau_{j} \p_j (t_{k+1}\theta_k + r) - t_{k+1}\eta_k}{\p_j} \norms{[\Gc_{\lambda} u^k]_j - \gamma_k [\Gc_{\lambda} u^{k-1}]_j}^2 \vspace{0ex}\\
				&& + {~} 2r \phi_k \Big[ \iprods{\Gc_{\lambda}u^{k-1}, u^{k-1} - \uopt} - \sum_{j=1}^m \tau_{j} \norms{[\Gc_{\lambda} u^{k-1}]_j}^2 \Big],
			\end{array}
		}
		where $\phi_k := t_k \eta_{k-1} - t_{k+1} \eta_k \gamma_k$.
	\end{lemma}
	
	Using this estimate, we can prove a quasi-descent property of $\Pc_k$ as follows.
	\begin{lemma}\label{lem:ch2:VACOG_lemma2}
		Suppose that $\zer{\Psi} \neq \emptyset$ and Assumption \ref{as:VACOG} holds for \eqref{eq:FNI}.
		Let $\sets{u^k}$ be generated by Algorithm~\ref{alg:A2} using the estimator $\widetilde{S}^k$ of $S^k$ satisfying Definition~\ref{def:VR_estimator}, and $\Pc_k$ be defined by \eqref{eq:VACOG_Lyapunov_funcs}.
		Assume that $\theta_k$, $\gamma_k$, and $\eta_k$ in Algorithm~\ref{alg:A2} are updated as
		\myeq{eq:lem2_param_conditions}{
			\arraycolsep=0.2em
			\begin{array}{lcl}
				\theta_k &:=& \frac{t_k - r - 1}{t_{k+1}}, \qquad \gamma_k := \frac{t_k - r - 1}{t_k - 1}, \quad \text{and} \quad  \eta_k := \frac{\eta (t_k - 1)}{t_{k+1}},
			\end{array}
		}
		where $r > 0$, $\eta > 0$, and $t_k > r + 1$ for all $k \geq 0$.
		Then, for all $k \geq 0$, we have
		\myeq{eq:VACOG_key_estimate2}{
			\hspace{-4ex}
			\arraycolsep=0.2em
			\begin{array}{lcl}
				\Expsk{k}{\Pc_{k+1}} &\leq& \big(1 + \frac{2r\eta^2 \bar{\p}}{\kappa \Lambda_0^2 }\delta_k\big)\Pc_k - (2t_k - r - 1)\norms{u^k - u^{k-1}}^2 \vspace{0ex}\\
				&& - {~} \eta (t_k - 1)^2 \sum\limits_{j=1}^m \frac{ 2\tau_{j} \p_j - \eta}{\p_j} \norms{[\Gc_{\lambda} u^k]_j - \gamma_k [\Gc_{\lambda} u^{k-1}]_j}^2 \vspace{0ex}\\
				&& \textls[-50]{$- {~} \frac{\eta (t_k - 1)^2}{n} \sum\limits_{i=1}^n \sum\limits_{j=1}^m \Big[ \frac{2\p_j(\bar{\tau}_{ij} - \tau_{j})(t_k - r - 1)}{t_k - 1} - \frac{\eta (\Theta + \hat{\Theta})}{\kappa} \Big] \frac{1}{\p_j}
					\norms{[\Gc_{i, \lambda} u^k]_j - [\Gc_{i, \lambda} u^{k-1}]_j}^2$}  \vspace{0ex}\\
				&& \textls[-40]{$- {~} 2r\eta (t_{k-1} - t_k + r) \Big[\iprods{\Gc_{\lambda} u^{k-1}, u^{k-1} - \uopt} - \sum\limits_{j=1}^m \tau_{j} \norms{[\Gc_{\lambda} u^{k-1}]_j}^2 \Big] + \frac{2\sigma_{\lambda, \star}^2r^2 \eta^2}{\kappa} \delta_k.$}
			\end{array}
			\hspace{-8ex}
		}
	\end{lemma}
	
	\vspace{0.5ex}
	\noindent\textbf{$\mathrm{(c)}$~Main result 2.}
	Now, for given $\kappa$, $\Theta$, and $\hat{\Theta}$ in Definition~\ref{def:VR_estimator}, $\p$ and $\bar{\p}$ in \eqref{eq:VACOG_prob}, $\Lambda_0$ in Subsection~\ref{subsec:ch2:FNI_reform}, $\bar{\tau}_{ij}$ in Lemma~\ref{lem:ch2:BFS_FNI_equivalence}, $0 < \tau_{j} < \bar{\tau}_{ij}$, $\eta > 0$, and $r > 2$, let us define the following constants:
	\myeq{eq:VACOG_constants1}{
		\hspace{-4ex}
		\arraycolsep=0.2em
		\begin{array}{lcllcllcl}
			\Lambda_1 &:=& \min\limits_{i,j} \{ (\bar{\tau}_{ij} - \tau_{j})\p_j \}, \quad&
			\Lambda_2 &:=& 2 \min\limits_{j} \{ \tau_{j} \p_j \}, \quad&
			\Lambda_3 &:=&  \frac{2\eta\Lambda_1}{r + 1} - \frac{\eta^2 (\Theta + \hat{\Theta})}{\kappa}, \\
			\Lambda_4 &:=& \textls[0]{$\eta \min\limits_{j}\sets{2\tau_{j} \p_j - \eta}$,} \quad
			&\bar{\eta} &:=& \textls[0]{$\min\set{ \frac{2 \kappa \Lambda_1}{(r+1)(\Theta + \hat{\Theta})}, \Lambda_2}$}, \quad
			&\underline{\bar{\tau}} &:=& \min\limits_{i, j} \{ \bar{\tau}_{ij} \}.
		\end{array}
		\hspace{-4ex}
	}
	We also denote $\underline{\tau} := \min_{j \in [m]}\tau_j$.
	Next, using \eqref{eq:VACOG_constants1}, $C_0 := r^2\left(\frac{\eta}{2\underline{\tau}} + 2\right)$,
	$\Sc_K := \sum_{k = 0}^K \delta_k$, $\hat{\Sc}_K := \frac{2r\eta^2 \bar{\p}}{\kappa \Lambda_0^2}\Sc_K$, we define
	\myeq{eq:VACOG_constants2}{
		\arraycolsep=0.2em
		\begin{array}{lcl}
			\Cc_{0,K}  &:=& \big(C_0 \norms{u^0 - \uopt}^2 + \frac{r\Lambda_0^2}{\bar{\p}} \sigma_{\lambda, \star}^2\big) \exp(\hat{\Sc}_K) - \frac{r\Lambda_0^2}{\bar{\p}} \sigma_{\lambda, \star}^2, \vspace{1ex}\\
			\Cc_{1,K} &:=& \frac{r^2 \eta^2}{\underline{\bar{\tau}}^2} \norms{u^0 - u^{\star}}^2 + \left[2 + \frac{\eta^2 (1 - \underline{\p})}{\Lambda_4} + \frac{\eta^2(\Theta + \hat{\Theta})}{\kappa \Lambda_3 }\right] \Cc_{0,K} + \hat{\Sc}_{K}\big(\Cc_{0,K} + \frac{r\Lambda_0^2}{\bar{\p}}\sigma_{\lambda, \star}^2\big), \vspace{1ex}\\
			\Cc_{2,K} &:=& \left(\frac{\eta^2}{\Lambda_4} + \frac{\eta^2 (\Theta + \hat{\Theta})}{\kappa \Lambda_3} + 8r \right)\Cc_{0,K} + \frac{8r}{r-2}\Cc_{1,K} + \hat{\Sc}_{K}\big(\Cc_{0,K} + \frac{r\Lambda_0^2}{\bar{\p}} \sigma_{\lambda, \star}^2\big).
		\end{array}
	}
	Then, we are ready to prove the convergence in expectation and almost-sure convergence of Algorithm~\ref{alg:A2} as stated in the following theorems.
	
	\begin{theorem}[Convergence in expectation]\label{th:VACOG_convergence}
		For  \eqref{eq:FNI}, suppose that $\zer{\Psi} \ne \emptyset$ and Assumption \ref{as:VACOG} holds.
		Let $\tau_{j}$ be such that $0 < \tau_{j} < \bar{\tau}_{ij}$.
		For given $r > 2$ and $0 < \eta < \bar{\eta}$, where $\bar{\eta}$ is defined in \eqref{eq:VACOG_constants1}, let $\sets{(x^k, u^k)}$ be generated by Algorithm~\ref{alg:A2} using:
		\myeq{eq:th1_param}{
			\arraycolsep=0.2em
			\begin{array}{lcl}
				\theta_k &:=& \frac{k + 1}{k + r + 3}, \hspace{1.5ex} \gamma_k := \frac{k + 1}{k + r + 1},
				\hspace{1.5ex} \text{and} \hspace{1.5ex} \eta_k := \frac{\eta (k + r + 1)}{k + r + 3}.
			\end{array}
		}
		Then, for $\Cc_{0,K}$, $\Cc_{1,K}$, and $\Cc_{2,K}$ defined in \eqref{eq:VACOG_constants2}, we have
		\myeq{eq:th1_summability_bounds}{
			\arraycolsep=0.2em
			\begin{array}{lcl}
				\sum\limits_{k = 1}^{K} (2k + r + 3)\ \Expk{\norms{u^k - u^{k-1}}^2} &\leq& \Cc_{0,K}, \vspace{0ex}\\
				\sum\limits_{k = 0}^K (k + r + 2)\ \Expk{\norms{\Gc_{\lambda}u^{k}}^2} &\leq& \frac{8}{\eta^2} \left(\Cc_{0,K} + \frac{\Cc_{1,K}}{r-2} \right).
			\end{array}
		}
		Moreover, let $\xi^k := \lambda^{-1}\circ(u^k - x^k)$, so that $\xi^k \in Tx^k$.
		Then, for all $K \geq 0$, we have
		\myeq{eq:th1_BigO_convergence_rates}{
			\arraycolsep=0.2em
			\begin{array}{lcl}
				\Expk{\norms{x^{K+1} - x^K}^2}  \leq   \frac{\Cc_{2,K}}{ { (K+r+3)^2 }} \quad \text{and} \quad
				\Expk{\norms{Gx^{K} + \xi^K}^2}  \leq   \frac{8(\Cc_{1,K} + \Cc_{2,K})}{\eta^2 { (K+r+3)(K+r) }}.
			\end{array}
		}
		In particular, if $\Cc_{1,K}$ and $\Cc_{2,K}$ are bounded uniformly in $K$, then both quantities in \eqref{eq:th1_BigO_convergence_rates} converge to zero at the rate $\BigOs{1/K^2}$.
		If, in addition, $\sum_{k = 0}^\infty \delta_k < +\infty$, we also have the following $\SmallOs{1/k^2}$ convergence results:
		\myeq{eq:th1_small_o_convergence_rates}{
			\arraycolsep=0.2em
			\begin{array}{lcl}
				\lim\limits_{k \to \infty} { k^2 } \Expk{\norms{x^{k+1} - x^k}^2} &=& 0
				\quad \text{and} \quad
				\lim\limits_{k \to \infty} { k^2 } \Expk{\norms{Gx^k + \xi^k}^2} = 0, \ \xi^k \in Tx^k.
			\end{array}
		}
		Finally, if $\Cc_{1,K} + \Cc_{2,K} = \BigOs{1}$, then for a given $\epsilon > 0$, Algorithm~\ref{alg:A2} requires at most $\BigO{\epsilon^{-1}}$ iterations to achieve $\Expk{\norms{Gx^K + \xi^K}^2} \leq \epsilon^2$. 
	\end{theorem}
	
	\begin{proof}
		With $t_k := k + r + 2$, the parameters in \eqref{eq:th1_param} coincide with those in \eqref{eq:lem2_param_conditions}.
		Since $\frac{k+1}{k+r+1} \geq \frac{1}{r+1}$ and $0 < \eta < \bar{\eta}$, the constants $\Lambda_3$ and $\Lambda_4$ in \eqref{eq:VACOG_constants1} are positive.
		We divide the proof into four steps.

		\vspace{0.5ex}
		\noindent\textit{\textbf{\underline{Step 1}.}~Lyapunov and variance bounds.}
		Let $\hat{\delta}_k := \frac{2r\eta^2 \bar{\p}}{\kappa \Lambda_0^2}\delta_k$.
		With $t_k = k + r + 2$, we have $t_{k-1} - t_k + r = r - 1 > 0$, and the bracket in the last line of \eqref{eq:VACOG_key_estimate2} is nonnegative by \eqref{eq:BFS_cocoercive}, Jensen's inequality, and $\tau_j < \bar{\tau}_{ij}$.
		Moreover, by the definitions of $\Lambda_1$ and $\Lambda_3$ in \eqref{eq:VACOG_constants1}, we have
		\myeqn{
			\arraycolsep=0.2em
			\begin{array}{lcl}
				\frac{2\p_j(\bar{\tau}_{ij} - \tau_{j})(k + 1)}{k + r + 1} - \frac{\eta (\Theta + \hat{\Theta})}{\kappa} &\geq& \frac{2\Lambda_1}{r + 1} - \frac{\eta (\Theta + \hat{\Theta})}{\kappa} = \frac{\Lambda_3}{\eta}.
			\end{array}
		}
		Using these facts, $\beta_j = \p_j^{-1}$, and the definition of $\Lambda_4$, \eqref{eq:VACOG_key_estimate2} implies
		\myeq{eq:th1_proof0}{
			\arraycolsep=0.2em
			\begin{array}{lcl}
				\Expsk{k}{\Pc_{k+1}} &\leq& (1 + \hat{\delta}_k)\Pc_k - \Dc_k + \frac{r\Lambda_0^2}{\bar{\p}}\sigma_{\lambda, \star}^2\hat{\delta}_k,
			\end{array}
		}
		where $\Dc_k := (2k + r + 3)\norms{u^k - u^{k-1}}^2 + (k + r + 1)^2 \big(\Lambda_3 U_k + \Lambda_4 \norms{\Gc_{\lambda}u^k - \gamma_k\Gc_{\lambda}u^{k-1}}_{\beta}^2\big)$.
		Let $\hat{\Pc}_k := \Pc_k + \frac{r\Lambda_0^2}{\bar{\p}}\sigma_{\lambda, \star}^2$, $A_0 := 1$, and $A_{k+1} := \prod_{l=0}^{k}(1+\hat{\delta}_l)$.
		Taking the total expectation of \eqref{eq:th1_proof0} and dividing the result by $A_{k+1}$, we obtain
		\myeqn{
			\arraycolsep=0.2em
			\begin{array}{lcl}
				\frac{\Expk{\hat{\Pc}_{k+1}}}{A_{k+1}} &\leq& \frac{\Expk{\hat{\Pc}_k}}{A_k} - \frac{\Expk{\Dc_k}}{A_{k+1}}.
			\end{array}
		}
		We sum this inequality from $k = 0$ to $k = K$ and multiply the result by $A_{K+1}$.
		Since $1 \leq A_{k+1} \leq A_{K+1} \leq \exp(\hat{\Sc}_K)$ for $0 \leq k \leq K$ and $\Pc_0 \leq C_0\norms{u^0 - \uopt}^2$ by \eqref{eq:th1_boundP0} of Lemma~\ref{lem:ch2:VACOG_P0_bound}, we get
		\myeq{eq:th1_proof1.1}{
			\arraycolsep=0.2em
			\begin{array}{lcl}
				\Expk{\Pc_{K+1}} + \sum_{k=0}^K \Expk{\Dc_k} &\leq& \exp(\hat{\Sc}_K)\big(\Pc_0 + \frac{r\Lambda_0^2}{\bar{\p}}\sigma_{\lambda, \star}^2\big) - \frac{r\Lambda_0^2}{\bar{\p}}\sigma_{\lambda, \star}^2 \leq \Cc_{0,K}.
			\end{array}
		}
		Since $\Cc_{0,K}$ is nondecreasing in $K$, \eqref{eq:th1_proof1.1} yields $\Expk{\Pc_k} \leq \Cc_{0,K}$ for $0 \leq k \leq K + 1$, and
		\myeq{eq:VACOG_summability_results1}{
			\arraycolsep=0.2em
			\begin{array}{lcl}
				\sum_{k = 0}^K (2k + r + 3) \Expk{\norms{u^k - u^{k-1}}^2} &\leq& \Cc_{0,K}, \vspace{1ex}\\
				\sum_{k = 0}^K (k + r + 1)^2 \Expk{U_k} &\leq& \frac{\Cc_{0,K}}{\Lambda_3}, \vspace{1ex}\\
				\sum_{k = 0}^K (k + r + 1)^2 \Expk{\norms{\Gc_{\lambda}u^k - \gamma_k\Gc_{\lambda}u^{k-1}}_{\beta}^2} &\leq& \frac{\Cc_{0,K}}{\Lambda_4}.
			\end{array}
		}
		Since $u^{-1} = u^0$, the first bound in \eqref{eq:VACOG_summability_results1} is the first inequality in \eqref{eq:th1_summability_bounds}.

		Next, we bound the variance terms.
		We take the total expectation of \eqref{eq:th1_proof2} in Lemma~\ref{lem:ch2:VACOG_Wk_bound} and sum it from $k = 1$ to $k = K$; the telescoping terms are nonpositive because $\Delta_0 = U_0 = 0$.
		By \eqref{eq:lemma2_proof3} and $\Expk{\Pc_k} \leq \Cc_{0,K}$, we also have $\sum_{k=1}^K \delta_k\Expk{V_{k-1}} \leq 2\Sc_K\big(\frac{\bar{\p}\Cc_{0,K}}{r\Lambda_0^2} + \sigma_{\lambda,\star}^2\big)$.
		Combining these two facts with the second bound in \eqref{eq:VACOG_summability_results1}, we obtain
		\myeq{eq:th1_bound_Delta_k}{
			\arraycolsep=0.2em
			\begin{array}{lcl}
				\sum_{k=0}^K (k + r + 1)^2 \Expk{\Delta_k} &\leq& B_K := \frac{\Theta + \hat{\Theta}}{\kappa\Lambda_3}\Cc_{0,K} + \frac{\hat{\Sc}_K}{\eta^2}\Big(\Cc_{0,K} + \frac{r\Lambda_0^2}{\bar{\p}}\sigma_{\lambda, \star}^2\Big).
			\end{array}
		}
		\vspace{0.5ex}
		\noindent\textit{\textbf{\underline{Step 2}.}~Summability of the residual.}
		Taking the total expectation of \eqref{eq:th1_proof4} in Lemma~\ref{lem:ch2:VACOG_Wk_bound} and summing it from $k = 0$ to $k = K$, we obtain
		\myeq{eq:th1_proof5}{
			\arraycolsep=0.2em
			\begin{array}{lcl}
				(K + r + 3)(K + r)\Expk{\norms{u^{K+1} - u^K + \eta_K\Gc_{\lambda}u^K}^2} && \vspace{1ex}\\
				\quad + {~} (r - 2)\sum_{k=0}^K (k + r + 2)\Expk{\norms{u^k - u^{k-1} + \eta_{k-1}\Gc_{\lambda}u^{k-1}}^2} &\leq& \Cc_{1,K}.
			\end{array}
		}
		Indeed, since $u^{-1} = u^0$, the initial term is $(r+2)(r-1)\eta_{-1}^2\norms{\Gc_{\lambda}u^0}^2$.
		Since $(r+2)(r-1)\eta_{-1}^2 \leq r^2\eta^2$, this term is at most $\frac{r^2\eta^2}{\underline{\bar{\tau}}^2}\norms{u^0 - \uopt}^2$ by \eqref{eq:th1_bound_Gu0}.
		Since $\frac{4}{r}\big(k + r + 2 - \frac{2}{r}\big) \leq 2(2k + r + 3)$ for $r > 2$, the remaining terms are bounded by $2\Cc_{0,K} + \frac{(1 - \underline{\p})\eta^2}{\Lambda_4}\Cc_{0,K} + \eta^2B_K$, owing to \eqref{eq:VACOG_summability_results1} and \eqref{eq:th1_bound_Delta_k}.
		The sum of these bounds is exactly $\Cc_{1,K}$.

		We now shift the indices in \eqref{eq:VACOG_summability_results1} and \eqref{eq:th1_proof5}.
		Minimizing the quadratic function $\Qc_{K+1}$ with respect to $u^K - \uopt$ gives $\Pc_{K+1} \geq \Qc_{K+1} \geq \frac{(K+r+3)^2}{r+1}\norms{u^{K+1} - u^K}^2$.
		Since $(r+1)(K+r+2) \leq (K+r+3)^2$, this bound and \eqref{eq:th1_proof1.1} imply the first estimate below.
		Similarly, since $(r-2)(K+r+2) \leq (K+r+3)(K+r)$, \eqref{eq:th1_proof5} implies the second one:
		\myeqn{
			\arraycolsep=0.2em
			\begin{array}{lcl}
				\sum_{k=0}^K (k + r + 2)\Expk{\norms{u^{k+1} - u^k}^2} &\leq& \Cc_{0,K}, \vspace{1ex}\\
				\sum_{k=0}^K (k + r + 2)\Expk{\norms{u^{k+1} - u^k + \eta_k\Gc_{\lambda}u^k}^2} &\leq& \frac{\Cc_{1,K}}{r - 2}.
			\end{array}
		}
		Finally, since $\eta_k \geq \frac{\eta}{2}$ and $\eta_k\Gc_{\lambda}u^k = (u^{k+1} - u^k + \eta_k\Gc_{\lambda}u^k) - (u^{k+1} - u^k)$, Young's inequality gives
		\myeq{eq:th1_proof5.5}{
			\arraycolsep=0.2em
			\begin{array}{lcl}
				\frac{\eta^2}{4}\norms{\Gc_{\lambda}u^k}^2 &\leq& 2\norms{u^{k+1} - u^k + \eta_k\Gc_{\lambda}u^k}^2 + 2\norms{u^{k+1} - u^k}^2.
			\end{array}
		}
		Multiplying \eqref{eq:th1_proof5.5} by $k + r + 2$ and summing the result, we obtain the second inequality in \eqref{eq:th1_summability_bounds}.

		\vspace{0.5ex}
		\noindent\textit{\textbf{\underline{Step 3}.}~Pointwise rates.}
		Taking the total expectation of \eqref{eq:th1_proof6} in Lemma~\ref{lem:ch2:VACOG_Zk}, summing it from $k = 0$ to $k = K$, and using $u^{-1} = u^0$, \eqref{eq:VACOG_summability_results1}, \eqref{eq:th1_bound_Delta_k}, the second inequality in \eqref{eq:th1_summability_bounds}, and $k + 1 \leq k + r + 2$, we obtain
		\myeqn{
			\arraycolsep=0.2em
			\begin{array}{lcl}
				(K + r + 3)^2\Expk{\norms{u^{K+1} - u^K}^2} &\leq& \frac{\eta^2}{\Lambda_4}\Cc_{0,K} + \eta^2B_K + 8r\Big(\Cc_{0,K} + \frac{\Cc_{1,K}}{r - 2}\Big) = \Cc_{2,K}.
			\end{array}
		}
		Combining this bound with \eqref{eq:th1_proof5} and \eqref{eq:th1_proof5.5}, we get
		\myeqn{
			\arraycolsep=0.2em
			\begin{array}{lcl}
				\Expk{\norms{\Gc_{\lambda}u^K}^2} &\leq& \frac{8}{\eta^2}\Big(\frac{\Cc_{1,K}}{(K+r+3)(K+r)} + \frac{\Cc_{2,K}}{(K+r+3)^2}\Big) \leq \frac{8(\Cc_{1,K} + \Cc_{2,K})}{\eta^2(K+r+3)(K+r)}.
			\end{array}
		}
		Since $J_{\lambda \circ T}$ is nonexpansive by Lemma~\ref{lem:ch2:BFS_FNI_equivalence}(ii), we have
		\myeq{eq:th1_proof8}{
			\arraycolsep=0.2em
			\begin{array}{lcl}
				\norms{x^{K+1} - x^K} &=& \norms{J_{\lambda \circ T}u^{K+1} - J_{\lambda \circ T}u^K} \leq \norms{u^{K+1} - u^K}.
			\end{array}
		}
		Moreover, $\norms{\Gc_{\lambda}u^K} = \norms{Gx^K + \xi^K}$ by Lemma~\ref{lem:ch2:BFS_FNI_equivalence}(i).
		These two facts prove \eqref{eq:th1_BigO_convergence_rates}.
		The iteration-complexity statement follows directly from the second bound in \eqref{eq:th1_BigO_convergence_rates}.

		\vspace{0.5ex}
		\noindent\textit{\textbf{\underline{Step 4}.}~Small-$o$ rates.}
		If $\sum_{k=0}^\infty \delta_k < +\infty$, then $\hat{\Sc}_K$ is bounded, and so are $\Cc_{0,K}$, $B_K$, $\Cc_{1,K}$, and $\Cc_{2,K}$ uniformly in $K$.
		After taking the total expectation of \eqref{eq:VACOG_Wk_bound} and \eqref{eq:VACOG_Zk_bound}, all nonnegative terms on their right-hand sides, except $\Expk{W_{k-1}}$ and $\Expk{Z_{k-1}}$, are summable by \eqref{eq:th1_summability_bounds}, \eqref{eq:VACOG_summability_results1}, \eqref{eq:th1_bound_Delta_k}, and the bound on $\sum_k\delta_k\Expk{V_{k-1}}$ in Step~1.
		By Lemma~\ref{lem:ch2:tech_convergence1}, both $\Expk{W_k}$ and $\Expk{Z_k}$ converge.
		Since $(k + r + 1)^2\Expk{\Delta_k} \to 0$ and $(k + r + 1)^2\Expk{U_k} \to 0$ by \eqref{eq:VACOG_summability_results1} and \eqref{eq:th1_bound_Delta_k}, the limits of $(k + r + 3)(k + r)\Expk{\norms{u^{k+1} - u^k + \eta_k\Gc_{\lambda}u^k}^2}$ and $(k + r + 3)^2\Expk{\norms{u^{k+1} - u^k}^2}$ also exist.
		Combining this with the weighted summability bounds in Step~2 and Lemma~\ref{lem:ch2:tech_convergence2}, both limits are zero.
		Then, \eqref{eq:th1_proof5.5} gives $k^2\Expk{\norms{\Gc_{\lambda}u^k}^2} \to 0$.
		As in Step~3, \eqref{eq:th1_proof8} and Lemma~\ref{lem:ch2:BFS_FNI_equivalence}(i) transfer these two limits to \eqref{eq:th1_small_o_convergence_rates}.
	\end{proof}
		
	\begin{theorem}[Almost-sure convergence]\label{th:VACOG_convergence_as}
		Under the hypotheses of Theorem~\ref{th:VACOG_convergence}, if additionally $\sum_{k = 0}^\infty \delta_k < +\infty$, then
		the following limits hold almost surely: \vspace{-1ex}
		\myeq{eq:th2_almost_sure_convergence}{
			\arraycolsep=0.2em
			\begin{array}{lcl}
				\lim\limits_{k\to \infty} k^2 {\norms{x^{k+1} - x^k}^2} &=& 0
				\qquad \text{and} \qquad
				\lim\limits_{k\to \infty} k^2 {\norms{Gx^k + \xi^k}^2} = 0, \ \xi^k \in Tx^k.
			\end{array}
		}
		Moreover,  $\sets{x^k}$ converges almost surely to a $\zer{\Psi}$-valued random variable $x^{\star}  \in \zer{\Psi}$, a solution of \eqref{eq:FNI}.
	\end{theorem}
	
	\begin{proof}
		Fix $\uopt \in \zer{\Gc_{\lambda}}$.
		Each statement below holds almost surely, and since only countably many of them are involved, they hold simultaneously with probability one.

		\vspace{0.5ex}
		\noindent\textit{\textbf{\underline{Step 1}.}~Summability.}
		Since $\sum_{k=0}^\infty \hat{\delta}_k < +\infty$, applying Lemma~\ref{lem:ch2:supermartingale} to \eqref{eq:th1_proof0}, we conclude that $\Pc_k$ converges and $\sum_{k=0}^\infty \Dc_k < +\infty$, that is,
		\myeq{eq:th1_proof9}{
			\arraycolsep=0.2em
			\begin{array}{lcl}
				\sum_{k=1}^\infty (2k + r + 3)\norms{u^k - u^{k-1}}^2 & < & +\infty, \vspace{1ex}\\
				\sum_{k=1}^\infty (k + r + 1)^2 U_k  & < &  +\infty, \vspace{1ex}\\
				\sum_{k=1}^\infty (k + r + 1)^2 \norms{\Gc_{\lambda}u^k - \gamma_k\Gc_{\lambda}u^{k-1}}_{\beta}^2 & < & +\infty.
			\end{array}
		}
		Moreover, \eqref{eq:lemma2_proof3} gives $V_{k-1} \leq \frac{2\bar{\p}}{r\Lambda_0^2}\Pc_k + 2\sigma_{\lambda,\star}^2$, so $\sup_k V_k < +\infty$ and $\sum_{k=0}^\infty \delta_k V_{k-1} < +\infty$.
		Since $\nu_k = \frac{k+r}{k+r+1}$ and $(1 - \gamma_k)^2(k + r + 1)^2 = r^2$, multiplying \eqref{eq:VR_estimator3} by $(k + r + 1)^2$ gives
		\myeqn{
			\arraycolsep=0.2em
			\begin{array}{lcl}
				(k + r + 1)^2 \Expsk{k}{\Delta_k} &\leq& (k+r)^2\Delta_{k-1} - \kappa(k+r)^2\Delta_{k-1} + \Theta (k + r + 1)^2 U_k \vspace{1ex}\\
				&& + {~} \hat{\Theta} (k + r)^2 U_{k-1} + r^2 \delta_k V_{k-1}.
			\end{array}
		}
		The last three terms are summable, so Lemma~\ref{lem:ch2:supermartingale} yields $\sum_{k=0}^\infty (k + r + 1)^2\Delta_k < +\infty$.
		In particular, $(k + r + 1)^2\Delta_k \to 0$ and $(k + r + 1)^2U_k \to 0$.

		\vspace{0.5ex}
		\noindent\textit{\textbf{\underline{Step 2}.}~Rates of the residual and the increment.}
		By \eqref{eq:th1_proof9} and Step~1, the last three terms on the right-hand side of \eqref{eq:VACOG_Wk_bound} are summable.
		By Lemma~\ref{lem:ch2:supermartingale}, $W_k$ converges and $\sum_{k} (k + r + 2)\norms{u^k - u^{k-1} + \eta_{k-1}\Gc_{\lambda}u^{k-1}}^2 < +\infty$.
		Since the variance terms in $W_k$ tend to zero, $(k + r + 3)(k + r)\norms{u^{k+1} - u^k + \eta_k\Gc_{\lambda}u^k}^2$ also converges, and Lemma~\ref{lem:ch2:tech_convergence2} shows that its limit is zero.
		Using \eqref{eq:th1_proof5.5} and \eqref{eq:th1_proof9}, we also get $\sum_{k} (k + 1)\norms{\Gc_{\lambda}u^k}^2 < +\infty$.
		Therefore, the last four terms on the right-hand side of \eqref{eq:VACOG_Zk_bound} are summable.
		By the same argument, $Z_k$ converges, $\sum_{k} (k + r)\norms{u^k - u^{k-1}}^2 < +\infty$, and $(k + r + 3)^2\norms{u^{k+1} - u^k}^2 \to 0$.
		Together with \eqref{eq:th1_proof5.5}, this gives $k^2\norms{\Gc_{\lambda}u^k}^2 \to 0$.
		Since $J_{\lambda\circ T}$ is nonexpansive and $\norms{\Gc_{\lambda}u^k} = \norms{Gx^k + \xi^k}$ by Lemma~\ref{lem:ch2:BFS_FNI_equivalence}, we obtain \eqref{eq:th2_almost_sure_convergence}.

		\vspace{0.5ex}
		\noindent\textit{\textbf{\underline{Step 3}.}~Convergence of the iterates.}
		Since $r\norms{u^{k-1} - \uopt}^2 \leq \Pc_k$ and $\Pc_k$ converges, the sequence $\sets{u^k}$ is bounded.
		With $t_k = k + r + 2$, we can expand $\Pc_k$ in \eqref{eq:VACOG_Lyapunov_funcs} as
		\myeqn{
			\arraycolsep=0.2em
			\begin{array}{lcl}
				\Pc_k &=& r(r+1)\norms{u^{k-1} - \uopt}^2 + 2rt_k \iprods{u^{k-1} - \uopt, u^k - u^{k-1}} + t_k^2 \norms{u^k - u^{k-1}}^2 \vspace{1ex}\\
				&& + {~} 2r t_k \eta_{k-1} \Big[\iprods{\Gc_{\lambda}u^{k-1}, u^{k-1} - \uopt} - \sum_{j=1}^m \tau_{j} \norms{[\Gc_{\lambda}u^{k-1}]_j}^2 \Big] \vspace{1ex}\\
				&& + {~} \frac{\eta^2}{\kappa}(k + r)^2\big[(1-\kappa)\Delta_{k-1} + \hat{\Theta}U_{k-1}\big].
			\end{array}
		}
	By Step~2, $t_k(u^k - u^{k-1}) \to 0$ and $t_k\Gc_{\lambda}u^{k-1} \to 0$.
	Since $\sets{u^k}$ is bounded and $\eta_{k-1} \leq \eta$, the Cauchy--Schwarz inequality shows that all terms on the right-hand side except the first one converge to zero.
	Since $\Pc_k$ converges, so does $\norms{u^{k-1} - \uopt}$.
	As $\uopt \in \zer{\Gc_{\lambda}}$ is arbitrary, for every $\uopt \in \zer{\Gc_{\lambda}}$, the limit $\lim_{k\to\infty}\norms{u^k - \uopt}$ exists almost surely.
	Since $\Gc_{\lambda}$ is continuous by Lemma~\ref{lem:ch2:BFS_FNI_equivalence}(iii) and $\norms{\Gc_{\lambda}u^k} \to 0$, Lemma~\ref{lem:ch2:Opial} implies that $\sets{u^k}$ converges almost surely to a $\zer{\Gc_{\lambda}}$-valued random variable $\bar{u}^{\star}$.
	Finally, by the continuity of $J_{\lambda\circ T}$ and Lemma~\ref{lem:ch2:BFS_FNI_equivalence}(i), $\sets{x^k}$ converges almost surely to $x^{\star} := J_{\lambda\circ T}\bar{u}^{\star} \in \zer{\Psi}$.
	\end{proof}
	{
		\beforesubsec
		\subsection{Oracle complexity with importance-based sampling strategy}\label{subsec:ch2:oracle_complexity_analysis}
		\aftersubsec
		We now estimate the oracle complexity of Algorithm~\ref{alg:A2} using the estimators in Section~\ref{sec:ch2:OG_operator}.
		Recall that $\bar{L}_j := \max_{i \in [n]} L_{ij}$, then we choose $\lambda_j := \frac{2}{\bar{L}_j} \in \Big[2\rho_j, \frac{2+2\sqrt{1-\bar{L}_j\rho_j}}{\bar{L}_j}\Big)$ due to $\bar{L}_j \rho_j < 1$.
		We can also show that $\underline{\bar{\tau}}_j = \frac{1}{\bar{L}_j}$.
		Let us fix $\tau_j := \frac{\underline{\bar\tau}_j}{3} = \frac{1}{3\bar L_j}$ for all $j \in [m]$, choose $r := 3$, and consider the importance-based sampling $\p_j := \frac{\bar{L}_j}{\norms{\bar{L}}_1}$ for $j \in [m]$.
		Then, we can compute $\bar{\p}$, $\Lambda_0$, the constants in \eqref{eq:VACOG_constants1}, and \eqref{eq:VACOG_constants2} as follows:
		\myeq{eq:co_VACOG_complexity_quantities}{
			\arraycolsep=0.2em
			\left\{
			\begin{array}{lcl}
				\bar{\p} &=& \frac{\norms{\bar{L}}_{\infty}}{\norms{\bar{L}}_1}, \quad
				\underline{\bar{\tau}} = \min_{i,j} \bar{\tau}_{ij} = \frac{1}{\norms{\bar{L}}_{\infty}}, \quad
				\underline{\tau} = \min_j \tau_j = \frac{1}{3\norms{\bar{L}}_{\infty}}, \quad
				\Lambda_0 = \frac{1}{\norms{\bar{L}}_1}, \vspace{1ex}\\
				\Lambda_1 &=& \frac{2}{3\norms{\bar{L}}_1}, \quad
				\Lambda_2 = \frac{2}{3\norms{\bar{L}}_1}, \quad
				\Lambda_3 = \frac{\eta}{3\norms{\bar{L}}_1} - \frac{\eta^2 (\Theta + \hat{\Theta})}{\kappa}, \quad
				\Lambda_4 = \eta\big(\frac{2}{3\norms{\bar{L}}_1} - \eta\big), \vspace{1ex}\\
				\bar{\eta} &=& \frac{1}{3\norms{\bar{L}}_1} \min \set{\frac{\kappa}{\Theta + \hat{\Theta}}, 2}, \quad
				C_0 = \frac{27\norms{\bar{L}}_{\infty}\eta}{2} + 18.
			\end{array}
			\right.
		}
		With these choices, the next lemma shows that $\Cc_{0,K}$, $\Cc_{1,K}$, and $\Cc_{2,K}$ in \eqref{eq:VACOG_constants2} are bounded uniformly in $K$.
		\begin{lemma}\label{lem:ch2:VACOG_bound_C}
			Let the parameters be chosen as in~\eqref{eq:co_VACOG_complexity_quantities} with $\eta := \frac{\bar{\eta}}{2}$.
			Then, for $\Cc_{0,K}$, $\Cc_{1,K}$, and $\Cc_{2,K}$ defined in \eqref{eq:VACOG_constants2}, the following bounds hold.
			\begin{compactitem}
			\itemsep=0.4em
			\item[$\mathrm{(a)}$] 
			If the mini-batch estimator \eqref{eq:minibatch_estimator} is used and $\hat{\Sc}_K \leq \frac{1}{12}$, then $\Cc_{0,K} \leq 20.8\Rc_{\lambda}^2$, $\Cc_{1,K} \leq 67.5\Rc_{\lambda}^2$, and $\Cc_{2,K} \leq 2145\Rc_{\lambda}^2$, where $\Rc_{\lambda}^2 := \norms{u^0 - u^{\star}}^2 + \frac{\sigma_{\lambda, \star}^2}{\norms{\bar{L}}_1 \norms{\bar{L}}_{\infty}}$.
			
			\item[$\mathrm{(b)}$] 
			If the \ref{eq:loopless_SVRG} or \ref{eq:JacSketch_estimator} estimator is used with $1 \leq b\q^2 \leq 32$, then $\Cc_{0,K} \leq 22.5\norms{u^0 - u^{\star}}^2$, $\Cc_{1,K} \leq 91\norms{u^0 - u^{\star}}^2$, and $\Cc_{2,K} \leq 2769\norms{u^0 - u^{\star}}^2$.
			\end{compactitem}
		\end{lemma}
		\begin{proof}
			\noindent\textbf{$\mathrm{(i)}$~Mini-batch estimator.}
			First, using $\kappa = 1$, $\Theta = 2$, $\hat{\Theta} = 0$, we get $\bar{\eta} = \frac{1}{6\norms{\bar{L}}_1}$.
			We choose $\eta := \frac{\bar{\eta}}{2} = \frac{1}{12\norms{\bar{L}}_1}$, leading to $\Lambda_3 = \frac{1}{72\norms{\bar{L}}_1^2}$ and $\Lambda_4  = \frac{7}{144\norms{\bar{L}}_1^2}$.
			Thus, $C_0 = \frac{9\norms{\bar{L}}_{\infty}}{8\norms{\bar{L}}_1} + 18 \leq \frac{9}{8} + 18 = 19.125$.
			Moreover, $\hat{\Sc}_K = \frac{2r\eta^2 \bar{\p}}{\kappa \Lambda_0^2}\Sc_K = \frac{\bar{\p}}{24}\Sc_K$.
			Therefore, using $\hat{\Sc}_K \leq \frac{1}{12}$, we can evaluate that
			\myeqn{
				\arraycolsep=0.2em
				\begin{array}{lcl}
					\Cc_{0,K}  &=& \big(C_0 \norms{u^0 - \uopt}^2 + \frac{3\sigma_{\lambda, \star}^2}{\norms{\bar{L}}_1 \norms{\bar{L}}_{\infty}}\big) \exp(\hat{\Sc}_K) - \frac{3\sigma_{\lambda, \star}^2}{\norms{\bar{L}}_1 \norms{\bar{L}}_{\infty}} \vspace{1ex}\\
					&\leq& 19.125\exp(\frac{1}{12})\norms{u^0 - u^{\star}}^2 + 3\big(\exp(\frac{1}{12}) - 1\big)\frac{\sigma_{\lambda, \star}^2}{\norms{\bar{L}}_1 \norms{\bar{L}}_{\infty}} \leq 20.8\Rc_{\lambda}^2, \vspace{1ex}\\
					\Cc_{1,K} &=& \frac{\norms{\bar{L}}_{\infty}^2}{16\norms{\bar{L}}_1^2}\norms{u^0 - u^{\star}}^2 + \left[3 + \frac{1 - \underline{\p}}{7} \right] \Cc_{0,K} + \hat{\Sc}_{K}\big(\Cc_{0,K} + \frac{3\sigma_{\lambda, \star}^2}{\norms{\bar{L}}_1 \norms{\bar{L}}_{\infty}}\big) \vspace{1ex}\\
					&\leq&\frac{1}{16}\norms{u^0 - u^{\star}}^2 + \frac{22}{7}\Cc_{0,K} + \frac{1}{12}\big(\Cc_{0,K} + \frac{3\sigma_{\lambda, \star}^2}{\norms{\bar{L}}_1 \norms{\bar{L}}_{\infty}} \big) \vspace{1ex}\\
					&\leq& \big[\frac{1}{16} + \big(\frac{22}{7} + \frac{1}{12}\big)20.8 + \frac{1}{4}\big]\Rc_{\lambda}^2 \leq 67.5 \Rc_{\lambda}^2,\vspace{1ex}\\
					\Cc_{2,K} &=& \left(\frac{1}{7} + 1 + 24 \right)\Cc_{0,K} + 24\Cc_{1,K} + \hat{\Sc}_{K}\big(\Cc_{0,K} + \frac{3\sigma_{\lambda, \star}^2}{\norms{\bar{L}}_1 \norms{\bar{L}}_{\infty}} \big) \vspace{1ex}\\
					&\leq& \frac{176}{7}\Cc_{0,K} + 24\Cc_{1,K} + \frac{1}{12}\big(\Cc_{0,K} + \frac{3\sigma_{\lambda, \star}^2}{\norms{\bar{L}}_1 \norms{\bar{L}}_{\infty}} \big) \vspace{1ex}\\
					&\leq& \big[\big(\frac{176}{7} + \frac{1}{12}\big)20.8 + 24\cdot 67.5 + \frac{1}{4}\big]\Rc_{\lambda}^2 \leq 2145\Rc_{\lambda}^2.
				\end{array}
			}
			\noindent\textbf{$\mathrm{(ii)}$~\ref{eq:loopless_SVRG} and \ref{eq:JacSketch_estimator}.}
			First, using $\kappa := \frac{\q}{2}$, $\Theta = \hat{\Theta} := \frac{4}{b\q}$, $\Lambda_1$ and $\Lambda_2$ in \eqref{eq:co_VACOG_complexity_quantities}, and $1 \leq b\q^2 \leq 32$, we can compute from \eqref{eq:VACOG_constants1} that $\bar{\eta} = \frac{b\q^2}{48\norms{\bar{L}}_1}$.
			We choose $\eta := \frac{\bar{\eta}}{2} = \frac{b\q^2}{96\norms{\bar{L}}_1}$, leading to $\Lambda_3 = \frac{b\q^2}{576\norms{\bar{L}}_1^2}$ and $\Lambda_4 = \frac{b\q^2 (64 - b\q^2)}{9216\norms{\bar{L}}_1^2}$.
			Using $1 \leq b\q^2 \leq 32$ and $\bar{\p} \leq 1$, we have $C_0 = \frac{9b\q^2 \bar{\p}}{64} + 18 \leq \frac{9}{2} + 18 = 22.5$.
			Moreover, since $\delta_k = 0$, we have $\Sc_K = \hat{\Sc}_K = 0$.
			Thus, we have $\Cc_{0,K} = C_0 \norms{u^0 - u^{\star}}^2 \leq 22.5\norms{u^0 - u^{\star}}^2$.
			Then, since $\frac{r^2 \eta^2}{\underline{\bar{\tau}}^2} = \frac{b^2\q^4 \bar{\p}^2}{1024} \leq 1$, $\frac{\eta^2(1-\underline{\p})}{\Lambda_4} \leq \frac{\eta^2}{\Lambda_4} = \frac{b\q^2}{64 - b\q^2} \leq 1$, and $\frac{\eta^2 (\Theta + \hat{\Theta})}{\kappa\Lambda_3} = 1$, we have
			\myeqn{
				\arraycolsep=0.2em
				\begin{array}{lcl}
					\Cc_{1,K} &\leq& \norms{u^0 - u^{\star}}^2 + 4\Cc_{0,K}  \leq 91\norms{u^0 - u^{\star}}^2, \vspace{1ex}\\
					\Cc_{2,K} &\leq& 26\Cc_{0,K} + 24\Cc_{1,K}
					\leq 2769\norms{u^0 - u^{\star}}^2,
				\end{array}
			}
			which completes the proof.
		\end{proof}
		
		\noindent\textbf{The resolvent cost of full refreshes.}
		In Algorithm~\ref{alg:A2}, each iteration only evaluates the block $j_k$ of the resolvent.
		However, each refresh of the L-SVRG snapshot requires all $m$ blocks of $Gw^k$ with $w^k = J_{\lambda\circ T}u^{k-1}$, and each refresh of the GJS table requires all blocks of $G_i$ at $J_{\lambda\circ T}u^k$ and $J_{\lambda\circ T}u^{k-1}$.
		Under the convention in Section~\ref{sec:ch2:OG_operator}, this costs one or two full evaluations of $J_{\lambda\circ T}$ per refresh.
		Counting one full evaluation as one unit and one block evaluation as $\frac{1}{m}$ unit, the expected resolvent cost over $K$ iterations is $\BigOs{1 + K(\frac{1}{m} + \q)}$.
		This additional cost does not arise in Algorithm~\ref{alg:A1}, which stores $x^k$ explicitly.
		
		For the mini-batch estimator, we fix $\zeta > 0$, define $C_{\zeta} := \sum_{k=0}^{\infty} \big[(k+2)[\log(k+2)]^{1+\zeta}\big]^{-1} < +\infty$, and use
		\myeq{eq:VACOG_batch_schedule}{
			\arraycolsep=0.2em
			\begin{array}{lcl}
				b_k &:=& \Big\lceil \frac{n}{1 + (n-1)/\big(\bar{\p}C_{\zeta}(k+2)[\log(k+2)]^{1+\zeta}\big)} \Big\rceil, \qquad k \geq 1.
			\end{array}
		}
		\begin{corollary}\label{co:VACOG_svrg_saga_complexity}
			For \eqref{eq:FNI}, suppose that $\zer{\Psi} \ne \emptyset$ and  Assumption \ref{as:VACOG} hold.
			Let $\sets{(u^k, x^k)}$ be generated by Algorithm~\ref{alg:A2} using $r := 3$ and $\eta := \frac{\bar{\eta}}{2}$.
			\vspace{0.5ex}
			\begin{compactitem}
			\itemsep=0.4em
			\item[\textbf{$\mathrm{(i)}$~$(\text{Mini-batch estimator})$}]
			For $\epsilon > 0$, if we use the batch sizes \eqref{eq:VACOG_batch_schedule}, then Algorithm~\ref{alg:A2} with the mini-batch estimator \eqref{eq:minibatch_estimator} requires at most
			$\Tc_{G_i} := n + \BigOlog{\min\sets{\frac{\norms{\bar{L}}_1 \norms{\bar{L}}_{\infty}}{m\epsilon^2}, \frac{n\norms{\bar{L}}_1}{m\epsilon}}}$ and $\Tc_{T} := \BigO{\frac{\norms{\bar{L}}_1}{m\epsilon}}$
			evaluations of $G_i$ and $J_{\lambda\circ T}$, respectively, to achieve $\Expk{\norms{Gx^K + \xi^K}^2} \leq \epsilon^2$.
			
			\item[\textbf{$\mathrm{(ii)}$~$($L-SVRG and Generalized JacSketch$)$}]
			Let $n \geq 2$ and $\epsilon > 0$.
			If we choose $\q := n^{-1/3}$ and $b := \lceil n^{2/3} \rceil$, then Algorithm~\ref{alg:A2} with \eqref{eq:loopless_SVRG} or \eqref{eq:JacSketch_estimator} requires at most
			$\Tc_{G_i} := \BigO{n + \frac{\norms{\bar{L}}_1 n^{2/3}}{\epsilon}}$ and $\Tc_{T} := \BigO{\big(\frac{1}{m} + n^{-1/3}\big)\frac{\norms{\bar{L}}_1}{\epsilon}}$
			evaluations of $G_i$ and $J_{\lambda \circ T}$, respectively, to achieve $\Expk{\norms{Gx^K + \xi^K}^2} \leq \epsilon^2$.
			In particular, $\Tc_T = \BigO{\frac{\norms{\bar{L}}_1}{m\epsilon}}$, whenever $m \leq n^{1/3}$.
			\end{compactitem}
		\end{corollary}
		
		\begin{proof}
			\noindent\textbf{$\mathrm{(i)}$~Mini-batch estimator.}
			For the mini-batch estimator \eqref{eq:minibatch_estimator}, we have $\kappa := 1$, $\Theta = 2$, $\hat{\Theta} = 0$, and $\delta_k = \frac{2(n-b_k)}{b_k(n-1)}$.
			Then, we have $\bar{\eta} = \frac{1}{6\norms{\bar{L}}_1}$.
			Let us choose $\eta := \frac{\bar{\eta}}{2} = \frac{1}{12\norms{\bar{L}}_1}$.
			Now, from \eqref{eq:th1_BigO_convergence_rates}, for a fixed tolerance $\epsilon > 0$, to guarantee $\Expk{\norms{G{x}^k + \xi^k}^2} \leq \epsilon^2$, we can impose $\frac{1152 \norms{\bar{L}}_1^2 (\Cc_{1,K} + \Cc_{2,K})}{(K+6)(K+3)} \leq \epsilon^2$.
			To get $\Cc_{1,K} + \Cc_{2,K} = \BigOs{1}$, by Lemma~\ref{lem:ch2:VACOG_bound_C}, it suffices to ensure $\sum_{k=0}^K \frac{n - b_k}{b_k(n-1)} \leq \frac{1}{\bar{\p}}$, which holds for the batch sizes \eqref{eq:VACOG_batch_schedule}, since they satisfy $b_k = \BigOs{k(\log k)^{1+\zeta}}$ and
			$\frac{n - b_k}{b_k(n-1)} \leq \frac{1}{\bar{\p}C_{\zeta} (k+2)[\log(k+2)]^{1+\zeta}}$ for all $k \geq 1$.
			Then, the number of iterations is $K = \mathcal{O}\big( {\frac{\norms{\bar{L}}_1}{\epsilon}} \big)$.
			Hence, the expected number of evaluations of $G_i$ is
			\myeqn{
				\arraycolsep=0.2em
				\begin{array}{lcl}
					\Expk{\Tc_{G_i}} & = & n + \frac{2}{m}\sum\limits_{k=0}^K b_k  = n + \mathcal{O}\Big({\min\set{\frac{\bar{\p}K^2[\log(K+2)]^{1+\zeta}}{m}, \frac{nK}{m}}}\Big)
					\vspace{1ex}\\
					&=& n + \widetilde{\mathcal{O}}\Big({\min\set{\frac{\norms{\bar{L}}_1 \norms{\bar{L}}_{\infty}}{m\epsilon^2}, \frac{n\norms{\bar{L}}_1}{m\epsilon}}}\Big).
				\end{array}
			}
			The expected number of evaluations of $J_{\lambda \circ T}$ is $\Expk{\Tc_{T}} = \frac{K}{m} = \BigO{\frac{\norms{\bar{L}}_1}{m\epsilon}}$.
			
			\vspace{0.5ex}
			\noindent\textbf{$\mathrm{(ii)}$~\ref{eq:loopless_SVRG} and \ref{eq:JacSketch_estimator}.}
			For the loopless SVRG estimator \eqref{eq:loopless_SVRG}, we have $\kappa := \frac{\q}{2}$, $\Theta = \hat{\Theta} := \frac{4}{b\q}$.
			Suppose that $1 \leq b\q^2 \leq 32$.
			Let us choose $\eta := \frac{\bar{\eta}}{2} = \frac{b\q^2}{96\norms{\bar{L}}_1}$.
			From Lemma~\ref{lem:ch2:VACOG_bound_C}, we can upper bound the quantities in \eqref{eq:co_VACOG_complexity_quantities} as $\Cc_{0,K} \leq 22.5\norms{u^0 - u^{\star}}^2$, $\Cc_{1,K} \leq 91\norms{u^0 - u^{\star}}^2$, and $\Cc_{2,K} \leq 2769\norms{u^0 - u^{\star}}^2$.
			Thus, the convergence rates in \eqref{eq:th1_BigO_convergence_rates} hold for \eqref{eq:loopless_SVRG} with these $\Cc_{1,K}$ and $\Cc_{2,K}$.
			
			Now, from \eqref{eq:th1_BigO_convergence_rates}, for a fixed tolerance $\epsilon > 0$, to guarantee $\Expk{\norms{G{x}^k + \xi^k}^2} \leq \epsilon^2$, we need to impose $\frac{210862080 \norms{\bar{L}}_1^2 \norms{u^0 - \uopt}^2}{b^2 \q^4 (K+6)(K+3)} \leq \epsilon^2$.
			This relation leads to $K = \BigO{\frac{\norms{\bar{L}}_1}{b\q^2 \epsilon}}$.
			Hence, the expected number of evaluations of $G_i$ is $\Expk{\Tc_{G_i}} = n + K \big(n\q + \frac{3b}{m}\big) = \BigO{n + \frac{\norms{\bar{L}}_1}{\epsilon}\left(\frac{n}{b\q} + \frac{1}{m\q^2}\right)}$.
			With $b := \lceil n^{2/3} \rceil$ and $\q := n^{-1/3}$, where $n \geq 2$, we have $1 \leq b\q^2 < 2$.
			Hence, Lemma~\ref{lem:ch2:VACOG_bound_C}(b) applies, and $\Expk{\Tc_{G_i}} = \BigO{n + \frac{\norms{\bar{L}}_1 n^{2/3}}{\epsilon}}$.
			By the discussion before this corollary, the expected number of evaluations of $J_{\lambda \circ T}$ is $\BigOs{1 + K(\frac{1}{m} + \q)} = \BigO{\big(\frac{1}{m} + n^{-1/3}\big)\frac{\norms{\bar{L}}_1}{\epsilon}}$.
			For \eqref{eq:JacSketch_estimator}, the expected component cost per iteration becomes $2n\q + \frac{2b}{m}$, which leads to the same orders.
		\end{proof}
		
		\begin{remark}[\textbf{Effect of block heterogeneity}]
			In a highly heterogeneous regime where $\norms{\bar{L}}_1 = \BigOs{\norms{\bar{L}}_{\infty}}$, the iteration complexity becomes $\BigOs{\epsilon^{-1}}$, up to fixed problem-dependent constants, and the \ref{eq:loopless_SVRG} and \ref{eq:JacSketch_estimator} oracle complexities reduce to $
			\BigOs{n + n^{2/3}\epsilon^{-1}}$.
			If instead $\bar{L}_1 = \dots = \bar{L}_m$, then $\norms{\bar{L}}_1 = m\norms{\bar{L}}_{\infty}$ and $\p_j = \frac{1}{m}$, so the iteration complexity is $\BigOs{m \epsilon^{-1}}$.
			For \ref{eq:loopless_SVRG} and \ref{eq:JacSketch_estimator}, the full-snapshot cost is not reduced by block sampling, leading to an oracle complexity of $\BigOs{m n^{2/3}\epsilon^{-1}}$.
			For the mini-batch estimator, however, the factor \(m\) in the iteration count is canceled by the $\frac{1}{m}$ cost of a selected block, yielding an oracle complexity of $\BigOs{\min\sets{\epsilon^{-2}, n\epsilon^{-1}}}$.
			Thus the nonuniform and uniform regimes generally do not have the same oracle complexity.
			Note that the iteration complexity $\BigOs{m \epsilon^{-1}}$ has been observed in \cite{Nesterov2012,fercoq2015accelerated,Alacaoglu2017} for Nesterov's accelerated randomized coordinate first-order methods.
			Here, we obtain similar results, but for a class of possibly nonmonotone problems covered by \eqref{eq:FNI}.
		\end{remark}
	}
	
	\beforesec
	\section{Numerical Experiments}\label{sec:ch2:num_experiments}
	\aftersec
	We consider a class of finite-sum quadratic minimax problems, adapted from the synthetic construction in \cite{TranDinh2024}, to illustrate the performance of our variance-reduced randomized block-coordinate methods. In particular, we examine the effects of the block partition and variance-reduction strategy on both non-accelerated and accelerated schemes.
	All algorithms and baselines are implemented in Python and executed on UNC’s Longleaf HPC cluster with NVIDIA A100 (40 GB) and L40 (48 GB) GPUs.
	
	\vspace{0.5ex}
	\noindent\textbf{$\mathrm{(a)}$~Mathematical model.}
	Consider the regularized saddle-point problem
	\myeq{eq:num_minimax}{
		\arraycolsep=0.2em
		\begin{array}{lcl}
			\min\limits_{u \in \R^{p_1}} \max\limits_{v \in \R^{p_2}}
			\set{f(u) + \frac{1}{n}\sum_{i = 1}^n \Hc_i(u,v) - g(v)},
		\end{array}
	}
	where
	$\Hc_i(u, v) := \frac{1}{2} u^\top A_i u + u^\top L_i v - \frac{1}{2} v^\top B_i v + b_i^\top u - c_i^\top v$,
	with $A_i\in\mathbb{R}^{p_1\times p_1}$, $B_i\in\mathbb{R}^{p_2\times p_2}$, $L_i\in\mathbb{R}^{p_1\times p_2}$, $b_i\in\mathbb{R}^{p_1}$, and $c_i\in\mathbb{R}^{p_2}$.
	The functions $f$ and $g$ are usually nonsmooth, used to characterize the set of constraints or regularizers.
	
	Let $x=(u,v) \in \R^p$, where $p := p_1 + p_2$.
	Then, the first-order optimality condition of \eqref{eq:num_minimax} can be written as the composite inclusion $0 \in \frac{1}{n} \sum_{i=1}^n G_ix + Tx$, a special case of \eqref{eq:FNI}, where $G_ix:=
	\begin{bmatrix}
		A_i u + L_i v + b_i\\
		-L_i^\top u + B_i v + c_i
	\end{bmatrix}$,
	$Tx:=
	\begin{bmatrix}
		\partial f(u)\\
		\partial g(v)
	\end{bmatrix}$.
	We can see that $G$ has a finite-sum structure, and each $G_i$ is affine and its restriction to any coordinate block can be evaluated without forming the full operator.
	
	\vspace{0.5ex}
	\noindent\textbf{$\mathrm{(b)}$~Data generation and experimental settings.}
	Following the construction in \cite{TranDinh2024}, we generate $
	A_i := Q_i D_i Q_i^\top$, where $Q_i$ is an orthogonal matrix obtained from a QR factorization of an independent standard Gaussian matrix, and $D_i$ is a diagonal matrix whose entries are independently sampled from the standard normal distribution $\Nc(0,1)$ and then clipped as $d_{ij} := \max\sets{d_{ij}, \underline{d}}$, with $\underline{d} = 0.1$.
	The matrices $B_i$ are sampled in the same way as $A_i$.
	The entries of $L_i$, $b_i$, and $c_i$ are independently drawn from $\Nc(0,1)$.
	This ensures that each $G_i$ is co-coercive.
	
	We perform two sets of experiments:
	\begin{compactitem}
		\item \textit{Experiment 1}: 
		The functions $f(u)$ and $g(v)$ are chosen as the $\ell_1$-norm regularizer, i.e., $f(u) = \varrho\norms{u}_1$ and $g(v) = \varrho\norms{v}_1$, promoting sparsity of the variables.
		In this case, the operator $T$ in \eqref{eq:FNI} is maximally monotone.
		\item \textit{Experiment 2}: The functions $f(u)$ and $g(v)$ are chosen as the nonconvex SCAD (smoothly clipped absolute deviation) regularizer \cite{fan2001variable}.
		In this case, the operator $T$ in \eqref{eq:FNI} is nonmonotone.
		In fact, $T$ is locally $\rho$-co-hypomonotone with $\rho \leq a - 1 = 2.7$.
		Similar to \cite{TranDinh2025a}, we can still apply our methods to solve \eqref{eq:num_minimax}.
	\end{compactitem}
	We choose $n = 5000$ and $p = 200$ ($p_1 = 133$, $p_2 = 67$).
	For the regularizers, we choose $\varrho = 0.04$, $\rho = 0.04$, and $a = 3.7$.
	For each experiment, we generate five independent problem instances and report the mean relative forward-backward residual $\frac{\norms{\Fc_{\lambda}x^k}}{\norms{\Fc_{\lambda}x^0}}$ over five independent runs, where $x^0 := \mbf{1} \in \R^p$ and $\lambda := \frac{1}{L}$, with $L$ being the Lipschitz constant of $G$.
	
	\vspace{0.5ex}
	\noindent\textbf{$\mathrm{(c)}$~Algorithms and parameters.}
	First, we compare \texttt{VCOG} and \texttt{VCFOG} using three estimators: increasing mini-batching (\texttt{MB}), loopless SVRG (\texttt{SVRG}), and generalized JacSketch (\texttt{GJS}).
	For each estimator, we consider $m \in \sets{10, 30, 50}$ coordinate blocks, sampled uniformly.
	The same parameter settings are used across all three block choices within each experiment.
	In particular, for \texttt{VCOG}, we choose $\gamma := 0.97$ for all variants, and $\eta := \frac{0.045}{L}$ in \textit{Experiment 1} and $\eta := \frac{0.05}{L}$ in \textit{Experiment 2}.
	For \texttt{VCFOG}, we choose $\eta := \frac{0.004}{L}$ in \textit{Experiment 1} and $\eta := \frac{0.0065}{L}$ in \textit{Experiment 2}, $r := 50$, and $\lambda := \frac{1}{L}$.
	The mini-batch size and refresh probability are chosen as follows:
	\begin{compactitem}
		\item For \texttt{MB}, we choose $b_k := \min\sets{n, \left\lceil 32 + 0.002(k+1)^{1.1} \right\rceil}$.
		\item For \texttt{SVRG} and \texttt{GJS}, we choose  $b := \left\lfloor 0.5n^{2/3} \right\rfloor$ and $\q := 0.05n^{-1/3}$.
	\end{compactitem}
	Next, we compare one representative (\texttt{SVRG}) of \texttt{VCOG} and \texttt{VCFOG} with three other methods: the variance-reduced forward-reflected-backward splitting (\texttt{VFRBS}) \cite{alacaoglu2022beyond}, variance-reduced extragradient (\texttt{VEG}) \cite{alacaoglu2021stochastic}, and the randomized block-coordinate primal-dual algorithm (\texttt{RB-PDA}) \cite{hamedani2026randomized}.
	The parameters are chosen as follows:
	\begin{compactitem}
		\item For \texttt{VCOG-SVRG}, we choose $\eta := \frac{0.04}{L}$ in \textit{Experiment 1} and $\eta := \frac{0.05}{L}$ in \textit{Experiment 2}, fixed $m := 50$ blocks, $b := \lfloor 0.25n^{2/3}\rfloor$, and $\q := 0.05n^{-1/3}$.
		\item For \texttt{VCFOG-SVRG}, we choose $\eta := \frac{0.004}{L}$ in \textit{Experiment 1} and $\eta := \frac{0.006}{L}$ in \textit{Experiment 2}, fixed $m := 50$ blocks, $b := \lfloor 0.25n^{2/3}\rfloor$, and $\q := 0.05n^{-1/3}$.
		\item For \texttt{VFRBS}, we use $\eta := \frac{0.95(1-\sqrt{1-\q})}{2L}$, $b := \lfloor 0.25n^{2/3}\rfloor$, and $\q := 0.05n^{-1/3}$.
		\item For \texttt{VEG}, we choose $b := \lfloor 0.25n^{2/3}\rfloor$,
		$\q := 0.05n^{-1/3}$, $\alpha := 1-\q$, and
		$\eta := \frac{0.95\sqrt{1 - \alpha}}{L}$.
		\item For \texttt{RB-PDA}, we use $M := 33$ primal and $N := 17$ dual blocks (leading to $m = 50$ blocks in total), $\theta := 0.05$, primal stepsize $\tau := \frac{0.75}{ML}$,
		and dual stepsize $\sigma := \frac{0.75}{NL}$.
	\end{compactitem}
	
\begin{figure}[!ht]
		\centering
		\includegraphics[width=\linewidth]{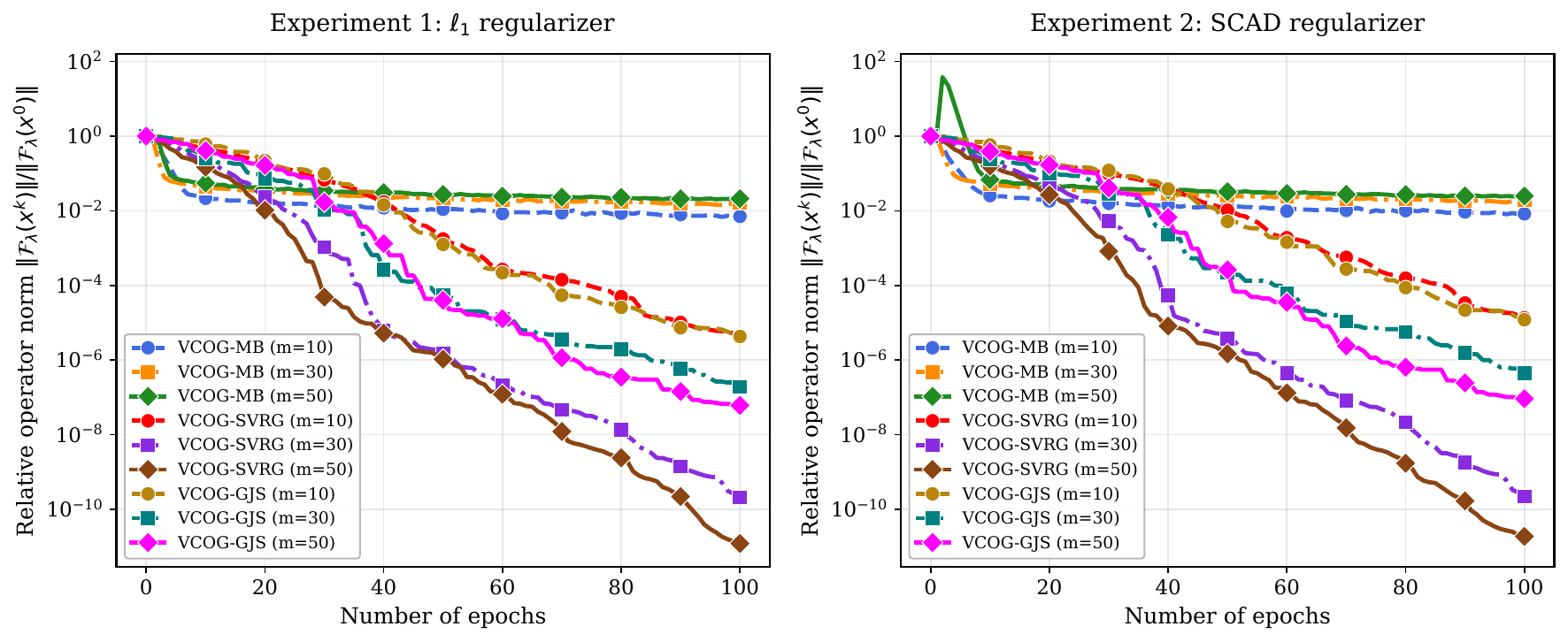}
		\vspace{-5ex}
		\caption{Performance of \texttt{VCOG} with three estimators for \eqref{eq:num_minimax} with $\ell_1$ and SCAD regularizers (averaged over 5 problem instances).
		}
		\label{fig:num_vcog}
	\end{figure}

\begin{figure}[!ht]
		\vspace{-2ex}
		\centering
		\includegraphics[width=\linewidth]{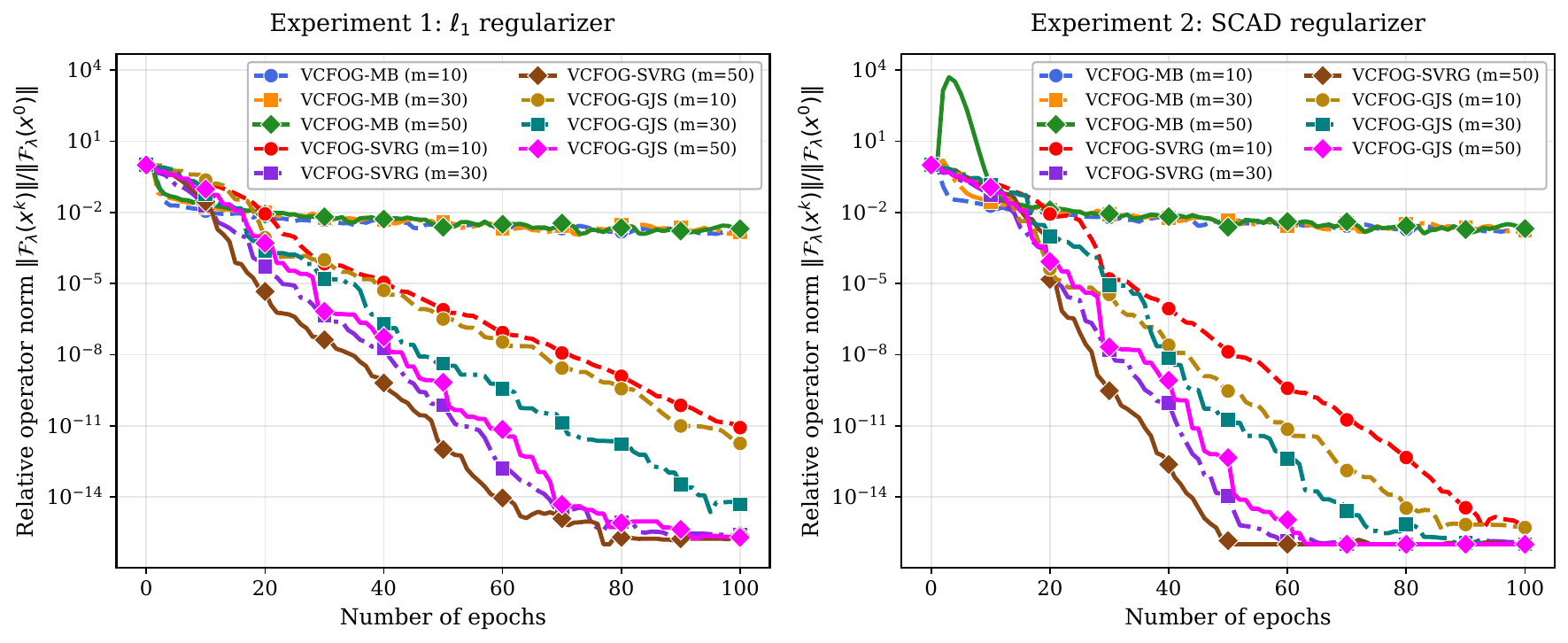}
		\vspace{-5ex}
		\caption{Performance of \texttt{VCFOG} with three estimators for \eqref{eq:num_minimax} with $\ell_1$ and SCAD regularizers (averaged over 5 problem instances).}
		\label{fig:num_vcfog}
	\end{figure}

	\vspace{0.5ex}
	\noindent\textbf{$\mathrm{(d)}$~Numerical results.}
	Figures~\ref{fig:num_vcog} and \ref{fig:num_vcfog} present the results after 100 epochs for several variants of \texttt{VCOG} and \texttt{VCFOG}, respectively.
	Here, the number of epochs represents the number of passes over all components $G_i$ and block-coordinates.
	
	First, the variance-reduced estimators using control variates substantially outperform increasing mini-batching in both experiments.
	While the \texttt{MB} curves decrease rapidly at the beginning, they subsequently exhibit a pronounced slowdown and remain around $10^{-2}$ within the prescribed epoch budget.
	In contrast, \texttt{L-SVRG} and \texttt{GJS} continue decreasing over several additional orders of magnitude.
	This behavior is observed for both the convex $\ell_1$ and nonconvex SCAD regularizers.
	
	Second, increasing the number of coordinate blocks considerably improves the performance of the variance-reduced schemes.
	For both \texttt{L-SVRG} and \texttt{GJS}, the choice $m = 50$ consistently yields faster convergence than $m = 30$ and $m = 10$.
	The effect is especially pronounced for \texttt{L-SVRG}, which provides the fastest convergence among the tested methods.
	However, we do not observe a similar trend in the \texttt{MB} variants.
	
	Third, \texttt{VCFOG} exhibits a substantial improvement over its non-accelerated counterpart under the tested settings.
	With $m = 50$, the accelerated variance-reduced methods reach residuals close to machine precision within 80 epochs in both experiments, whereas \texttt{VCOG} decreases more gradually.
	The SCAD results also show that the accelerated schemes remain effective with a nonmonotone regularization operator.
	
\begin{figure}[!ht]
		\centering
		\includegraphics[width=\linewidth]{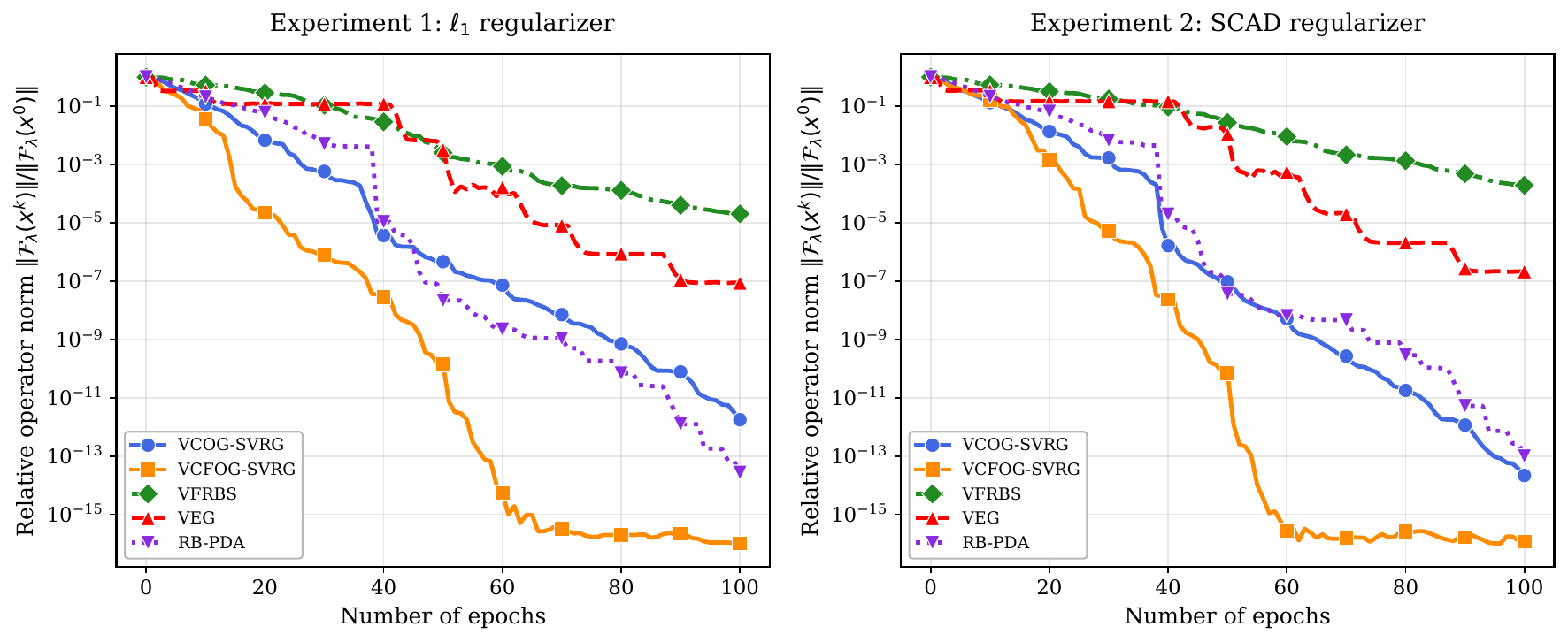}
		\vspace{-5ex}
		\caption{Performance of five different methods for \eqref{eq:num_minimax} with $\ell_1$ and SCAD regularizers (averaged over 5 problem instances).}
		\label{fig:num_compare}
	\end{figure}

	Finally, Figure~\ref{fig:num_compare} presents the numerical results for five different methods for solving \eqref{eq:num_minimax}.
	The comparison results show that \texttt{VCFOG-SVRG} consistently achieves the fastest convergence in both experiments, reaching residuals near machine precision within approximately 65 epochs.
	\texttt{VCOG-SVRG} and \texttt{RB-PDA} also exhibit strong convergence, attaining residuals of order $10^{-12}$--$10^{-14}$ after 100 epochs, although both converge more slowly than \texttt{VCFOG-SVRG}.
	In contrast, \texttt{VFRBS} and \texttt{VEG} exhibit substantially slower residual decay, particularly under SCAD regularization.
	Overall, these experiments demonstrate the complementary benefits of randomized coordinate updates, variance reduction, and acceleration within the proposed frameworks \texttt{VCOG} and \texttt{VCFOG}.
	
	\beforesec
	\section*{Acknowledgments}
	\aftersec
	This work was partially supported by NSF RTG grant DMS-2134107 and ONR grant N00014-23-1-2588.
	
	\appendix
	\renewcommand{\thesubsection}{\thesection.\arabic{subsection}}
	\renewcommand{\thesubsubsection}{\thesubsection.\arabic{subsubsection}}
	
	\beforesec
	\section{Technical Lemmas}\label{apdx:sec:ch2:appendix1}
	\aftersec
	We recall here some technical lemmas necessary to prove our results in the main text.
	
	\begin{lemma}[\cite{Bauschke2011}]\label{lem:ch2:tech_convergence1}
		Let $\sets{ a_k }$, $\sets{ b_k }$, and $\sets{ d_k }$ be sequences of real numbers.
		Assume that $\sets{ a_k }$ is bounded from below, $\sets{b_k }$ and $\sets{ d_k }$ are nonnegative sequences such that $\sum_{k=0}^{\infty} d_k < \infty$, and
		\myeqn{
			\arraycolsep=0.2em
			\begin{array}{lcl}
				a_{k+1} &\leq& a_k - b_k + d_k, \quad \forall k \geq 0.
			\end{array}
		}
		Then, the following statements hold:
		\begin{compactitem}
			\item[$\mathrm{(a)}$] $\sets{b_k}$ is summable, namely $\sum_{k=0}^{\infty} b_k < \infty$.
			\item[$\mathrm{(b)}$] $\sets{ a_k }$ is convergent.
		\end{compactitem}
	\end{lemma}
	
	\begin{proof}
		Let $D_k := \sum_{l=k}^{\infty} d_l$ and $c_k := a_k + D_k$.
		Then $c_{k+1} \leq c_k - b_k$, and $\sets{c_k}$ is bounded from below because $\sets{a_k}$ is bounded from below and $D_k \geq 0$.
		Hence, $\sets{c_k}$ converges, and summing $b_k \leq c_k - c_{k+1}$ over $k$ gives $\sum_{k=0}^{\infty} b_k < +\infty$.
		Since $D_k \to 0$, the sequence $a_k = c_k - D_k$ also converges.
	\end{proof}
	
	\begin{lemma}[\cite{TranDinh2025a}]\label{lem:ch2:tech_convergence2}
		Let $\sets{ a_k }_{k\geq 0}$ be a sequence of nonnegative numbers and $\alpha > 0$ be a constant.
		Assume that $\lim_{k \to \infty} k^{\alpha + 1} a_k$ exists and $\sum_{k=0}^\infty k^\alpha a_k < \infty$. Then, $\lim_{k\to \infty} k^{\alpha + 1} a_k = 0$.
	\end{lemma}
	
	\begin{proof}
		Let $\ell := \lim_{k\to\infty} k^{\alpha+1}a_k \geq 0$.
		If $\ell > 0$, then $k^{\alpha}a_k \geq \frac{\ell}{2k}$ for all sufficiently large $k$, which contradicts $\sum_{k=0}^{\infty} k^{\alpha}a_k < \infty$.
		Therefore, $\ell = 0$.
	\end{proof}
	
	\begin{lemma}[Supermartingale theorem \citep{robbins1971convergence}]\label{lem:ch2:supermartingale}
		Let $\{X_k\}$, $\{ \alpha_k \}$, $\{V_k\}$, and $\{R_k\}$ be sequences of nonnegative integrable random variables on some arbitrary probability space and adapted to the filtration $\{\mathcal F_k\}$ with $\sum_{k=0}^\infty \alpha_k < +\infty$ and $\sum_{k=0}^\infty R_k < +\infty$ almost surely, and
		\myeqn{
			\arraycolsep=0.2em
			\begin{array}{lcl}
				\Exp{X_{k+1} \mid \mathcal F_k} &\leq& (1 + \alpha_k) X_k - V_k + R_k, \quad \forall k \geq 0 \text{ almost surely.}
			\end{array}
		}
		Then, $\{X_k\}$ almost surely converges and $\sum_{k=0}^\infty V_k < +\infty$ almost surely.
	\end{lemma}
	
	\begin{lemma}[\cite{tran2025vfog}]\label{lem:ch2:Opial}
		Let $\Psi: \R^p \rightrightarrows \R^p$ have a closed graph and a nonempty zero set.
		Let $\sets{x^k}$ be a sequence of random vectors such that, for every $x^{\star} \in \zer{\Psi}$, the sequence $\sets{\norms{x^k - x^{\star}}}$ converges almost surely.
		Assume further that there exist random vectors $w^k$ such that $(x^k, w^k) \in \gra{\Psi}$ and $\norms{w^k} \to 0$ almost surely.
		Then, $\sets{x^k}$ converges almost surely to a $\zer{\Psi}$-valued random variable.
	\end{lemma}
	
	\beforesec
	\section{The Proofs of Lemmas in Section~\ref{sec:ch2:OG_operator}}\label{apdx:sec:ch2:OG_operator}
	\aftersec
	\noindent\textbf{\textit{Proof of Lemma~\ref{le:imb_estimator}.}}
	First, since $\Bc_k$ is sampled uniformly without replacement and independently of $\Fc_k$, we have $\Expsk{k}{\widetilde{S}^k} = S^k$.
	
	Second, let $\hat{e}_i^k := G_i x^k - Gx^k$ and $d_i^k := \hat{e}_i^k - \gamma^k \circ \hat{e}_i^{k-1}$, $i \in [n]$.
	Then, $\frac{1}{n} \sum_{i=1}^n d_i^k = 0$, and thus the finite-population variance identity \cite[Theorem 2.2]{cochran1977sampling} yields
	\myeqn{
		\arraycolsep=0.2em
		\begin{array}{lcl}
			\Expsk{k}{\norms{\widetilde{S}^k - S^k}_{\beta}^2} &=& \Expsk{k}{\norms{\frac{1}{b_k}\sum_{i \in \Bc_k} d_i^k}_{\beta}^2} = \vartheta_k \frac{1}{n}\sum_{i = 1}^n \norms{d_i^k}_{\beta}^2, \qquad \vartheta_k := \frac{n - b_k}{b_k(n-1)}.
		\end{array}
	}
	Third, from $d_i^k = (\hat{e}_i^k - \hat{e}_i^{k-1}) + (1-\gamma^k) \circ \hat{e}_i^{k-1}$, using Young's inequality, we have
	$\norms{d_i^k}_{\beta}^2 \leq 2\norms{\hat{e}_i^k - \hat{e}_i^{k-1}}_{\beta}^2 + 2\norms{1-\gamma^k}_{\infty}^2\norms{\hat{e}_i^{k-1}}_{\beta}^2$.
	Moreover,
	\myeqn{
		\arraycolsep=0.2em
		\begin{array}{lcl}
			\frac{1}{n}\sum_{i = 1}^n \norms{\hat{e}_i^k - \hat{e}_i^{k-1}}_{\beta}^2 &=& \frac{1}{n} \sum_{i = 1}^n \norms{(G_ix^k - G_ix^{k-1}) - (Gx^k - Gx^{k-1})}_{\beta}^2 \vspace{1ex}\\
			&=& \frac{1}{n} \sum_{i = 1}^n \norms{G_ix^k - G_ix^{k-1}}_{\beta}^2 - \norms{Gx^k - Gx^{k-1}}_{\beta}^2 \vspace{1ex}\\
			&\leq& \frac{1}{n} \sum_{i = 1}^n \norms{G_ix^k - G_ix^{k-1}}_{\beta}^2 = U_k.
		\end{array}
	}
	Combining the last three relations, $V_k = \frac{1}{n} \sum_{i = 1}^n \norms{\hat{e}_i^k}_{\beta}^2$, and $\vartheta_k \leq 1$, we obtain
	\myeqn{
		\arraycolsep=0.2em
		\begin{array}{lcl}
			\Expsk{k}{\norms{\widetilde{S}^k - S^k}_{\beta}^2} &\leq& 2\vartheta_k U_k + 2\norms{1-\gamma^k}_{\infty}^2 \vartheta_k V_{k-1} \leq 2U_k + 2\norms{1-\gamma^k}_{\infty}^2 \vartheta_k V_{k-1}.
		\end{array}
	}
	Therefore, if we set $\Delta_k := 2\vartheta_k U_k + 2\vartheta_k\norms{1-\gamma^k}_{\infty}^2 V_{k-1}$ for $k \geq 1$ and $\Delta_0 := 0$ (the initial estimator is exact), then the mini-batch estimator \eqref{eq:minibatch_estimator} satisfies Definition~\ref{def:VR_estimator} with $\kappa = 1$, $\Theta = 2$, $\hat{\Theta} = 0$, and $\delta_k = 2\vartheta_k$.
	\Eproof
	
	\vspace{0.75ex}
	\noindent\textbf{\textit{Proof of Lemma~\ref{lem:ch2:SVRG_estimator}.}}
	For $k \geq 1$, we define
	\myeqn{
		\arraycolsep=0.2em
		\begin{array}{lcl}
			a_i^k &:=& G_ix^k - \gamma^k \circ G_ix^{k-1} - (1 - \gamma^k) \circ G_iw^k \quad\text{and}\quad \Delta_k := \frac{1}{bn}\sum_{i=1}^n \norms{a_i^k}_{\beta}^2 .
		\end{array}
	}
	Then, $\widetilde{S}^k - S^k = \frac{1}{b}\sum_{i \in \Bc_k} a_i^k - \frac{1}{n}\sum_{i=1}^n a_i^k$.
	Conditioned on $\Fc_k$ and on the refreshed snapshot $w^k$, the mini-batch $\Bc_k$ is still uniform.
	Hence, the variance bound for sampling without replacement \cite[Theorem 2.2]{cochran1977sampling} gives
	\myeqn{
		\arraycolsep=0.2em
		\begin{array}{lcl}
			\Expsk{k}{\widetilde{S}^k \mid w^k} &=& S^k \quad\text{and}\quad \Expsk{k}{\norms{\widetilde{S}^k - S^k}_{\beta}^2 \mid w^k} \leq \Delta_k .
		\end{array}
	}
	Taking the conditional expectation with respect to the snapshot proves the first two conditions of Definition~\ref{def:VR_estimator}.
	
	For the third condition, we consider two cases.
	If the snapshot is refreshed, i.e., $w^k = x^{k-1}$, then $a_i^k = G_ix^k - G_ix^{k-1}$.
	Otherwise, $w^k = w^{k-1}$, and the identity $1 - \gamma^k = \nu^k \circ (1 - \gamma^{k-1})$ gives
	\myeqn{
		\arraycolsep=0.2em
		\begin{array}{lcl}
			a_i^k &=& \nu^k \circ a_i^{k-1} + (G_ix^k - G_ix^{k-1}) - \nu^k\gamma^{k-1} \circ (G_ix^{k-1} - G_ix^{k-2}).
		\end{array}
	}
	We apply Young's inequality with parameter $\frac{\q}{2(1-\q)}$ to this decomposition and use $\norms{\nu^k\gamma^{k-1} \circ v}_{\beta} \leq \norms{\nu^k}_{\infty}\norms{v}_{\beta}$.
	Averaging over the refresh event, we obtain
	\myeqn{
		\arraycolsep=0.2em
		\begin{array}{lcl}
			\Expsk{k}{\Delta_k}
			&\leq& \big(1 - \frac{\q}{2}\big)\norms{\nu^k}_{\infty}^2\Delta_{k-1} + \frac{\q}{b}U_k + \frac{2(1-\q)(2-\q)}{b\q}\big(U_k + \norms{\nu^k}_{\infty}^2U_{k-1}\big) \vspace{1ex}\\
			&\leq& \big(1 - \frac{\q}{2}\big)\norms{\nu^k}_{\infty}^2\Delta_{k-1} + \frac{4}{b\q}U_k + \frac{4}{b\q}\norms{\nu^k}_{\infty}^2U_{k-1},
		\end{array}
	}
	where the last inequality follows from $\q^2 + 2(1-\q)(2-\q) = 4 - 6\q + 3\q^2 \leq 4$ for $\q \in (0,1)$.
	This proves the claim with $\kappa = \frac{\q}{2}$, $\Theta = \hat{\Theta} = \frac{4}{b\q}$, and $\delta_k = 0$.
	\Eproof
	
	\vspace{0.75ex}
	\noindent\textbf{\textit{Proof of Lemma~\ref{le:GJS_estimator}.}}
	First, conditioning on the table $\{H_i^k\}_{i=1}^n$ and using the uniform sampling of $\mathcal{B}_k$, we have
	\myeqn{
		\arraycolsep=0.2em
		\begin{array}{lcl}
			\Expsk{k}{\widetilde S^k\mid \sets{H_i^k}_{i=1}^n}
			&=&	\bar{H}^k + \frac{1}{n} \sum_{i=1}^n (S_i^k-H_i^k) = S^k.
		\end{array}
	}
	Hence, taking the conditional expectation $\Expsk{k}{\cdot}$ of both sides, we obtain $\Expsk{k}{\widetilde S^k-S^k} = 0$, proving the first condition of \eqref{eq:VR_estimator}.
	Next, let $d_i^k := S_i^k - H_i^k$, $\bar{d}^k := \frac{1}{n} \sum_{i=1}^n d_i^k$.
	Then $\widetilde{S}^k-S^k = \frac{1}{b} \sum_{i\in\Bc_k} d_i^k - \bar{d}^k$.
	Using \cite[Theorem 2.2]{cochran1977sampling}, we have
	\myeqn{
		\arraycolsep=0.2em
		\begin{array}{lcl}
			\Expsk{k}{\norms{\widetilde{S}^k-S^k}_{\beta}^2 \mid \sets{H_i^k}_{i=1}^n}
			&=&
			\frac{n-b}{b(n-1)}
			\frac{1}{n} \sum_{i=1}^n
			\norms{d_i^k - \bar{d}^k}_{\beta}^2
			\leq
			\frac{1}{bn}\sum_{i=1}^n\norms{d_i^k}_{\beta}^2
			:= \Delta_k.
		\end{array}
	}
	Consequently, $\Expsk{k}{\norms{\widetilde{S}^k-S^k}_{\beta}^2} \leq \Expsk{k}{\Delta_k}$, proving the second condition of \eqref{eq:VR_estimator}.
	
	Finally, with probability $\q$, $H_i^k = S_i^k$ for all $i\in[n]$, and therefore $\Delta_k=0$.
	With probability $1 - \q$, $H_i^k = \nu^k \circ H_i^{k-1}$ and thus $S_i^k - H_i^k = \nu^k \circ (S_i^{k-1} - H_i^{k-1}) + S_i^k - \nu^k \circ S_i^{k-1}$.
	Then, we obtain
	\myeqn{
		\arraycolsep=0.2em
		\begin{array}{lcl}
			\Expsk{k}{\Delta_k} &=& \frac{1-\q}{bn}\sum_{i=1}^n
			\norms{\nu^k \circ (S_i^{k-1} - H_i^{k-1})	+ S_i^k - \nu^k \circ S_i^{k-1}}_{\beta}^2 \vspace{0ex}\\
			
			&\stackrel{\tiny\textcircled{1}}{\leq}& \textls[-70]{$\frac{1-\q}{bn}\big[1 + \frac{\q}{2(1-\q)}\big]\norms{\nu^k}_{\infty}^2 \sum\limits_{i = 1}^n \norms{S_i^{k-1} - H_i^{k-1}}_{\beta}^2  + \frac{1-\q}{bn}\big[1 + \frac{2(1-\q)}{\q}\big]\sum\limits_{i = 1}^n \norms{S_i^k - \nu^k \circ S_i^{k-1}}_{\beta}^2$}
			\vspace{0ex}\\
			
			&\stackrel{\tiny\textcircled{2}}{\leq}& \textls[-50]{$\left(1-\frac{\q}{2}\right) \norms{\nu^k}_{\infty}^2 \Delta_{k-1} + \frac{2}{b\q}
				\frac{1}{n}\sum\limits_{i=1}^n
				\norms{(G_i x^k - G_i x^{k-1})	- \nu^k \gamma^{k-1} \circ  (G_i x^{k-1} - G_i x^{k-2})}_{\beta}^2$} \vspace{0ex}\\
			
			&\stackrel{\tiny\textcircled{3}}{\leq}&\left(1-\frac{\q}{2}\right)\norms{\nu^k}_{\infty}^2 \Delta_{k-1} + \frac{4}{b\q}U_k + \frac{4}{b\q} \norms{\nu^k}_{\infty}^2 U_{k-1},
		\end{array}
	}
	where we have used Young's inequality in {\tiny\textcircled{1}} and {\tiny\textcircled{3}}, $(1-\q)\big[1 + \frac{2(1-\q)}{\q}\big] = \frac{(1 - \q)(2 - \q)}{\q} \leq\frac{2}{\q}$ in {\tiny\textcircled{2}}, and $\gamma_j^k \leq 1$ in {\tiny\textcircled{3}}.
	Thus, the generalized JacSketch-type estimator \eqref{eq:JacSketch_estimator} satisfies Definition~\ref{def:VR_estimator} with
	$\kappa=\frac{\q}{2}$, $\Theta = \widehat\Theta = \frac{4}{b\q}$, and $\delta_k=0$.
	\Eproof
	
	We conclude this section with the following lemma, which is used in Section~\ref{sec:ch2:VACOG}.
	\begin{lemma}[Common-shift compatibility of control-variate estimators]\label{le:control_variate_compability}
		Consider a control-variate construction rule $\mathfrak{C}$ which, for component directions $\sets{R_i^k}_{i=1}^n$, generates $\widetilde{R}^k = \bar{H}^k + R^k_{\Bc_k} - H_{\Bc_k}^k$, where $\bar{H}^k := \frac{1}{n} \sum_{i = 1}^n H_i^k$, $R_{\Bc_k}^k := \frac{1}{|\Bc_k|} \sum_{i \in \Bc_k} R_i^k$, and $H_{\Bc_k}^k := \frac{1}{|\Bc_k|} \sum_{i \in \Bc_k} H_i^k$.
		Assume that $\mathfrak{C}$ has the following properties.
		\begin{itemize}
			\item[$\mathrm{(i)}$] \textbf{$($Variance-reduction admissibility$)$}
			When $\mathfrak{C}$ is applied to an arbitrary finite-sum family $\sets{F_i}_{i=1}^n$, the resulting estimator satisfies Definition~\ref{def:VR_estimator} with $U_k$ and $V_k$ defined in terms of $\sets{F_i}_{i=1}^n$.
			
			\item [$\mathrm{(ii)}$] \textbf{$($Common-shift compatibility$)$}
			If another family $\sets{\Rc_i^k}_{i=1}^n$ satisfies $\Rc_i^k = R_i^k + d^k$ for all $i \in [n]$, where $d^k$ is independent of $i$, then under the same realization of the sampling and internal randomness, the corresponding control variates generated by $\mathfrak{C}$ satisfy $\Hc_i^k = H_i^k + c^k$ for some $c^k$ independent of $i$, $i \in [n]$.
		\end{itemize}
		In addition, assume that the shifted directions $\sets{\Rc_i^k}_{i=1}^n$ are the components of a finite-sum family to which property $\mathrm{(i)}$ applies.
		Then, the estimator $\widehat{\Rc}^k$ generated by $\mathfrak{C}$ from $\sets{\Rc_i^k}_{i=1}^n$ satisfies Definition~\ref{def:VR_estimator} with $U_k$ and $V_k$ defined from this shifted family.
		Moreover, $\widehat{\Rc}^k$ is pathwise identical to $\widetilde{\Rc}^k := \widetilde{R}^k + d^k$.
	\end{lemma}
	
	\begin{proof}
		Since $\widehat{\Rc}^k$ is generated by the construction rule $\mathfrak{C}$, property $\mathrm{(i)}$ ensures that $\widehat{\Rc}^k$ is an estimator of $\Rc^k$ that satisfies Definition~\ref{def:VR_estimator} with $U_k$ and $V_k$ defined from the shifted family.
		
		Next, from $\Rc_i^k = R_i^k + d^k$ and $d^k$ is independent of $i$, we have $\Rc_{\Bc_k}^k = R_{\Bc_k}^k + d^k$.
		Moreover, from $\Hc_i^k = H_i^k + c^k$ and $c^k$ is independent of $i$, we also have $\bar{\Hc}^k = \bar{H}^k + c^k$ and ${\Hc}_{\Bc_k}^k = H_{\Bc_k}^k + c^k$.
		Hence,
		\myeqn{
			\arraycolsep=0.2em
			\begin{array}{lcl}
				\widehat{\Rc}^k &:=& \bar{\Hc}^k + \Rc_{\Bc_k}^k - \Hc_{\Bc_k}^k = (\bar{H}^k + c^k) + (R_{\Bc_k}^k + d^k) - (H_{\Bc_k}^k + c^k) = \widetilde{R}^k + d^k =: \widetilde{\Rc}^k,
			\end{array}
		}
		which completes the proof.
	\end{proof}
	\begin{remark}
		The mini-batch, \eqref{eq:minibatch_estimator}, \eqref{eq:loopless_SVRG}, and \eqref{eq:JacSketch_estimator} estimators satisfy the common-shift compatibility condition in Lemma~\ref{le:control_variate_compability}. For the mini-batch estimator, there is no control variate, so we can take $c^k = 0$.
		If $\chi_k \sim \texttt{Bernoulli}(\q)$ denotes the refresh indicator, then $c^0 := d^0 = (1-\gamma_0) A_{\lambda\circ T} u^0$, and for \ref{eq:loopless_SVRG}, we have $c_k = \chi_k (1 - \gamma_k) A_{\lambda \circ T}u^{k-1} + (1 - \chi_k) \nu_k c_{k-1}$, and for \ref{eq:JacSketch_estimator}, we have $c_k = \chi_k d^k + (1 - \chi_k) \nu_k c_{k-1}$.
	\end{remark}
	
	\beforesubsec
	\subsection{The proofs of technical lemmas in Section~\ref{sec:ch2:VCOG_alg}}\label{apdx:VCOG_convergence}
	\aftersubsec
\textbf{\textit{Proof of Lemma~\ref{lem:ch2:VCOG_key_lemma1}.}}
First, since $\widetilde{S}^k = S^k + e^k$, the notations in \eqref{eq:VCOG_w_notations} give $S^k + (1-\gamma) \circ \xi^k = w^k - \gamma \circ \hat{w}^k$.
Hence, the first line of \eqref{eq:VCOG4} can be written as
\myeq{eq:VCOG_keylemma_proof1}{
	\arraycolsep=0.2em
	\begin{array}{lcl}
		y^{k+1} - y^k &=& -\frac{\eta_{j_k}}{\p_{j_k}}\big(w^k - \gamma \circ \hat{w}^k + e^k\big)_{[j_k]}.
	\end{array}
}
Since $y^k = x^k + \gamma\eta \circ \xi^k$, adding $\gamma\eta \circ Gx^k - \xopt$ to both sides gives $x^{k+1} - \xopt + \gamma\eta \circ \hat{w}^{k+1} = x^k - \xopt + \gamma\eta \circ w^k + (y^{k+1} - y^k)$.
Since $j_k$ is independent of the estimator randomness, $\Prob{j_k = j} = \p_j$, and $\Expsk{k}{e^k} = 0$, taking the squared norm of this identity and then the conditional expectation $\Expsk{k}{\cdot}$, we get
\myeqn{
	\arraycolsep=0.1em
	\begin{array}{lcl}
		\Expsk{k}{\norms{x^{k {\!} + {\!} 1} {\!\!} - \xopt {\!} + {\!} \gamma\eta \circ \hat{w}^{k{\!}+{\!}1}}^2 }  & = & \norms{x^k {\!} - {\!} \xopt {\!} + {\!} \gamma\eta \circ w^k}^2 - 2\eta \circ \iprods{w^k {\!} - {\!} \gamma \circ \hat{w}^k, x^k {\!} - {\!} \xopt {\!} + {\!} \gamma\eta \circ w^k} \vspace{1ex}\\
		&& + {~} \sum_{j=1}^m \frac{\eta_j^2}{\p_j}\big(\norms{[w^k]_j - \gamma_j[\hat{w}^k]_j}^2 + \Expsk{k}{\norms{[e^k]_j}^2}\big).
	\end{array}
}
We expand the first two terms on the right-hand side and subtract the result from $\norms{x^k - \xopt + \gamma\eta \circ \hat{w}^k}^2 = \norms{x^k - \xopt}^2 + 2\gamma\eta \circ \iprods{\hat{w}^k, x^k - \xopt} + \gamma^2\eta^2 \circ \norms{\hat{w}^k}^2$.
The terms $\norms{x^k - \xopt}^2$ and $2\gamma\eta \circ \iprods{\hat{w}^k, x^k - \xopt}$ cancel, and we obtain
\myeq{eq:VCOG_keylemma_proof2}{
	\arraycolsep=0.2em
	\begin{array}{lcl}
		\multicolumn{3}{l}{\norms{x^k - \xopt + \gamma\eta \circ \hat{w}^k}^2 - \Expsk{k}{\norms{x^{k+1} - \xopt + \gamma\eta \circ \hat{w}^{k+1}}^2}} \vspace{1ex}\\
		\qquad &=& 2(1-\gamma)\eta \circ \iprods{w^k, x^k - \xopt} + \gamma(2 - \gamma)\eta^2 \circ \norms{w^k}^2 + \gamma^2\eta^2 \circ \norms{\hat{w}^k}^2 \vspace{1ex}\\
		&& - {~} 2\gamma^2\eta^2 \circ \iprods{w^k, \hat{w}^k} - \sum_{j=1}^m \frac{\eta_j^2}{\p_j}\big(\norms{[w^k]_j - \gamma_j[\hat{w}^k]_j}^2 + \Expsk{k}{\norms{[e^k]_j}^2}\big).
	\end{array}
}
Second, by \eqref{eq:VCOG_Lipschitz_resolvent}, we have $\norms{x^{k+1} - x^k} = \norms{J_{\gamma\eta \circ T}(y^{k+1}) - J_{\gamma\eta \circ T}(y^k)} \leq \mu\norms{y^{k+1} - y^k}$.
Taking the conditional expectation of the square of \eqref{eq:VCOG_keylemma_proof1} as above, we get
\myeqn{
	\arraycolsep=0.2em
	\begin{array}{lcl}
		\frac{1}{\mu^2}\Expsk{k}{\norms{x^{k+1} - x^k}^2} &\leq& \sum_{j=1}^m \frac{\eta_j^2}{\p_j}\big(\norms{[w^k]_j - \gamma_j[\hat{w}^k]_j}^2 + \Expsk{k}{\norms{[e^k]_j}^2}\big).
	\end{array}
}
Adding $\frac{1}{\mu^2}\norms{x^k - x^{k-1}}^2$ to both sides of \eqref{eq:VCOG_keylemma_proof2}, subtracting this inequality, and using $\Lc_k$ in \eqref{eq:VCOG_potential_func}, we obtain
\myeq{eq:VCOG_keylemma_proof3}{
	\arraycolsep=0.2em
	\begin{array}{lcl}
		\Lc_k - \Expsk{k}{\Lc_{k+1}} &\geq& \frac{1}{\mu^2}\norms{x^k - x^{k-1}}^2 + 2(1-\gamma)\eta \circ \iprods{w^k, x^k - \xopt} \vspace{1ex}\\
		&& + {~} (2 - \gamma)\gamma\eta^2 \circ \norms{w^k}^2 + \gamma^2\eta^2 \circ \norms{\hat{w}^k}^2 - 2\gamma^2\eta^2 \circ \iprods{w^k, \hat{w}^k} \vspace{1ex}\\
		&& - {~} \sum_{j=1}^m \frac{2\eta_j^2}{\p_j}\norms{[w^k]_j - \gamma_j[\hat{w}^k]_j}^2 - \sum_{j=1}^m \frac{2\eta_j^2}{\p_j}\Expsk{k}{\norms{[e^k]_j}^2}.
	\end{array}
}
Third, since $[\hat{w}^k]_j = [w^k]_j - [Gx^k - Gx^{k-1}]_j$, a direct expansion gives
\myeqn{
	\arraycolsep=0.2em
	\begin{array}{lcl}
		\multicolumn{3}{l}{(2-\gamma_j)\gamma_j\norms{[w^k]_j}^2 + \gamma_j^2\norms{[\hat{w}^k]_j}^2 - 2\gamma_j^2\iprods{[w^k]_j, [\hat{w}^k]_j}} \vspace{1ex}\\
		\qquad &=& 2\gamma_j(1-\gamma_j)\norms{[w^k]_j}^2 + \gamma_j^2\norms{[Gx^k - Gx^{k-1}]_j}^2, \vspace{1ex}\\
		\norms{[w^k]_j - \gamma_j[\hat{w}^k]_j}^2 &=& \norms{(1-\gamma_j)[w^k]_j + \gamma_j[Gx^k - Gx^{k-1}]_j}^2.
	\end{array}
}
Substituting these identities into \eqref{eq:VCOG_keylemma_proof3} and using $(1-\gamma_j)\eta_j = \tau$ ($i \in [m]$), we obtain
\myeq{eq:VCOG_keylemma_proof4}{
	\arraycolsep=0.2em
	\begin{array}{lcl}
		\Lc_k - \Expsk{k}{\Lc_{k+1}} &\geq& \frac{1}{\mu^2}\norms{x^k - x^{k-1}}^2 + 2\tau\iprods{w^k, x^k - \xopt} - \sum_{j=1}^m \frac{2\eta_j^2}{\p_j}\Expsk{k}{\norms{[e^k]_j}^2} \vspace{1ex}\\
		&& + {~} \sum_{j=1}^m \frac{2(1-\gamma_j)\eta_j^2[\gamma_j\p_j - (1-\gamma_j)]}{\p_j}\norms{[w^k]_j}^2 \vspace{1ex}\\
		&& - {~} \sum_{j=1}^m \frac{4(1-\gamma_j)\gamma_j\eta_j^2}{\p_j}\iprods{[w^k]_j, [Gx^k - Gx^{k-1}]_j} \vspace{1ex}\\
		&& - {~} \sum_{j=1}^m \frac{(2-\p_j)\gamma_j^2\eta_j^2}{\p_j}\norms{[Gx^k - Gx^{k-1}]_j}^2.
	\end{array}
}
Finally, \eqref{eq:VCOG_Lipschitz_continuity} and Jensen's inequality give $\frac{1}{2\mu^2}\norms{x^k - x^{k-1}}^2 \geq \sum_{j=1}^m \frac{1}{2\mu^2L_j^2}\norms{[Gx^k - Gx^{k-1}]_j}^2$.
Hence, half of the first term in \eqref{eq:VCOG_keylemma_proof4} and its last two lines are bounded from below by $\sum_{j=1}^m\big(\phi_j\norms{[Gx^k - Gx^{k-1}]_j}^2 - \frac{4(1-\gamma_j)\gamma_j\eta_j^2}{\p_j}\iprods{[w^k]_j, [Gx^k - Gx^{k-1}]_j}\big)$.
Since $\phi_j > 0$ by \eqref{eq:VCOG_para_cond1}, Young's inequality gives
\myeqn{
	\arraycolsep=0.2em
	\begin{array}{lcl}
		\phi_j\norms{[Gx^k - Gx^{k-1}]_j}^2 - \frac{4(1-\gamma_j)\gamma_j\eta_j^2}{\p_j}\iprods{[w^k]_j, [Gx^k - Gx^{k-1}]_j} &\geq& -\frac{4(1-\gamma_j)^2\gamma_j^2\eta_j^4}{\p_j^2\phi_j}\norms{[w^k]_j}^2.
	\end{array}
}
Substituting this bound and $\iprods{w^k, x^k - \xopt} \geq -\rho \circ \norms{w^k}^2$ from Assumption~\ref{as:VCOG}(i) into \eqref{eq:VCOG_keylemma_proof4}, we obtain exactly \eqref{eq:VCOG_key_lemma1}.
\Eproof
	
	\vspace{2ex}
	\noindent\textbf{\textit{Proof of Lemma~\ref{lem:ch2:VCOG_key_lemma2}.}}
	Since $\beta_j = L_j^{-2}$, the second condition in \eqref{eq:VR_estimator} gives
	\myeq{eq:VCOG_key_lemma2_proof1}{
		\arraycolsep=0.2em
		\begin{array}{lcl}
			\sum_{j=1}^m \frac{\eta_j^2}{\p_j}\Expsk{k}{\norms{[e^k]_j}^2} &\leq& M_1\Expsk{k}{\norms{e^k}_{\beta}^2} \leq M_1\Expsk{k}{\Delta_k}.
		\end{array}
	}
	Moreover, rearranging \eqref{eq:VR_estimator2}, we can show that
	\myeqn{
		\arraycolsep=0.2em
		\begin{array}{lcl}
			\Expsk{k}{\Delta_k} &\leq& \frac{1-\kappa}{\kappa}\big(\Delta_{k-1} - \Expsk{k}{\Delta_k}\big) + \frac{\hat{\Theta}}{\kappa}\big(U_{k-1} - U_k\big) + \frac{\Theta + \hat{\Theta}}{\kappa}U_k + \frac{\norms{1-\gamma}_{\infty}^2\delta_k}{\kappa}V_{k-1}.
		\end{array}
	}
	We substitute these two bounds into \eqref{eq:VCOG_key_lemma1}, use $U_k \leq \norms{x^k - x^{k-1}}^2$ from \eqref{eq:VCOG_Lipschitz_continuity}, and recall the definition of $\Ec_k$ in \eqref{eq:VCOG_potential_func}.
	This leads to
	\myeq{eq:VCOG_key_lemma2_proof3}{
		\arraycolsep=0.2em
		\begin{array}{lcl}
			\Ec_k - \Expsk{k}{\Ec_{k+1}} &\geq& \Big[\frac{1}{2\mu^2} - \frac{2M_1(\Theta + \hat{\Theta})}{\kappa}\Big]\norms{x^k - x^{k-1}}^2 + \sum_{j=1}^m \psi_j\norms{[w^k]_j}^2 \vspace{1ex}\\
			&& - {~} \frac{2M_1\norms{1-\gamma}_{\infty}^2}{\kappa}\delta_kV_{k-1}.
		\end{array}
	}
	It remains to bound $V_{k-1}$.
	Let $M_2 := \max_{j \in [m]}\gamma_j\eta_jL_j$.
	Since $\xi^k = \hat{w}^k - Gx^{k-1} \in Tx^k$ and $-G\xopt \in T\xopt$, we have $x^k = J_{\gamma\eta\circ T}\big(x^k + \gamma\eta \circ (\hat{w}^k - Gx^{k-1})\big)$ and $\xopt = J_{\gamma\eta\circ T}(\xopt - \gamma\eta \circ G\xopt)$.
	Therefore, \eqref{eq:VCOG_Lipschitz_resolvent} and \eqref{eq:VCOG_Lipschitz_continuity} yield
	\myeq{eq:VCOG_key_lemma2_proof2}{
		\arraycolsep=0.2em
		\begin{array}{lcl}
			\norms{x^k - \xopt} &\leq& \mu\norms{x^k - \xopt + \gamma\eta \circ \hat{w}^k} + \mu\norms{\gamma\eta \circ (Gx^{k-1} - G\xopt)} \vspace{1ex}\\
			&\leq& \mu\norms{x^k - \xopt + \gamma\eta \circ \hat{w}^k} + \mu M_2\norms{x^{k-1} - \xopt}.
		\end{array}
	}
	Combining this inequality with $\norms{x^{k-1} - \xopt} \leq \norms{x^k - x^{k-1}} + \norms{x^k - \xopt}$, we get
	\myeqn{
		\arraycolsep=0.2em
		\begin{array}{lcl}
			(1 - \mu M_2)\norms{x^{k-1} - \xopt} &\leq& \mu\norms{x^k - \xopt + \gamma\eta \circ \hat{w}^k} + \norms{x^k - x^{k-1}}.
		\end{array}
	}
	By \eqref{eq:VCOG_para_cond1}, $\mu^2L_j^2\gamma_j^2\eta_j^2 < \frac{\p_j}{2(2-\p_j)} \leq \frac{1}{2}$, so $\mu M_2 < \frac{1}{\sqrt{2}}$ and $(1 - \mu M_2)^{-2} < 16$.
	Hence,
	\myeq{eq:VCOG_key_lemma2_proof3.1}{
		\arraycolsep=0.2em
		\begin{array}{lcl}
			\norms{x^{k-1} - \xopt}^2 &\leq& \frac{2\mu^2}{(1 - \mu M_2)^2}\Big[\norms{x^k - \xopt + \gamma\eta \circ \hat{w}^k}^2 + \frac{1}{\mu^2}\norms{x^k - x^{k-1}}^2\Big] \leq 32\mu^2\Lc_k.
		\end{array}
	}
	Since $\Lc_k \leq \Ec_k$, \eqref{eq:VCOG_Vk_bound} and \eqref{eq:VCOG_key_lemma2_proof3.1} imply
	\myeq{eq:VCOG_key_lemma2_proof4}{
		\arraycolsep=0.2em
		\begin{array}{lcl}
			V_{k-1} &\leq& 2\norms{x^{k-1} - \xopt}^2 + 2\sigma_{\star}^2 \leq 64\mu^2\Ec_k + 2\sigma_{\star}^2.
		\end{array}
	}
	Substituting \eqref{eq:VCOG_key_lemma2_proof4} into \eqref{eq:VCOG_key_lemma2_proof3} and rearranging the result, we obtain \eqref{eq:VCOG_key_lemma2} with $\hat{\delta}_k = \frac{128\mu^2M_1\norms{1-\gamma}_{\infty}^2}{\kappa}\delta_k$.
	\Eproof
	
	\beforesec
	\section{The Proofs of Technical Results in Section~\ref{sec:ch2:VACOG}}\label{apdx:sec:ch2:VACOG}
	\aftersec
	
	\subsection{The proofs of technical lemmas in Section~\ref{sec:ch2:VACOG}}\label{apdx:subsec:ch2:VACOG_technical_lemmas}
	\aftersubsec
	
	\noindent\textbf{\textit{Proof of Lemma~\ref{lem:ch2:BFS_FNI_equivalence}.}}
	(i) For any $u \in \R^p$, let $x := J_{\lambda \circ T}u$ and $\xi := \lambda^{-1} \circ (u - x) \in Tx$.
	Then, $\Gc_{\lambda}u = Gx + \xi$, which proves the identity $\norms{\Gc_{\lambda}u} = \norms{Gx + \xi}$ and the inclusion $J_{\lambda\circ T}(\zer{\Gc_{\lambda}}) \subseteq \zer{G+T}$.
	Conversely, if $0 \in Gx + Tx$, let $u := x - \lambda \circ Gx$.
	Since $-Gx \in Tx$, we have $J_{\lambda\circ T}u = x$ and $\Gc_{\lambda}u = Gx + \lambda^{-1} \circ (u - x) = 0$, which proves the reverse inclusion.
	
	\noindent(ii) Let $x := J_{\lambda \circ T}u$, $y := J_{\lambda \circ T}v$, $a := A_{\lambda \circ T}u$, and $b := A_{\lambda \circ T}v$, so that $a \in Tx$ and $b \in Ty$.
	By the $\rho$-co-hypomonotonicity of $T$, for each block $j$, we have
	\myeqn{
		\arraycolsep=0.2em
		\begin{array}{lcl}
			\iprods{a_j - b_j, u_j - v_j} &=& \iprods{a_j - b_j, x_j - y_j} + \lambda_j\norms{a_j - b_j}^2 \geq (\lambda_j - \rho_j)\norms{a_j - b_j}^2.
		\end{array}
	}
	Hence, $A_{\lambda_j T_j}$ is $(\lambda_j - \rho_j)$-co-coercive.
	Moreover, since $u_j - v_j = x_j - y_j + \lambda_j(a_j - b_j)$, we get $\norms{u_j - v_j}^2 \geq \norms{x_j - y_j}^2 + \lambda_j(\lambda_j - 2\rho_j)\norms{a_j - b_j}^2$.
	If $\lambda_j \geq 2\rho_j$ for all $j \in [m]$, summing this inequality over the blocks shows that $J_{\lambda \circ T}$ is nonexpansive.
	
	\noindent(iii) Let $x := J_{\lambda \circ T}u$ and $y := J_{\lambda \circ T}v$, and denote $D_{ij} := [G_i x - G_i y]_j$ and $F_{ij} := [\Gc_{i, \lambda} u - \Gc_{i, \lambda} v]_j$.
	Then, $F_{ij} - D_{ij} = [A_{\lambda\circ T}u - A_{\lambda\circ T}v]_j$ and $u_j - v_j = x_j - y_j + \lambda_j(F_{ij} - D_{ij})$.
	Using the co-coercivity of $A_{\lambda_jT_j}$ from part (ii), we obtain
	\myeqn{
		\arraycolsep=0.2em
		\begin{array}{lcl}
			\iprods{F_{ij}, u_j - v_j} &\geq& (\lambda_j - \rho_j)\norms{F_{ij}}^2 - (\lambda_j - 2\rho_j)\iprods{F_{ij}, D_{ij}} + \iprods{D_{ij}, x_j - y_j} - \rho_j\norms{D_{ij}}^2.
		\end{array}
	}
	We sum this inequality over $j \in [m]$, average it over $i \in [n]$, and apply \eqref{eq:FNI_assumption_cocoercive}.
	With $c_{ij} := L_{ij}^{-1} - \rho_j > 0$, each resulting quadratic form satisfies
	\myeqn{
		\arraycolsep=0.2em
		\begin{array}{lcl}
			(\lambda_j - \rho_j)\norms{F_{ij}}^2 - (\lambda_j - 2\rho_j)\iprods{F_{ij}, D_{ij}} + c_{ij}\norms{D_{ij}}^2
			&=& \bar{\tau}_{ij}\norms{F_{ij}}^2 + c_{ij}\norms{D_{ij} - \frac{\lambda_j - 2\rho_j}{2c_{ij}}F_{ij}}^2 \vspace{1ex}\\
			&\geq& \bar{\tau}_{ij}\norms{F_{ij}}^2,
		\end{array}
	}
	which proves \eqref{eq:BFS_cocoercive}.
	The interval prescribed for $\lambda_j$ is exactly the set on which $4(\lambda_j - \rho_j) - \lambda_j^2\bar{L}_j > 0$, so $\bar{\tau}_{ij} > 0$ for all $i$ and $j$.
	Finally, \eqref{eq:BFS_cocoercive} and Jensen's inequality give $\underline{\bar{\tau}}\norms{\Gc_{\lambda}u - \Gc_{\lambda}v}^2 \leq \iprods{\Gc_{\lambda}u - \Gc_{\lambda}v, u - v}$, where $\underline{\bar{\tau}} := \min_{i,j}\bar{\tau}_{ij}$, and the Cauchy-Schwarz inequality shows that $\Gc_{\lambda}$ is $\underline{\bar{\tau}}^{-1}$-Lipschitz continuous.
	\Eproof
	
\vspace{0.75ex}
\noindent\textbf{\textit{Proof of Lemma~\ref{le:VACOG_key_est1}.}}
We split the proof of this lemma into the following steps.

\noindent\textit{\underline{Step 1.} Expanding $\Qc_k - \Qc_{k+1}$.}
By \eqref{eq:VACOG}, we have $r(u^k - \uopt) + t_{k+1}(u^{k+1} - u^k) = r(u^{k-1} - \uopt) + (r + t_{k+1}\theta_k)(u^k - u^{k-1}) - t_{k+1}\eta_k\p_{j_k}^{-1}\widetilde{\Sc}_{\lambda, [j_k]}^k$.
Expanding the squares in $\Qc_k$ and $\Qc_{k+1}$ defined by \eqref{eq:VACOG_Lyapunov_funcs} as in \cite{tran2022accelerated}, we get
\myeq{eq:lem1_proof1}{
\arraycolsep=0.2em
\hspace{-4ex}
\begin{array}{lcl}
	\Qc_k - \Qc_{k+1} &=& [(t_k - r)^2 - t_{k+1}^2 \theta_k^2 + r]\norms{u^k - u^{k-1}}^2 - t_{k+1}^2 \eta_k^2 \norms{\p_{j_k}^{-1} \widetilde{\Sc}_{\lambda, [j_k]}^k}^2 \vspace{1ex}\\
	&& + {~} 2r t_{k+1} \eta_k \iprods{\p_{j_k}^{-1}\widetilde{\Sc}_{\lambda, [j_k]}^k, u^k - \uopt} + 2t_{k+1}^2 \theta_k \eta_k \iprods{\p_{j_k}^{-1}\widetilde{\Sc}_{\lambda, [j_k]}^k, u^k - u^{k-1}} \vspace{1ex}\\
	&& + {~} 2r (t_k - r - 1 - t_{k+1} \theta_k)\iprods{u^k - u^{k-1}, u^k - \uopt}.
\end{array}
\hspace{-3ex}
}
\noindent\textit{\underline{Step 2.} Taking conditional expectations.}
Let $e^k := \widetilde{\Sc}_{\lambda}^k - \Sc_{\lambda}^k$.
Since $j_k$ is independent of the estimator randomness, $\Expsk{k}{e^k} = 0$, and \eqref{eq:VACOG_prob} holds, for every $\Fc_k$-measurable random vector $v$, we have
\myeq{eq:VACOG_conditional_expectations}{
	\arraycolsep=0.2em
	\begin{array}{lcl}
		\Expsk{k}{\norms{\p_{j_k}^{-1} \widetilde{\Sc}_{\lambda, [j_k]}^k}^2} = \norms{\Sc_{\lambda}^k}_{\beta}^2 + \Expsk{k}{\norms{e^k}_{\beta}^2} \quad \text{and} \quad
		\Expsk{k}{\iprods{\p_{j_k}^{-1}\widetilde{\Sc}_{\lambda, [j_k]}^k, v}} = \iprods{\Sc_{\lambda}^k, v}.
	\end{array}
}
By the update rules of $\theta_k$ and $\gamma_k$, we have $t_k - r - 1 - t_{k+1}\theta_k = 0$ and $t_{k+1}\theta_k(1 - \gamma_k) - r\gamma_k = 0$.
We take the conditional expectation of \eqref{eq:lem1_proof1}, use \eqref{eq:VACOG_conditional_expectations} with $\Sc_{\lambda}^k = \Gc_{\lambda}u^k - \gamma_k\Gc_{\lambda}u^{k-1}$, and write $u^k - \uopt = (u^{k-1} - \uopt) + (u^k - u^{k-1})$ and $\Sc_{\lambda}^k = (\Gc_{\lambda}u^k - \Gc_{\lambda}u^{k-1}) + (1 - \gamma_k)\Gc_{\lambda}u^{k-1}$ in the two inner products.
Then, the coefficient of $\iprods{\Gc_{\lambda}u^{k-1}, u^k - u^{k-1}}$ is $2t_{k+1}\eta_k[t_{k+1}\theta_k(1 - \gamma_k) - r\gamma_k] = 0$, and the last line of \eqref{eq:lem1_proof1} vanishes, which leads to
\myeq{eq:lem1_proof2}{
\arraycolsep=0.2em
\hspace{-3ex}
\begin{array}{lcl}
	\Qc_k - \Expsk{k}{\Qc_{k+1}} &=& [(t_k - r)^2 - t_{k+1}^2 \theta_k^2 + r]\norms{u^k - u^{k-1}}^2 - t_{k+1}^2 \eta_k^2 \Expsk{k}{\norms{e^k}_{\beta}^2} \vspace{1ex}\\
	&& + {~} 2t_{k+1}^2 \theta_k \eta_k \iprods{\Gc_{\lambda}u^k - \Gc_{\lambda}u^{k-1}, u^k - u^{k-1}} \vspace{1ex}\\
	&& - {~} t_{k+1}^2 \eta_k^2 \sum_{j=1}^m \frac{1}{\p_j} \norms{[\Gc_{\lambda}u^k]_j - \gamma_k [\Gc_{\lambda}u^{k-1}]_j}^2 \vspace{1ex}\\
	&& + {~} 2r t_{k+1} \eta_k \iprods{\Gc_{\lambda}u^k, u^k - \uopt} - 2r t_{k+1} \eta_k \gamma_k \iprods{\Gc_{\lambda}u^{k-1}, u^{k-1} - \uopt}.
\end{array}
\hspace{-3ex}
}
\noindent\textit{\underline{Step 3.} Lower bounding $\Qc_k - \Expsk{k}{\Qc_{k+1}}$.}
By \eqref{eq:BFS_cocoercive} and Jensen's inequality, we have
\myeq{eq:lem1_proof3}{
	\arraycolsep=0.2em
	\begin{array}{lcl}
		\iprods{\Gc_{\lambda}u^k - \Gc_{\lambda}u^{k-1}, u^k - u^{k-1}} &\geq& \frac{1}{n} \sum\limits_{i = 1}^n \sum\limits_{j = 1}^m (\bar{\tau}_{ij} - \tau_j) \norms{[\Gc_{i, \lambda}u^k]_j - [\Gc_{i, \lambda} u^{k-1}]_j}^2 \vspace{0ex}\\
		&& + {~} \sum\limits_{j = 1}^m \tau_j \norms{[\Gc_{\lambda}u^k]_j - [\Gc_{\lambda} u^{k-1}]_j}^2.
	\end{array}
}
We substitute \eqref{eq:lem1_proof3} into \eqref{eq:lem1_proof2}, and add and subtract $2rt_{k+1}\eta_k\sum_{j=1}^m\tau_j\norms{[\Gc_{\lambda}u^k]_j}^2$ and $2rt_{k+1}\eta_k\gamma_k\sum_{j=1}^m\tau_j\norms{[\Gc_{\lambda}u^{k-1}]_j}^2$.
Since $\gamma_k = \frac{t_{k+1}\theta_k}{t_{k+1}\theta_k + r}$, a direct expansion shows that
\myeqn{
	\arraycolsep=0.2em
	\begin{array}{lcl}
		\multicolumn{3}{l}{t_{k+1}\theta_k\norms{[\Gc_{\lambda}u^k]_j - [\Gc_{\lambda}u^{k-1}]_j}^2 + r\norms{[\Gc_{\lambda}u^k]_j}^2 - r\gamma_k\norms{[\Gc_{\lambda}u^{k-1}]_j}^2} \vspace{1ex}\\
		\qquad &=& (t_{k+1}\theta_k + r)\norms{[\Gc_{\lambda}u^k]_j - \gamma_k[\Gc_{\lambda}u^{k-1}]_j}^2.
	\end{array}
}
Hence, the terms with $\tau_j$ combine with the third line of \eqref{eq:lem1_proof2}, and the total coefficient of $\norms{[\Gc_{\lambda}u^k]_j - \gamma_k[\Gc_{\lambda}u^{k-1}]_j}^2$ becomes $t_{k+1}\eta_k\frac{2\tau_j\p_j(t_{k+1}\theta_k + r) - t_{k+1}\eta_k}{\p_j}$.
Finally, by the definition of $\Lc_k$ in \eqref{eq:VACOG_Lyapunov_funcs}, the bracket $2rt_{k+1}\eta_k\big[\iprods{\Gc_{\lambda}u^k, u^k - \uopt} - \sum_{j=1}^m\tau_j\norms{[\Gc_{\lambda}u^k]_j}^2\big]$ cancels with the correction term in $\Lc_{k+1}$.
The bracket at $u^{k-1}$ is multiplied by $2r(t_k\eta_{k-1} - t_{k+1}\eta_k\gamma_k) = 2r\phi_k$, and we obtain \eqref{eq:VACOG_key_estimate1}.
\Eproof
	
	\vspace{0.5ex}
	\noindent\textbf{\textit{Proof of Lemma~\ref{lem:ch2:VACOG_lemma2}.}}
	By \eqref{eq:lem2_param_conditions}, we have $1 - \gamma_k = \frac{r}{t_k - 1}$ and $\nu_k = \frac{t_{k-1} - 1}{t_k - 1}$.
	Hence, $(t_k - 1)^2\nu_k^2 = (t_{k-1} - 1)^2$, $(t_k - 1)^2(1 - \gamma_k)^2 = r^2$, and $t_{k+1}^2\eta_k^2 = \eta^2(t_k - 1)^2$.
	Multiplying \eqref{eq:VR_estimator3} by $(t_k - 1)^2$ and rearranging, we can write it as
	\myeq{eq:lemma2_proof1}{
		\arraycolsep=0.2em
		\begin{array}{lcl}
			(t_k - 1)^2\Expsk{k}{\Delta_k} &\leq& \frac{1-\kappa}{\kappa}\Big[(t_{k-1} - 1)^2\Delta_{k-1} - (t_k - 1)^2\Expsk{k}{\Delta_k}\Big] \vspace{1ex}\\
			&& + {~} \frac{\hat{\Theta}}{\kappa}\Big[(t_{k-1} - 1)^2U_{k-1} - (t_k - 1)^2U_k\Big] \vspace{1ex}\\
			&& + {~} \frac{\Theta + \hat{\Theta}}{\kappa}(t_k - 1)^2U_k + \frac{r^2}{\kappa}\delta_kV_{k-1}.
		\end{array}
	}
	We multiply \eqref{eq:lemma2_proof1} by $\eta^2$ and use $\Expsk{k}{\norms{e^k}_{\beta}^2} \leq \Expsk{k}{\Delta_k}$ in \eqref{eq:VACOG_key_estimate1}.
	The first two differences on the right-hand side of \eqref{eq:lemma2_proof1} are the correction terms in $\Pc_k$.
	Moreover, \eqref{eq:lem2_param_conditions} gives $(t_k - r)^2 - t_{k+1}^2\theta_k^2 + r = 2t_k - r - 1$, $t_{k+1}^2\theta_k\eta_k = \eta(t_k - 1)(t_k - r - 1)$, and $\phi_k = \eta(t_{k-1} - t_k + r)$.
	Together with the definition of $U_k$ (with $\beta_j = \p_j^{-1}$), these identities yield all descent terms in \eqref{eq:VACOG_key_estimate2}, up to the remainder $\frac{r^2\eta^2}{\kappa}\delta_kV_{k-1}$.
	Finally, since $\Pc_k \geq r\norms{u^{k-1} - \uopt}^2$, we can derive from \eqref{eq:VACOG_Vk_bound} that
	\myeq{eq:lemma2_proof3}{
		\arraycolsep=0.2em
		\begin{array}{lcl}
			V_{k-1} &\leq& \frac{2\bar{\p}}{\Lambda_0^2}\norms{u^{k-1} - \uopt}^2 + 2\sigma_{\lambda,\star}^2 \leq \frac{2\bar{\p}}{r\Lambda_0^2}\Pc_k + 2\sigma_{\lambda,\star}^2.
		\end{array}
	}
	Substituting \eqref{eq:lemma2_proof3} into the remainder, we obtain \eqref{eq:VACOG_key_estimate2}.
	\Eproof
	
	\beforesubsec
	\subsection{Supporting lemmas used in Theorems~\ref{th:VACOG_convergence} and \ref{th:VACOG_convergence_as}}
	\aftersubsec
	\begin{lemma}\label{lem:ch2:VACOG_P0_bound}
		Under the hypotheses of Theorem~\ref{th:VACOG_convergence}, $\Pc_k$ defined by \eqref{eq:VACOG_Lyapunov_funcs} satisfies
		\myeq{eq:th1_boundP0}{
			\arraycolsep=0.2em
			\begin{array}{lcl}
				\Pc_0 &\leq& \left[ \frac{r^2 \eta}{2\underline{\tau}} + r^2 + r \right]\norms{u^0 - \uopt}^2 \leq C_0 \norms{u^0 - \uopt}^2,
			\end{array}
		}
		where $C_0 := r^2\left(\frac{\eta}{2\underline{\tau}} + 2\right)$ as defined in \eqref{eq:VACOG_constants2}.
	\end{lemma}
	
	\begin{proof}
		By \eqref{eq:BFS_cocoercive}, Jensen's inequality, and $\Gc_{\lambda}\uopt = 0$, we have $\underline{\bar{\tau}}\norms{\Gc_{\lambda}u^0}^2 \leq \iprods{\Gc_{\lambda}u^0, u^0 - \uopt}$.
		The Cauchy--Schwarz inequality then gives
		\myeq{eq:th1_bound_Gu0}{
			\arraycolsep=0.2em
			\begin{array}{lcl}
				\norms{\Gc_{\lambda}u^0}^2 &\leq& \frac{1}{\underline{\bar{\tau}}^2}\norms{u^0 - \uopt}^2.
			\end{array}
		}
		Next, since $u^{-1} = u^0$, $\Delta_{-1} = U_{-1} = 0$, and $t_0\eta_{-1} = r\eta$, the definition \eqref{eq:VACOG_Lyapunov_funcs} gives
		\myeqn{
			\arraycolsep=0.2em
			\begin{array}{lcl}
				\Pc_0 &=& (r^2 + r)\norms{u^0 - \uopt}^2 + 2r^2\eta\sum_{j=1}^m\Big[\iprods{[\Gc_{\lambda}u^0]_j, u^0_j - \uopt_j} - \tau_j\norms{[\Gc_{\lambda}u^0]_j}^2\Big].
			\end{array}
		}
		By Young's inequality, $\iprods{[\Gc_{\lambda}u^0]_j, u^0_j - \uopt_j} \leq \tau_j\norms{[\Gc_{\lambda}u^0]_j}^2 + \frac{1}{4\underline{\tau}}\norms{u^0_j - \uopt_j}^2$.
		Therefore, $\Pc_0 \leq \big[r^2 + r + \frac{r^2\eta}{2\underline{\tau}}\big]\norms{u^0 - \uopt}^2 \leq C_0\norms{u^0 - \uopt}^2$, where the last inequality uses $r \leq r^2$.
	\end{proof}
	
	\begin{lemma}\label{lem:ch2:VACOG_Wk_bound}
		Under the hypotheses of Theorem~\ref{th:VACOG_convergence}, let
		\myeq{eq:VACOG_Wk}{
			\arraycolsep=0.2em
			\hspace{-1ex}
			\begin{array}{lcl}
				\textls[-10]{$W_k := (k+r+3)(k+r) \norms{u^{k+1} - u^k + \eta_k \Gc_{\lambda}u^k}^2 + \frac{\eta^2(k+r+1)^2}{\kappa} [(1-\kappa)\Delta_k + \hat{\Theta}U_k],$}
			\end{array}
			\hspace{-1ex}
		}
		where $\Delta_k$ and $U_k$ are defined in \eqref{eq:VR_estimator}.
		Then, we have
		\myeq{eq:VACOG_Wk_bound}{
			\arraycolsep=0.2em
			\hspace{-8ex}
			\begin{array}{lcl}
				\Expsk{k}{W_k} &\leq& \textls[-70]{$W_{k-1} - (r-2)(k + r + 2){\norms{u^k - u^{k-1} + \eta_{k-1} \Gc_{\lambda} u^{k-1}}^2} + \frac{4(2k + r + 3)}{r} {\norms{u^k - u^{k-1}}^2}$} \vspace{1ex}\\
				&& + {~} \frac{\eta^2 (\Theta + \hat{\Theta}) (k + r + 1)^2}{\kappa} U_{k} + (1-\underline{\p})\eta^2(k + r + 1)^2 \norms{\Sc_{\lambda}^k}_{\beta}^2 + \frac{r^2\eta^2}{\kappa}\delta_k V_{k-1}.
			\end{array}
			\hspace{-10ex}
		}
	\end{lemma}
	
	\begin{proof}
		First, multiplying \eqref{eq:VR_estimator3} by $(k+r+1)^2$ and using $\nu_k = \frac{k+r}{k+r+1}$ and $(1-\gamma_k)^2(k+r+1)^2 = r^2$, we obtain
		\myeq{eq:th1_proof2}{
			\arraycolsep=0.2em
			\hspace{-4ex}
			\begin{array}{lcl}
				(k + r + 1)^2 \Expsk{k}{\Delta_{k}} &\leq& \frac{1 - \kappa}{\kappa}\left[(k + r)^2 \Delta_{k-1} - (k + r + 1)^2 \Expsk{k}{\Delta_{k}} \right] + \frac{r^2 \delta_k}{\kappa}V_{k-1} \vspace{1ex}\\
				&& + {~} \frac{\hat{\Theta}}{\kappa} \left[ (k + r)^2 U_{k-1} - (k + r + 1)^2 U_{k} \right] + \frac{(\Theta + \hat{\Theta}) (k + r + 1)^2}{\kappa} U_{k}.
			\end{array}
			\hspace{-5ex}
		}
		Next, let us define $\bar{u}^{k+1} := u^k + \theta_k (u^k - u^{k-1}) - \eta_k (\Gc_{\lambda} u^k - \gamma_k \Gc_{\lambda} u^{k-1})$.
		Denoting $\omega_k := \frac{\eta_k \gamma_k}{\eta_{k-1}}$, then we can derive that
		\myeqn{
			\arraycolsep=0.2em
			\begin{array}{lcl}
				\bar{u}^{k+1} - u^k + \eta_k \Gc_{\lambda} u^k &=& \omega_k (u^k - u^{k-1} + \eta_{k-1} \Gc_{\lambda} u^{k-1}) + (1 - \omega_k) \frac{\omega_k - \theta_k}{1 - \omega_k} (u^{k-1} - u^k).
			\end{array}
		}
		Using \eqref{eq:th1_param} and $r > 2$, we obtain
		\myeq{eq:th1_tau_k}{
			\arraycolsep=0.2em
			\begin{array}{lcl}
				\omega_k &=& \frac{(k+1)(k+r+2)}{(k+r+3)(k+r)} \in (0,1)
				\quad \text{and} \quad
				\frac{\omega_k - \theta_k}{1 - \omega_k} = \frac{2(k+1)}{rk + r^2 + 2r - 2} \leq \frac{2}{r}.
			\end{array}
		}
		Hence, by the convexity of $\norms{\cdot}^2$, we have
		\myeq{eq:th1_proof3}{
			\arraycolsep=0.2em
			\hspace{-2ex}
			\begin{array}{lcl}
				\norms{\bar{u}^{k+1} {\!\!} - u^k + \eta_k \Gc_{\lambda} u^k}^2 &\leq& \omega_k \norms{u^k - u^{k-1} {\!\!} + \eta_{k-1} \Gc_{\lambda} u^{k-1}}^2 + \frac{4(1 - \omega_k)}{r^2} \norms{u^k - u^{k-1}}^2.
			\end{array}
			\hspace{-2ex}
		}
		From \eqref{eq:VACOG}, we have $u^{k+1} = u^k + \theta_k (u^k - u^{k-1}) - \frac{\eta_k}{\p_{j_k}} \widetilde{\Sc}^k_{\lambda, [j_k]}$.
		Using the conditional expectation identity \eqref{eq:VACOG_conditional_expectations}, the second condition of \eqref{eq:VR_estimator} in Definition~\ref{def:VR_estimator}, and  $\sum_{j = 1}^m \frac{1 - \p_j}{\p_j} \norms{\Sc_{\lambda, [j]}^k }^2 = \sum_{j = 1}^m \frac{1 - \p_j}{\p_j} \norms{[\Sc_{\lambda}^k]_j }^2 \leq (1-\underline{\p})\norms{\Sc_{\lambda}^k}_{\beta}^2$, for every $\Fc_k$-measurable $v^k$, we obtain
		\myeqn{
			\arraycolsep=0.2em
			\begin{array}{lcl}
				\Expsk{k}{\norms{u^{k+1} - v^k}^2}
				&\leq& \norms{\bar{u}^{k+1} - v^k}^2 + (1-\underline{\p}) \eta_k^2\norms{\Sc_{\lambda}^k}_{\beta}^2 + \eta_k^2 \Expsk{k}{\Delta_k}.
			\end{array}
		}
		Taking $v^k := u^k - \eta_k \Gc_{\lambda} u^k$, then using \eqref{eq:th1_proof3}, \eqref{eq:th1_param}, and \eqref{eq:th1_tau_k}, we obtain
		\myeqn{
		\arraycolsep=0.2em
		\begin{array}{lcl}
			\breve{\Tc}_{[1]} &:=& \Expsk{k}{\norms{u^{k+1} - u^k + \eta_k \Gc_{\lambda} u^k}^2} \vspace{1ex}\\
			&\leq& \norms{\bar{u}^{k+1} - u^k + \eta_k \Gc_{\lambda} u^k}^2 + (1-\underline{\p}) \eta_k^2 \norms{\Sc_{\lambda}^k}_{\beta}^2 + \eta_k^2 \Expsk{k}{\Delta_k} \vspace{1ex}\\
			&\leq& \omega_k \norms{u^k - u^{k-1} + \eta_{k-1} \Gc_{\lambda} u^{k-1}}^2 + \frac{4(1 - \omega_k)}{r^2} \norms{u^k - u^{k-1}}^2 \vspace{1ex}\\
			&& + {~} (1-\underline{\p}) \eta_k^2 \norms{\Sc_{\lambda}^k}_{\beta}^2 + \eta_k^2 \Expsk{k}{\Delta_k} \vspace{1ex}\\
			&\leq& \frac{(k+1)(k+r+2)}{(k+r+3)(k+r)} \norms{u^k - u^{k-1} + \eta_{k-1} \Gc_{\lambda} u^{k-1}}^2 + \frac{4(rk + r^2 + 2r - 2)}{r^2(k+r+3)(k+r)} \norms{u^k - u^{k-1}}^2 \vspace{1ex}\\
			&& + {~} \frac{\eta^2 (k + r + 1)^2}{(k+r+3)(k+r)} \Expsk{k}{\Delta_k} + \frac{(1-\underline{\p}) \eta^2 (k + r + 1)^2}{(k+r+3)(k+r)} \norms{\Sc_{\lambda}^k}_{\beta}^2,
		\end{array}
		}
		where we have used the fact that $k+r+3 \geq k+r$ in the last inequality.
		Multiplying both sides of the last inequality by $(k+r+3)(k+r)$, we obtain
		\myeq{eq:th1_proof4}{
		\arraycolsep=0.2em
		\begin{array}{lcl}
			\breve{\Tc}_{[2]} &:=& (k+r+3)(k+r)\Expsk{k}{\norms{u^{k+1} - u^k + \eta_k \Gc_{\lambda} u^k}^2} \vspace{1ex}\\
			&\leq& (k + r + 2)(k + r -1) \norms{u^k - u^{k-1} + \eta_{k-1} \Gc_{\lambda} u^{k-1}}^2 \vspace{1ex}\\
			&& - {~} (r-2)(k + r + 2) \norms{u^k - u^{k-1} + \eta_{k-1} \Gc_{\lambda} u^{k-1}}^2 \vspace{1ex}\\
			&& + {~} \frac{4}{r}\big(k + r + 2 - \frac{2}{r}\big) \norms{u^k - u^{k-1}}^2 + (1-\underline{\p}) \eta^2 (k + r + 1)^2 \norms{\Sc_{\lambda}^k}_{\beta}^2 \vspace{1ex}\\
			&& + {~} \eta^2 (k + r + 1)^2 \Expsk{k}{\Delta_k}.
		\end{array}
		}
		Substituting \eqref{eq:th1_proof2} into \eqref{eq:th1_proof4}, and then rearranging the result, we obtain \eqref{eq:VACOG_Wk_bound}.
	\end{proof}
	
	\begin{lemma}\label{lem:ch2:VACOG_Zk}
		Under the hypotheses of Theorem~\ref{th:VACOG_convergence}, let
		\myeq{eq:VACOG_Zk}{
			\arraycolsep=0.2em
			\begin{array}{lcl}
				Z_k &:=& (k+r+3)^2 \norms{u^{k+1} - u^k}^2 + \frac{\eta^2 (k + r + 1)^2 }{\kappa} [(1-\kappa)\Delta_{k} + \hat{\Theta} U_{k}],
			\end{array}
		}
		where $\Delta_k$ and $U_k$ are defined in \eqref{eq:VR_estimator}.
		Then, we have
		\myeq{eq:VACOG_Zk_bound}{
			\arraycolsep=0.2em
			\hspace{-4ex}
			\begin{array}{lcl}
				\Expsk{k}{Z_k} &\leq& Z_{k-1} - (r+2)(k+r) \norms{u^k - u^{k-1}}^2  \vspace{1ex}\\
				&& - {~} \frac{2\eta (k+1)^2}{n} {\sum\limits_{i=1}^n \sum\limits_{j=1}^m } \bar{\tau}_{ij} \norms{[\Gc_{i, \lambda} u^k]_j - [\Gc_{i, \lambda} u^{k-1}]_j}^2 + \frac{\eta^2(\Theta + \hat{\Theta}) (k + r + 1)^2}{\kappa} U_{k} \vspace{1ex}\\
				&& + {~} r \eta^2 (k+1) \norms{\Gc_{\lambda}u^k}^2 +  \eta^2 (k + r + 1)^2 \norms{\Gc_{\lambda}u^k - \gamma_k \Gc_{\lambda}u^{k-1}}_{\beta}^2 + \frac{r^2\eta^2}{\kappa}\delta_k V_{k-1}.
			\end{array}
			\hspace{-8ex}
		}
	\end{lemma}
	
	\begin{proof}
		From \eqref{eq:VACOG}, $u^{k+1} - u^k = \theta_k (u^k - u^{k-1}) - \eta_k \p_{j_k}^{-1} \widetilde{\Sc}_{\lambda, [j_k]}^k$, we have
		\myeqn{
			\arraycolsep=0.2em
			\begin{array}{lcl}
				\norms{u^{k+1} - u^k}^2 = \theta_k^2 \norms{u^k - u^{k-1}}^2 - 2\theta_k \eta_k \iprods{\p_{j_k}^{-1} \widetilde{\Sc}_{\lambda, [j_k]}^k, u^k - u^{k-1}} + \eta_k^2 \norms{\p_{j_k}^{-1} \widetilde{\Sc}_{\lambda, [j_k]}^k }^2.
			\end{array}
		}
		Taking the conditional expectation $\Expsk{k}{\cdot}$ of this expression, then using \eqref{eq:VACOG_conditional_expectations} and $\Expsk{k}{\norms{e^k}^2_{\beta}} \leq \Expsk{k}{\Delta_k}$ from the second condition of \eqref{eq:VR_estimator}, we obtain
		\myeqn{
			\arraycolsep=0.2em
			\begin{array}{lcl}
				\Expsk{k}{\norms{u^{k+1} - u^k}^2}
				&\leq& \textls[-15]{$\theta_k^2 \norms{u^k - u^{k-1}}^2 + \eta_k^2\norms{\Gc_{\lambda} u^k - \gamma_k \Gc_{\lambda} u^{k-1}}_{\beta}^2 + \eta_k^2\Expsk{k}{\Delta_k}$} \vspace{1ex}\\
				&& \textls[-85]{$- {~} 2\theta_k \eta_k \gamma_k \iprods{\Gc_{\lambda} u^k - \Gc_{\lambda} u^{k-1}, u^k - u^{k-1}} - 2\theta_k \eta_k (1 - \gamma_k) \iprods{\Gc_{\lambda} u^k, u^k - u^{k-1}}.$}
			\end{array}
		}
		Using \eqref{eq:BFS_cocoercive} and Young's inequality in the last expression, we can show that
		\myeqn{
				\arraycolsep=0.2em
				\begin{array}{lcl}
					\Expsk{k}{\norms{u^{k+1} - u^k}^2} &\leq& \left(\theta_k^2 + \frac{\theta_k \eta_k (1-\gamma_k)}{\eta}\right) \norms{u^k - u^{k-1}}^2 + \eta_k^2\norms{\Gc_{\lambda} u^k - \gamma_k \Gc_{\lambda} u^{k-1}}_{\beta}^2 + \eta_k^2 \Expsk{k}{\Delta_k}\vspace{0ex}\\
					&& - {~} \frac{2\theta_k \eta_k \gamma_k}{n} \sum\limits_{i=1}^n \sum\limits_{j=1}^m \bar{\tau}_{ij} \norms{[\Gc_{i, \lambda} u^k]_j - [\Gc_{i, \lambda} u^{k-1}]_j}^2 + \eta \theta_k \eta_k (1 - \gamma_k) \norms{\Gc_{\lambda} u^k}^2.
				\end{array}
		}
		By \eqref{eq:th1_param}, we have $\theta_k^2 + \frac{\theta_k \eta_k (1-\gamma_k)}{\eta} = \frac{(k+1)(k+r+1)}{(k+r+3)^2}$, $\theta_k \eta_k \gamma_k = \frac{\eta (k + 1)^2}{(k+r+3)^2}$, and $\theta_k \eta_k (1 - \gamma_k) = \frac{\eta r (k + 1)}{(k+r+3)^2}$.
		Substituting these expressions into the last estimate and multiplying both sides of the result by $(k+r+3)^2$, we can derive that
		\myeq{eq:th1_proof6}{
			\arraycolsep=0.2em
			\hspace{-2ex}
			\begin{array}{lcl}
				\breve{\Tc}_{[3]} & := & (k+r+3)^2\Expsk{k}{\norms{u^{k+1} - u^k}^2} \vspace{1ex}\\
				&\leq& (k+r+2)^2 {\norms{u^k - u^{k-1}}^2} - (r+2)(k+r) \norms{u^k - u^{k-1}}^2 \vspace{1ex}\\
				&& + {~} \eta^2 (k + r + 1)^2 \norms{\Gc_{\lambda} u^k - \gamma_k \Gc_{\lambda} u^{k-1}}_{\beta}^2 + \eta^2 (k + r + 1)^2 \Expsk{k}{\Delta_k} \vspace{1ex}\\
				&& + {~} r \eta^2 (k+1) {\norms{\Gc_{\lambda} u^k}^2} - \frac{2\eta}{n} (k+1)^2 \sum_{i=1}^n \sum_{j=1}^m \bar{\tau}_{ij} {\norms{[\Gc_{i, \lambda} u^k]_j - [\Gc_{i, \lambda} u^{k-1}]_j}^2}.
			\end{array}
			\hspace{-8ex}
		}
		Here, we have used the fact that $(k+1)(k+r+1) - (k+r+2)^2 = - (r+2)(k+r) - r - 3 \leq - (r+2)(k+r)$ for any $r \geq 1$.
		Substituting \eqref{eq:th1_proof2} into \eqref{eq:th1_proof6} and rearranging the result, we obtain \eqref{eq:VACOG_Zk_bound}.
	\end{proof}

	\renewcommand{\doi}[1]{doi:~\href{https://doi.org/#1}{\nolinkurl{#1}}}
	\bibliographystyle{plainnat}
	\bibliography{VCOG_refs}
	
\end{document}